\documentclass[12pt]{amsart}
\usepackage{amsfonts,amsthm,latexsym,amsmath,amssymb,amscd,amsmath,epsf,a4wide} 
\usepackage[margin=1in]{geometry}

\usepackage{tikz,tkz-graph}
\usetikzlibrary{calc,shapes,fit,positioning,arrows,automata,decorations.pathreplacing}
\usepackage{hyperref}
\usepackage{subcaption}
\usepackage{pdflscape}

\makeatletter
\def\l@subsection{\@tocline{2}{0pt}{2.5pc}{5pc}{}}
\def\l@subsubsection{\@tocline{2}{0pt}{5pc}{7.5pc}{}}
\makeatother

\theoremstyle{plain}
\newtheorem{theorem}{\bf Theorem}[section]
\newtheorem{lemma}[theorem]{\bf Lemma}
\newtheorem{proposition}[theorem]{\bf Proposition}
\newtheorem{corollary}[theorem]{\bf Corollary}

\theoremstyle{definition}
\newtheorem{definition}[theorem]{Definition}

\theoremstyle{remark}
\newtheorem{remark}[theorem]{Remark}

\numberwithin{equation}{section}

\newcommand{\NN}{\mathbb{N}}

\newcommand{\E}{\mathcal{E}}
\newcommand{\sym}{\mathfrak{S}}

\newcommand{\antiexcl}{\rotatebox[origin=c]{180}{!}}

\renewcommand{\P}{\mathcal{P}}
\newcommand{\Perm}{\mathcal{P}erm}
\newcommand{\Com}{\mathcal{C}om}
\newcommand{\Prelie}{\mathcal{P}re\mathcal{L}ie}
\newcommand{\Lie}{\mathcal{L}ie}

\newcommand{\cover}{\lessdot}

\newcommand{\Flyn}{\mathcal{FL}yn}

\newcommand{\clr}{\mathbf{color}}

\newcommand{\ISF}{\mathcal{ISF}}
\newcommand{\SF}{\mathcal{SF}}

\newcommand{\bpi}{\pmb{\pi}}

\newcommand{\balpha}{\pmb{\alpha}}

\newcommand{\tre}[1]{{\color{red}{#1}}}
\newcommand{\tbl}[1]{{\color{blue}{#1}}}

\def\newop#1{\expandafter\def\csname #1\endcsname{\mathop{\rm #1}\nolimits}}

\newop{sort}
\newop{sgn}

\author[R. S. Gonz\'alez D'Le\'on]{Rafael S. Gonz\'alez D'Le\'on }
\address[R.\ S.\ Gonz\'alez D'Le\'on]{Department of Mathematics and Statistics \\Loyola University Chicago\\ Chicago, IL 60660, USA} 
\email{rgonzalezdleon@luc.edu}
\urladdr{\url{http://dleon.combinatoria.co}}

\author[J. Hallam]{Joshua Hallam}
\address[J. Hallam]{Department of Mathematics, Loyola Marymount University,  Los Angeles, CA 90045, USA}
\email{Joshua.Hallam@lmu.edu}
\urladdr{\url{https://jhallam.lmu.build/}}

\author[Y. Quiceno]{Yeison A. Quiceno D. }
\address[Y. Quiceno]{Department of Statistics, University of Florida, Gainesville, FL 32611, USA}
\email{yeison.quicenodu@ufl.edu}

\begin{document}
\title{Whitney twins, Whitney duals, and operadic partition posets}
\maketitle

\begin{abstract}
We say that a pair of nonnegative integer sequences $(\{a_k\}_{k\ge 0},\{b_k\}_{k\ge 0})$ is Whitney-realizable if there exists a poset $P$ for which (the absolute values) of the Whitney numbers of the first and second kind are given by the numbers $a_k$ and $b_k$ respectively. The pair is said to be Whitney-dualizable if, in addition, there exists another poset $Q$ for which their Whitney numbers of the first and second kind are instead given by $b_k$ and $a_k$ respectively. In this case, we say that $P$ and $Q$ are Whitney duals. We use results on Whitney duality,  recently developed by the first two authors, to exhibit a family of sequences which allows for multiple realizations and Whitney-dual realizations. More precisely, we study edge labelings for the families of posets of pointed partitions $\Pi_n^{\bullet}$ and weighted partitions $\Pi_n^{w}$ which are associated to the operads $\Perm$ and $\Com^2$ respectively. The first author and Wachs proved that these two families of posets share the same pair of Whitney numbers. We find EW-labelings for $\Pi_n^{\bullet}$ and $\Pi_n^{w}$ and use them to show that they also share multiple nonisomorphic Whitney dual posets.

In addition to EW-labelings, we also find two new EL-labelings for $\Pi_n^\bullet$ answering a question of Chapoton and Vallette. Using these EL-labelings of $\Pi_n^\bullet$, and an EL-labeling of $\Pi_n^w$ introduced by the first author and Wachs, we  give combinatorial descriptions of bases for the operads $\Prelie, \Perm,$ and $\Com^2$. We also show that the bases for $\Perm$ and $\Com^2$ are PBW bases.

\end{abstract}

\setcounter{secnumdepth}{3}

\tableofcontents

\section{Introduction}

To a finite graded \emph{poset}  (partially ordered set) $P$ with a minimal element (denoted $\hat{0}$ throughout) we can associate a pair of sequences of integers $\{w_k(P)\}_{k\ge 0}$ and $\{W_k(P)\}_{k\ge 0}$ known as the \emph{Whitney numbers of the first and second kind} respectively. These two sequences  are poset 
 invariants and encode relevant information in areas where partially ordered structures arise naturally. For example, Whitney showed in \cite{Whitney1932} that the coefficients of the  chromatic polynomial of a graph are the Whitney numbers of the first kind of a poset one can associate to a graph (its bond lattice). The Whitney numbers of the first kind  keep track of the M\"obius function at each rank level and the Whitney numbers of the second kind keep track of the number of elements at each rank level.

 \subsection{Whitney-realizable and dualizable sequences}
 In~\cite{GonzalezHallam2021}, the first and second authors introduced the concept of a Whitney dual of a graded poset $P$ with a $\hat{0}$.
 We say that two graded posets $P$ and $Q$ are \emph{Whitney duals} if, after taking absolute values, the sequences of Whitney numbers of the first and second kind of $P$ are equal to the sequences of Whitney numbers of the second and first kind of $Q$. That is, the Whitney numbers of $P$ and $Q$ are swapped with respect to one another. In  \cite{GonzalezHallam2021}, the authors also defined a new type of poset edge labeling, which is called an EW-labeling (or Whitney labeling). The authors show that these labelings provide a sufficient condition for the existence of a Whitney dual for any graded poset $P$ admitting such a labeling.  Moreover, they describe an explicit construction of the Whitney dual associated to a given EW-labeling.

 One can readily observe from the definition, that nothing prevents the existence of multiple Whitney duals to a graded poset $P$. Hence, the concept of Whitney duality is more precisely a duality between the sequences of numbers involved rather than a duality between posets. We say that a pair of nonnegative integer sequences $(\{a_k\}_{k\ge 0},\{b_k\}_{k\ge 0})$ is \emph{Whitney-realizable} if there exists a poset $P$ such that $(\{|w_k(P)|\}_{k\ge 0},\{W_k(P)\}_{k\ge 0})=(\{a_k\}_{k\ge 0},\{b_k\}_{k\ge 0})$. We will call two posets $P$ and $Q$ \emph{Whitney twins} if they realize the same pair of sequences. We say that a Whitney-realizable pair is \emph{Whitney-dualizable} if $(\{b_k\}_{k\ge 0},\{a_k\}_{k\ge 0})$ is also Whitney-realizable.

 Determining which pairs of nonnegative integer sequences $(\{a_k\}_{k\ge 0},\{b_k\}_{k\ge 0})$ are Whitney-realizable or Whitney-dualizable both seem to be challenging questions. In this article, we present results related to the non-uniqueness of Whitney realizations and dualizations of a pair $(\{a_k\}_{k\ge 0},\{b_k\}_{k\ge 0})$ by finding and exploring the algebraic and combinatorial consequences of EW-labelings on two families of posets which come from the theory of symmetric operads. The two particular families of posets are associated to the permutative operad $\Perm$ and to the double commutative operad $\Com^2$.

 \subsection{Operadic posets and EL/CL-labelings}
In \cite{Vallette2007},  Vallete  defined a family of partition posets $\Pi_n^{\P}$ associated to a basic-set quadratic operad $\P$. 
 These posets are an operadic generalization of the poset of set partitions $\Pi_n$ ordered by refinement. There, the author shows  that the top cohomology $H^{top}(\Pi_n^{\P})$ $\sym_n$-modules are, up to tensoring with the sign representation, equal to the Koszul dual cooperad $\P^{\antiexcl}$ to $\P$. He also shows that the Cohen-Macaulay property of the maximal intervals of $\Pi_n^{\P}$ is equivalent to the Koszul property of $\P$ and $\P^{\antiexcl}$. Hence, the application of combinatorial techniques on the family $\Pi_n^{\P}$ is relevant in determining the algebraic properties of $\P$ and $\P^{\antiexcl}$. One such technique is the theory of lexicographic shellability for posets introduced by Bj\"orner \cite{Bjorner1980} and further developed by Bj\"orner and Wachs in \cite{BjornerWachs1982, BjornerWachs1983} (see also \cite{BjornerWachs1996,BjornerWachs1997}). The main idea behind the theory of lexicographic shellability is that the maximal intervals of a poset $P$ which admit a type of edge labeling, known as an EL-labeling (or a CL-labeling in more generality), are Cohen-Macaulay. Finding an EL or CL-labeling for a poset  $\Pi_n^{\P}$ then implies under Vallette's relation that  $\P$ and $\P^{\antiexcl}$ are Koszul operads.  As an application of EL and CL-labelings for partition posets,   Bellier-Mill{\`e}s, Delcroix-Oger, and Hoffbeck~\cite{BelliermillesDelcroixogerHoffbeck2021} showed that if an EL or CL-labeling of $\Pi_n^{\P}$ satisfies a certain condition that they call being \emph{isomorphism-compatible}, then the operad $\P$ has a Poincaré–Birkhoff–Witt (PBW) basis determined by the labeling.  PBW bases are useful because they imply that the operads are Koszul as was shown by Hoffbeck \cite{Hoffbeck2010} for totally ordered PBW bases and  in more generality for partially ordered PBW bases in~\cite{BelliermillesDelcroixogerHoffbeck2021}.
 
  We note that the posets $\Pi_n^{\P}$ 
 have appeared before in a different but related context. They are relevant in finding compositional (or substitutional) inverses to species within Joyal's theory of combinatorial species (see\cite{Joyal1981,BergeronLabelleLeroux1998}) as was shown  by M\'endez and Yang in \cite{MendezYang1991}.

\subsection{Pointed and weighted partition posets}

  Vallette \cite{Vallette2007} showed that the pointed partition poset $\Pi_{n}^{\bullet}$ is isomorphic to the operadic poset $\Pi_n^{\Perm}$  associated to the operad $\mathcal{P}erm$. In Section~\ref{sec:EWLabOfPointPart}, we give an EW-labeling of $\Pi_n^\bullet$ and give an explicit description of its Whitney dual in terms of pointed Lyndon forests in Section~\ref{subsection:whitney_dual_pointed}.

  In \cite{ChapotonVallete2006}, Chapoton and Vallette show that the maximal intervals of  $\Pi_{n}^{\bullet}$ are totally semimodular. By the results in \cite{BjornerWachs1983}, this implies that they are also CL-shellable and hence Cohen-Macaulay. By the result in \cite{Vallette2007} this in turn implies that  $\mathcal{P}erm$, and its Koszul-dual operad $\mathcal{P}re\mathcal{L}ie$, are Koszul.  The authors in \cite{ChapotonVallete2006} leave  open the question of whether or not $\Pi_{n}^{\bullet}$ admits the more restrictive property of being EL-shellable. EL-shellability and CL-shellability have been shown recently by Li \cite{Tiansi2020} to  not be equivalent in general for posets. The authors in \cite{BelliermillesDelcroixogerHoffbeck2021} propose a possible EL-labeling of $\Pi_{n}^{\bullet}$ and claim that this labeling has the additional property of being isomorphism-compatible. We show in Section \ref{section:El-labeling_pointed_partition} that the proposed labeling does not satisfy the requirements for being an EL-labeling. We then provide a new EL-labeling which answers the open question in \cite{ChapotonVallete2006}. This labeling has the same set of labels as our EW-labeling for $\Pi_{n}^{\bullet}$, but differ in how these labels are partially ordered. We show this EL-labeling is isomorphism-compatible which in turn gives a PBW basis for the $\Perm$ operad using the results in \cite{BelliermillesDelcroixogerHoffbeck2021}. Although our EW-labeling for $\Pi_{n}^{\bullet}$ is not directly an EL-labeling, we show that reversing the order on the labels gives an EL-labeling for the order dual poset $(\Pi_{n}^{\bullet})^{*}$. This provides a second answer to the open question in \cite{ChapotonVallete2006}.  We also show that the former EL-labeling for $\Pi_{n}^{\bullet}$  is isomorphism-capatible, giving us a PBW bases for $\Perm$.

  In \cite{DotsenkoKhoroshkin2007}, Dotsenko and Khoroshkin introduced the weighted partition poset $\Pi_n^w$. They  showed that $\Pi_n^w$ is isomorphic to the poset $\Pi_n^{\Com^2}$ associated to the operad $\Com^2$ of algebras with two totally commutative products.  The combinatorial and homological properties of $\Pi_n^w$ were extensively studied by Gonz\'alez D'Le\'on and Wachs in \cite{DleonWachs2016}. In their study, the authors introduced an EL-labeling for $\Pi_n^w$. In Section~\ref{sec:EWforWeighted} we prove that this labeling is an EW-labeling and hence $\Pi_n^w$ has a Whitney dual.   In Section~\ref{section:WhitneyDualOfWeighted}, we give an explicit description of this Whitney dual in terms of bicolored Lyndon forests.  We also show in Section~\ref{sec:algConseq} that this labeling is isomorphism-compatible which gives a PBW basis for $\Com^2$.
  
  \subsection{Nonuniqueness of Whitney realizations}
  In \cite[Section 2.4]{DleonWachs2016} the authors show that $\Pi_n^w$ and $\Pi_n^{\bullet}$ are Whitney twins (though they do not use this terminology). Indeed as a consequence of their Theorem 2.8, Proposition 2.1, and the follow up discussion in Section 2.4 in \cite{DleonWachs2016},  the Whitney numbers of the first and second kind are given for all $k\ge 0$ by the sequences
  \begin{align*}
w_k(\Pi_n^{\bullet})=w_k(\Pi_n^{w})&=(-1)^{k}\binom{n-1}{k}n^{k}\\
W_k(\Pi_n^{\bullet})=W_k(\Pi_n^{w})&=\binom{n}{k}(n-k)^{k}.
  \end{align*}

  This already implies the nonuniqueness of realizations for a Whitney-realizable sequence. We show that the Whitney duals constructed with the EW-labelings for $\Pi_n^\bullet$ and $\Pi_n^w$ are not isomorphic for $n\geq 4$.  Since they have the same Whitney numbers of both kinds, we get multiple non-isomorphic Whitney duals for both $\Pi_n^\bullet$ and $\Pi_n^w$, implying further the nonuniqueness of dual realizations of Whitney-dualizable sequences. We also show that there is a third family $\SF_n$ of Whitney duals to $\Pi_n^\bullet$ and $\Pi_n^w$ which for $n\ge 3$ is not isomorphic to any of the Whitney duals discussed before. The family  $\SF_n$ is also shown in future work by the first two authors to be associated with a more general type of Whitney labeling. The three nonisomorphic families of Whitney dual posets to $\Pi_n^\bullet$ and $\Pi_n^w$ also constitute a new example of the nonuniqueness of Whitney realizations.

\subsection{Organization of this article}

The rest of the article is structured as follows. In Section~\ref{sec:EWLab} we review EW-labelings and EL-labelings, and we describe the labelings of $\Pi_n^w$ and $\Pi_n^\bullet$.  In Section~\ref{sec:CombDescriptOfWhitneyDual}, we give explicit descriptions of the Whitney duals of $\Pi_n^w$ and $\Pi_n^\bullet$.  In Section~\ref{sec:algConseq} we consider the algebraic consequences of these labelings. Specifically we use these labelings to describe bases for $\Prelie$, $\Perm$, and $\Com^2$, the latter two in particular being PBW bases.  In Section~\ref{sec:twins}, we discuss  the nonuniqueness  of Whitney realizations using our results for $\Pi_n^w$ and $\Pi_n^\bullet$ and their associated Whitney duals. 

Some results in this work have been announced as part of the third author's master's thesis in \cite{QuicenoDuran2020}.

\section{EW-labelings}\label{sec:EWLab}
In this section we describe three edge labelings: one for the weighted partition poset, which was introduced already in \cite{DleonWachs2016}, and two new edge labelings for the pointed partition poset. The edge labeling for the weighted partition poset, was shown in \cite{DleonWachs2016} to be an EL-labeling and here we show that is also an EW-labeling. Of the two labelings for the pointed partition poset, one is an EW, which we also show is a dual EL-labeling, and the second is an EL-labeling (but not an EW-labeling). We show that the two labelings have  the same sets of words of labels, however the labels come from two different partial orders.   Our main use of these labelings is three-fold: constructing Whitney duals for the two posets, understanding their homotopy type and cohomology of the respective order complexes, and finding PBW bases of the corresponding operads and bases for their dual (co)operads.    We start with a brief review of Whitney numbers, Whitney duals, and  edge labelings.

\subsection{Whitney numbers and Whitney duals}

We will assume some familiarity with posets. For a more in-depth review of posets as well as any undefined terms, see~\cite[Chapter 3]{Stanley2012}.  
For a review of poset topology see~\cite{Wachs2007}.
All the posets we consider in this article will be finite, graded, and contain a minimum element which we denote by $\hat{0}$.  We will use $\rho(x)$ for the rank of an element $x$.

 The (one-variable) \emph{M\"obius function} of a poset $P$, denoted by $\mu$,  is defined recursively by
 $$
 \mu(\hat{0}) =1
 $$
and for $x\neq \hat{0}$,
$$
\mu(x) = -\sum_{y<x} \mu(y).
$$

\begin{figure}
    \begin{tikzpicture}
    \node (0) at (0,0) {$\hat{0}$};  
      \node (a) at (-2,2) {$a$};  
    \node (b) at (0,2) {$b$};  
    \node (c) at (2,2) {$c$};  
    \node (1) at (0,4) {$\hat{1}$}; 
    
      \node at (0,-1) {$P$};
         \draw[thick]  (0)--(a) node [midway,sloped, below] {$\textcolor{blue}{a}$};
        \draw[thick]  (a)--(1) node [midway,sloped, above] {$\textcolor{blue}{b}$};
        
         \draw[thick]  (0)--(b) node [midway, sloped,above] {$\textcolor{blue}{b}$};
        \draw[thick]  (b)--(1) node [midway, sloped, above] {$\textcolor{blue}{a}$};

         \draw[thick]  (0)--(c) node [midway,sloped, below] {$\textcolor{blue}{c}$};
        \draw[thick]  (c)--(1) node [midway,sloped, above] {$\textcolor{blue}{a}$};

     \node    at (0-.5,0) {$\textcolor{red}{+1}$};  
      \node  at (-2-.5,2)  {$\textcolor{red}{-1}$};  
    \node  at (0-.5,2)   {$\textcolor{red}{-1}$};  
    \node   at (2-.5,2)    {$\textcolor{red}{-1}$};  
    \node  at (0-.5 ,4)  {$\textcolor{red}{+2}$};

    \begin{scope}[shift={(6,0)}]
    \node (0) at (0,0) {$(\hat{0},\emptyset)$};  
      \node (a) at (-2,2) {$(a,a)$};  
    \node (b) at (0,2) {$(b,b)$};  
    \node (c) at (2,2) {$(c,c)$};  
    \node (1a) at (-1,4) {$(\hat{1},ba)$}; 
    \node (1b) at (2,4) {$(\hat{1},ca)$};

      \node at (0,-1) {$Q$};
         \draw[thick]  (0)--(a) node [midway,sloped, below] {$\textcolor{blue}{a}$};
        \draw[thick]  (a)--(1a) node [midway,sloped, above] {$\textcolor{blue}{b}$};
        
         \draw[thick]  (0)--(b) node [midway, sloped,above] {$\textcolor{blue}{b}$};
        \draw[thick]  (b)--(1a) node [midway, sloped, above] {$\textcolor{blue}{a}$};

         \draw[thick]  (0)--(c) node [midway,sloped, below] {$\textcolor{blue}{c}$};
        \draw[thick]  (c)--(1b) node [midway,sloped, above] {$\textcolor{blue}{a}$};

       \node    at (0-.75,0) {$\textcolor{red}{+1}$};  
       \node  at (-2-.75,2)  {$\textcolor{red}{-1}$};  
     \node  at (0-.75,2)   {$\textcolor{red}{-1}$};  
     \node   at (2-.75,2)    {$\textcolor{red}{-1}$};  
     \node  at (-1-.85,4)  {$\textcolor{red}{+1}$};  
    \node  at (2-.85,4)  {$\textcolor{red}{0}$};  
    
    \end{scope}
    
    \end{tikzpicture}
    \caption{Two posets which are Whitney duals.  Values (in red) besides elements correspond to the M\"obius function and those beside edges (in blue) correspond to edge labels.  The edge labels are ordered alphabetically.}\label{fig:WhitneyDualExample}
\end{figure}
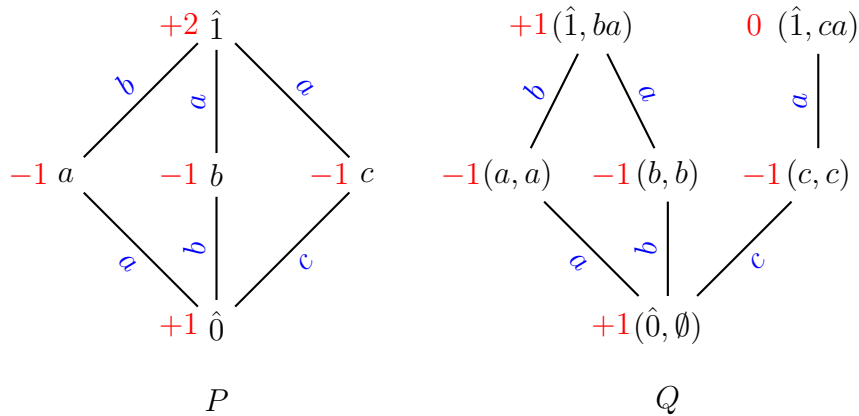

Note that the one-variable M\"obius function coincides with the classical two-variable M\"obius function $\mu(\hat{0},x)$ on the interval $[\hat{0},x]$. See Figure~\ref{fig:WhitneyDualExample}  for examples of the M\"obius function.    The \emph{$k^{th}$ Whitney number of the first kind}, denoted by $w_k(P)$,  is defined by
$$
w_k(P) = \sum_{\substack{x\in P\\ \rho(x) = k}} \mu(x).
$$
For the poset  $P$ in Figure~\ref{fig:WhitneyDualExample}, the Whitney numbers of the first kind are given by the sequence $(1, -3 ,2)$ and for the poset $Q$ these are given by the sequence $(1,-3,1)$.  

 The \emph{$k^{th}$ Whitney number of the second kind}, denoted by $W_k$, is defined by 
 $$
 W_k(P) =  \#\{x\in P \mid \rho(x)=k\}.
 $$
In Figure~\ref{fig:WhitneyDualExample}, the Whitney numbers of the second  kind of $P$ are given by $(1, 3, 1)$ and of $Q$ are given by $(1, 3 , 2)$.  

By comparing the Whitney numbers of the first and second kind of $P$ and $Q$ in Figure~\ref{fig:WhitneyDualExample}, the reader may notice a peculiar phenomenon. The Whitney numbers of $P$ and $Q$  switch (up to a sign).  It turns out that this phenomenon, which was first described in~\cite{GonzalezHallam2017} and further studied in~\cite{GonzalezHallam2021}, occurs for many other pairs of posets and motivates the next definition.

\begin{definition}
Let $P$ and $Q$ be ranked posets.  We say $P$ and $Q$ are  \emph{Whitney duals} if for all $k$,
$$
|w_k(P)| = W_k(Q) \mbox{ and }  |w_k(Q)| = W_k(P).
$$
\end{definition}

\noindent Using this definition, we can see that $P$ and $Q$ in Figure~\ref{fig:WhitneyDualExample} are Whitney duals.

\subsection{Consequences of ER, EL, and EW-labelings}

To approach Whitney duality, the first two authors in ~\cite{GonzalezHallam2021} used the poset topology technique of edge labelings.    Here, we review  a few concepts related to edge labelings, but for further details the reader can visit \cite{Wachs2007,Stanley2012} or some of the classical articles \cite{Stanley1974,Bjorner1980,BjornerWachs1983,BjornerWachs1996,BjornerWachs1997}

Let $P$ be a poset. We use $\mathcal{E}(P)$ to denote the set of edges in the Hasse diagram of $P$ (which is  in bijection with the set of cover relations in $P$).  An \emph{edge labeling} of $P$ is a map $\lambda:  \mathcal{E}(P)\rightarrow \Lambda$ where $\Lambda$ is a set of partially ordered labels.    See Figure~\ref{fig:WhitneyDualExample} for   examples of an edge labeling where the set of labels is $\{a,b,c\}$ which is ordered alphabetically.  Recall that a chain $x_0\cover x_1\cover \cdots\cover x_n$ is said to be saturated if it is maximal in the interval $[x_0,x_n]$. Given an edge labeling $\lambda$, we say  that a saturated chain, $x_0\cover x_1\cover \cdots\cover x_n$,  is \emph{increasing} if $\lambda(x_{i-1}\cover x_{i}) < \lambda(x_{i}\cover x_{i+1})$ for all $1\leq i\leq n-1$.  Similarly, $x_0\cover x_1\cover \cdots\cover x_n$,  is \emph{ascent-free} if $\lambda(x_{i-1}\cover x_{i}) \not < \lambda(x_{i}\cover x_{i+1})$ for all $1\leq i\leq n-1$.  Returning to our example    in   Figure~\ref{fig:WhitneyDualExample}, we see that among maximal chains of $P$, the chain $\hat{0}\cover a\cover \hat{1}$ is increasing (since $ab$ is an increasing sequence).  On the other hand, the maximal chains $\hat{0}\cover b \cover \hat{1}$ and $\hat{0}\cover c\cover \hat{1}$ are ascent-free.  We want to remark that the example in  Figure~\ref{fig:WhitneyDualExample} is rather small and in general there are saturated chains that are neither increasing nor ascent-free, but these two particular types of chains are the ones of interest in the following discussion.

\subsubsection{ER and EL-labelings}
We say an edge labeling is an \emph{ER-labeling} if every interval has a unique increasing maximal chain.  Moreover, we say an ER-labeling is an \emph{EL-labeling} if in each interval, the unique increasing maximal chain also precedes every other chain in lexicographic order.  One can check that the labeling of $P$ in   Figure~\ref{fig:WhitneyDualExample} is both an ER and an EL-labeling. Indeed, the lexicographic requirement holds trivially on rank $0$ and $1$ intervals, so the only interval to check the  lexicographic condition is on the full poset (which is also an interval in this case).  The increasing chain is labeled $ab$ and this precedes both $ba$ and $ca$ in lexicographic order.   One of the main reasons we are interested in ER and EL-labelings is because of the topological and combinatorial consequences given by the following two theorems. 

\begin{theorem}[Stanley~\cite{Stanley1974}]
\label{theorem:stanley_mobius}
  Let $P$ be a graded poset  with an ER-labeling $\lambda: \mathcal{E}(P) \rightarrow \Lambda$. Then  for  every $x<y$ in $P$  we have that
     $$
     \mu(x,y) = (-1)^{\rho([x,y])} |\{c\mid c \mbox{ an ascent-free maximal chain in  }[x,y]\}|.
     $$
\end{theorem}

\begin{theorem}  [Bj\"orner and Wachs~\cite{BjornerWachs1983}]\label{theorem:bjorner_wachs}
  Let $P$ be a graded poset  with an EL-labeling $\lambda: \mathcal{E}(P) \rightarrow \Lambda$. Then for  every $x<y$ in $P$  we have that:
 \begin{enumerate}
   \item   The order complex $\Delta((x,y))$ is shellable. Moreover, it has the homotopy type of a wedge of $|\{c\mid c\mbox{ an ascent-free maximal chain in }[x,y]\}|$ many spheres each of dimension $\rho([x,y])-2$. As a consequence, $[x,y]$ is Cohen-Macaulay.

\item The set
$$ \{c\setminus \{x,y\} \mid c \mbox{ an ascent-free maximal  chain in  }[x,y]\}$$ forms a basis for the top reduced cohomology $\widetilde H^{\rho([x,y])-2}((x,y))$ of  $\Delta((x,y))$.
 \end{enumerate}
\end{theorem}

\begin{remark}
   In this work we will be particularly interested on the consequences of Theorems \ref{theorem:stanley_mobius} and \ref{theorem:bjorner_wachs} for intervals of the form $[\hat{0},x]$ for all $x$ in a poset $P$.
\end{remark}

\subsubsection{EW-labelings}
In order to construct a Whitney dual, we need to impose two additional conditions on an ER-labeling.   Note that in the following definition, we do not require the labeling to be an EL-labeling.  

 \begin{definition}\label{definition:EW-labeling}
 Let $\lambda$ be an edge labeling of $P$.  We say $\lambda$ is an \emph{EW-labeling} if the following hold.
 \begin{enumerate}
    \item  $\lambda$ is an ER-labeling.
     \item ({\bf The rank two switching property}) For every interval $[x,y]$ with $\rho(y)-\rho(x)=2$, if the increasing chain is labeled $ab$, there exists a unique chain in $[x,y]$ labeled  $ba$.  
     \item ({\bf Injectivity of ascent-free chains}) For every $x<y\in P$, every ascent-free maximal chain in $[x,y]$ has a unique sequence of labels.
\end{enumerate}
\end{definition}

We already noted that the labeling of  $P$ in Figure~\ref{fig:WhitneyDualExample} is an ER-labeling.     In fact, it is an EW-labeling too.  Clearly we have injectivity of ascent-free chains.     Moreover, in the only rank two interval, the increasing chain is labeled by $ab$ and there is exactly one other chain in that interval labeled by $ba$.  As we saw, the poset $P$ has a Whitney dual (namely $Q$ in Figure~\ref{fig:WhitneyDualExample}).  This is no coincidence, rather it is a a consequence of the following theorem.

\begin{theorem}[Theorem 1.6 \cite{GonzalezHallam2021}] 
  Let $P$ be a poset with an EW-labeling $\lambda$. Then $P$ has a Whitney dual. Moreover, we can construct a Whitney dual $Q$ to $P$ that depends on $\lambda$.  
\end{theorem}

In Section~\ref{sec:constructingWD}, we describe a specific construction of such a Whitney dual $Q$ using $\lambda$. 

\subsection{An EW-labeling of the weighted partition poset}\label{sec:EWforWeighted}

In this subsection, we describe an EW-labeling of the weighted partition poset.  First, we briefly discuss the weighted partition poset.

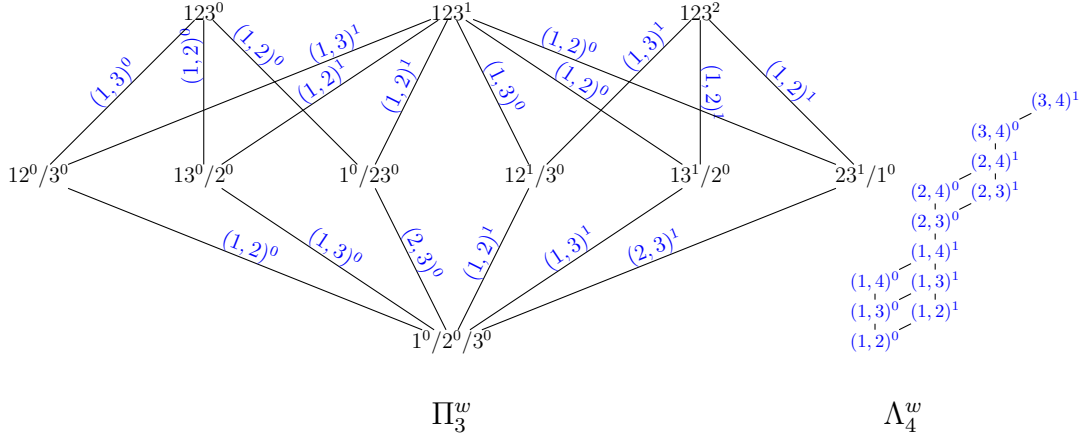
\begin{figure}
\begin{center} 
\begin{tikzpicture}
\node at (0,-1) {$\Pi_3^w$};
\node at (6,-1) {$\Lambda^w_4$};

\begin{scope}[xshift=0,scale=1.1]
\tikzstyle{every node}=[inner sep=0pt, scale=0.7, minimum width=4pt]
\node (n1232) at (3,4) {$123^{2}$};
  \node (n13020) at (-3,2) {$13^{ 0}/ 2^{ 0}$};
  \node (n102030) at (0,0)  {$1^{0}/ 2^{0}/ 3^{0}$};
  \node (n1231) at (0,4) {$123^{1}$};
  \node (n12030) at (-5,2) {$12^{0}/ 3^{0}$};
  \node (n13120) at (3,2)  {$13^{1}/ 2^{0}$};
  \node (n1230) at (-3,4){$123^ {0}$};
  \node (n10230) at (-1,2)  {$1^{0}/ 23^{0}$};
  \node (n12130) at (1,2)  {$12^{ 1}/ 3^{0}$};
  \node (n10231) at (5,2) {$23^ {1}/ 1^{0}$};


  \draw (n1231) -- (n10230)  node [midway,sloped, above] {$\textcolor{blue}{(1,2)^1}$};
  \draw [] (n13020) -- (n102030)   node [midway,sloped, above] {$\textcolor{blue}{(1,3)^0}$};
  \draw [] (n1232) -- (n13120)  node [midway,sloped, above] {$\textcolor{blue}{(1,2)^1}$};
  \draw [] (n1231)-- (n13020)node [midway,sloped, above] {$\textcolor{blue}{(1,2)^1}$};
  \draw [] (n10230)--(n102030) node [midway,sloped, above] {$\textcolor{blue}{(2,3)^0}$};
  \draw [] (n1230) -- (n10230) node [near start, sloped, above] {$\textcolor{blue}{(1,2)^0}$};
  \draw [] (n1231) -- (n13120) node [midway,sloped, above] {$\textcolor{blue}{(1,2)^0}$};
  \draw [] (n12030)-- (n102030) node [midway,sloped, above] {$\textcolor{blue}{(1,2)^0}$};
  \draw [] (n1231) --(n12130)node [midway,sloped, above] {$\textcolor{blue}{(1,3)^0}$};
  \draw [] (n1232) -- (n12130) node [near start,sloped, above] {$\textcolor{blue}{(1,3)^1}$};
  \draw [] (n13120) --(n102030) node [midway,sloped, above] {$\textcolor{blue}{(1,3)^1}$};
  \draw [] (n1231) --(n10231) node [near start,sloped, above] {$\textcolor{blue}{(1,2)^0}$};
  \draw [] (n1230) -- (n13020) node [near start,sloped, above] {$\textcolor{blue}{(1,2)^0}$};
  \draw [] (n1230)  -- (n12030) node [midway,sloped, above] {$\textcolor{blue}{(1,3)^0}$};
  \draw [] (n12130) --  (n102030) node [midway,sloped, above] {$\textcolor{blue}{(1,2)^1}$};
  \draw [] (n1232)  --  (n10231) node [midway,sloped, above] {$\textcolor{blue}{(1,2)^1}$};
  \draw [] (n10231)  --  (n102030)node [midway,sloped, above] {$\textcolor{blue}{(2,3)^1}$};
  \draw [] (n1231) -- (n12030) node [near start,sloped, above] {$\textcolor{blue}{(1,3)^1}$};

\tikzstyle{every node}= [scale=0.6]

\end{scope}
\begin{scope}[xshift=160,scale=0.4]
 \tikzstyle{every node}=[inner sep=1pt, minimum width=14pt,scale=0.7, font=\footnotesize]
\draw (0,0) node (n120) {\color{blue}$(1,2)^0$};
\draw (0,1) node (n130) {\color{blue}$(1,3)^0$};
\draw (0,2) node (n140) {\color{blue}$(1,4)^0$};
\draw (2,1) node (n121) {\color{blue}$(1,2)^1$};
\draw (2,2) node (n131) {\color{blue}$(1,3)^1$};
\draw (2,3) node (n141) {\color{blue}$(1,4)^1$};

\draw (n141) -- (n140) ;
\draw (n131) -- (n130) ;
\draw (n121) -- (n120) ;
\draw (n140)-- (n130) -- (n120) ;

\draw (2,4) node (n230) {\color{blue}$(2,3)^0$};
\draw (2,5) node (n240) {\color{blue}$(2,4)^0$};
\draw (4,5) node (n231) {\color{blue}$(2,3)^1$};
\draw (4,6) node (n241) {\color{blue}$(2,4)^1$};
\draw (4,7) node (n340) {\color{blue}$(3,4)^0$};
\draw (6,8) node (n341) {\color{blue}$(3,4)^1$};
\draw (n241) -- (n240) ;
\draw (n231) -- (n230) ;

\draw (n240)-- (n230)  -- (n141) -- (n131) -- (n121);

\draw (n341) -- (n340) ;

\draw (n340) -- (n241) -- (n231);

\end{scope}

\end{tikzpicture}

\end{center}
\caption[]{Edge labeling of  $\Pi_3^w$ and the poset $\Lambda_4^w$ }\label{figure:weightedn3k2EL}
\end{figure}

A \emph{weighted set} is a pair $(A,v)$ where $v\in\{0,\dots, |A|-1\}$. We will also denote weighted sets with the simpler notation $A^v$. A \emph{weighted partition} of $[n]$ is a collection of weighted sets $\bpi=B_1^{v_1}/B_2^{v_2}/\cdots/B_t^{v_t}$  such that  
$B_1/B_2/\cdots/B_t$ is a partition of $[n]$. The \emph{poset of weighted partitions}, $\Pi_n^{w}$,  is the set of weighted 
partitions of $[n]$ with cover order relation given by 
$$\bpi=A_1^{w_1}/A_2^{w_2}/\dots /A_s^{w_s}\lessdot B_1^{v_1}/ B_2^{v_2}/ \dots /B_{s-1}^{v_{s-1}}=\bpi'$$ 
if the following conditions hold:
\begin{itemize}
 \item $A_1/A_2/\cdots /A_s \lessdot B_1/B_2/ \cdots /B_{s-1} $ in $\Pi_n$
 \item if $B_k=A_{i}\cup A_{j}$, where $i \ne j$, then $v_k-(w_{i} + w_{j}) \in \{0,1\}$
 \item if $B_k = A_i$ then $v_k = w_i$
 \end{itemize}
See Figure~\ref{figure:weightedn3k2EL} for a depiction of $\Pi_3^w$. As was noted in the introduction, $\Pi_n^w$ is (isomorphic to) the poset of partitions for the operad $\Com^2$.

In \cite{DleonWachs2016}, Gonz\'alez D'Le\'on and Wachs gave an EL-labeling for $\Pi_n^w$. Here we show that this labeling is in fact an EW-labeling.  We now review the definition of their labeling. 

Let us start by defining the set of edge labels, $\Lambda^w_n$. For each $a \in [n]$, let $\Gamma_a:= \{(a,b)^u \mid   
a<b \le n, \,\, u \in \{0,1\} \}$.  
We partially order $\Gamma_a$  by letting $(a,b)^u \le  (a,c)^v$ if $b\le  c$ and $u \le v$. Note that $\Gamma_a$ is isomorphic to the direct product of the chain $a+1< a+2 <\dots < n $ and
the chain $0 < 1$.  Now define $\Lambda^w_n$ to be the 
ordinal sum
$\Lambda^w_n := \Gamma_1 \oplus  \Gamma_2  \oplus \cdots \oplus \Gamma_{n-1}$. See Figure~\ref{figure:weightedn3k2EL} for a depiction of the Hasse diagram of $\Lambda^w_4$.

We are now ready to describe the edge labeling.   The map $\lambda_w:\E(\Pi_n^w)\rightarrow \Lambda_n^w$ is defined as follows: let $\bpi \lessdot \bpi' $ in 
$\Pi^w_n$ so that $\bpi'$ is obtained from $\bpi$ by merging two blocks $A^{w_A}$ and $B^{w_B}$ of $\bpi$ into a 
new block $(A \cup B)^{w_A + w_B+u}$ of $\bpi'$, where $u \in \{0,1\}$ and where we assume without loss of 
generality that $\min A < \min B$.  We 
define then
$$
\lambda_w(\bpi \cover \bpi') = (\min A, \min B)^u.
$$
See Figure~\ref{figure:weightedn3k2EL}  for an example of this labeling on $\Pi_3^w$.   

The following theorem was proved in \cite{DleonWachs2016}.

\begin{theorem}[\cite{DleonWachs2016} Theorem 3.2] $\lambda_w$ is an EL-labeling (and hence also an ER-labeling).
\end{theorem}

According to Definition \ref{definition:EW-labeling}, to show that $\lambda_w$ is an EW-labeling, we  need to check the rank two switching property and the injectivity condition on ascent-free chains.  To see that the latter is satisfied, note that the information contained in the collection of labels of the form $(\min A, \min B)^u$ is enough to trace which blocks of a weighted partition are merged at each step which is enough to recover any saturated chain starting at any particular weighted partition $\bpi$. Hence the sequence of labels in each interval 
uniquely determines a chain. 

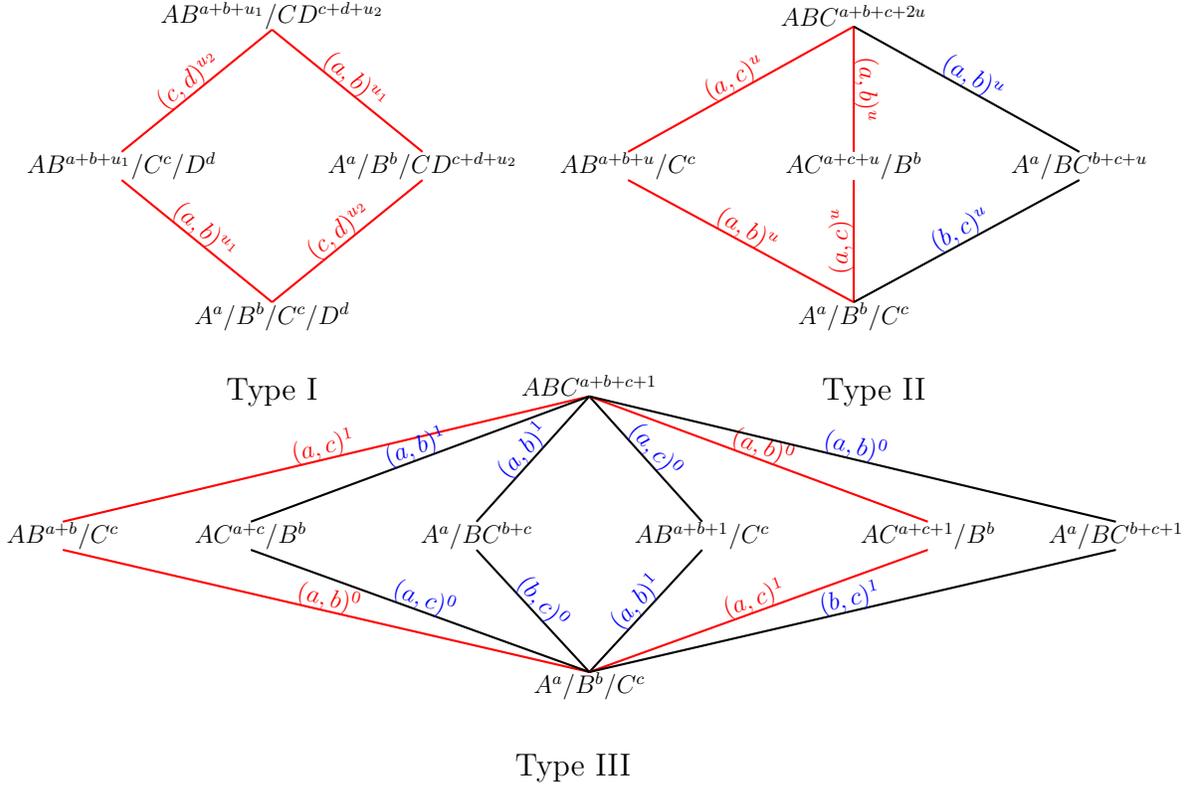
\begin{figure}
\centering
\begin{tikzpicture}[line join=bevel,scale=1]
  \node at (0,-1)  {Type I};
  \node at (8,-1)  {Type II};
  \node at (4,-6)  {Type III};

\begin{scope}
  \tikzstyle{every node}=[inner sep=0pt, scale=0.8, minimum width=4pt]
  \node (v1-2-3-4) at (0,0)  {$A^{a}/ B^{b}/ C^{c}/ D^{d}$};
  \node (v12-34) at (0,4)  {$AB^ {a+b+u_1}/CD^ {c+d+u_2}$};
  \node (v12-3-4) at (-2,2)  {$AB^{a+b+u_1}/ C^{c}/ D^{d}$};
  \node (v1-2-34) at (2,2)  {$A^{a}/ B^{b}/ CD^{c+d+u_2}$};
  \path[thick](v1-2-3-4.north)                                            edge[red] node[sloped,yshift=0.5em,draw=none,fill=none,color=red]{{$\Large{(a,b)^{u_1}}$}}
(v12-3-4.south)edge[red] node[sloped,yshift=0.5em,draw=none,fill=none,color=red]{{$\Large{(c,d)^{u_2}}$}}
(v1-2-34.south);
\path[thick](v12-34.south)                                            edge[red] node[sloped,yshift=0.5em,draw=none,fill=none,color=red]{{$\Large{(c,d)^{u_2}}$}}
(v12-3-4.north)edge[red] node[sloped,yshift=0.5em,draw=none,fill=none,color=red]{{$\Large{(a,b)^{u_1}}$}}
(v1-2-34.north);

 \end{scope}

\begin{scope}[xshift=220]
  \tikzstyle{every node}=[inner sep=0pt, scale=0.8, minimum width=4pt]
  \node (v1-2-3) at (0,0)  {$A^{a}/ B^{b}/ C^{c}$};
  \node (v123) at (0,4)  {$ABC^ {a+b+c+2u}$};
  \node (v12-3) at (-3,2)  {$AB^{a+b+u}/ C^{c}$};
  \node (v1-23) at (3,2)  {$A^{a}/ BC^{b+c+u} $};
 \node (v13-2) at (0,2)  {$AC^{a+c+u}/ B^{b}$};
 \path[thick](v1-2-3.north)                                            edge [red] node[sloped,yshift=0.5em,draw=none,fill=none,color=red]{{$\Large{(a,b)^u}$}}
(v12-3.south)edge[red] node[sloped,yshift=0.5em,draw=none,fill=none,color=red]{{$\Large{(a,c)^u}$}}
(v13-2.south)edge node[sloped,yshift=0.5em,draw=none,fill=none,color=blue]{{$\Large{(b,c)^u}$}}
(v1-23.south);
\path[thick](v123.south)                                            edge[red] node[sloped,yshift=0.5em,draw=none,fill=none,color=red]{{$\Large{(a,c)^u}$}}
(v12-3.north)edge[red] node[sloped,yshift=0.5em,draw=none,fill=none,color=red]{{$\Large{(a,b)^u}$}}
(v13-2.north)edge node[sloped,yshift=0.5em,draw=none,fill=none,color=blue]{{$\Large{(a,b)^u}$}}
(v1-23.north);

 \end{scope}

\begin{scope}[xshift=120,yshift=-140]
\tikzstyle{every node}=[ inner sep=0pt, scale=0.8, minimum width=4pt]
  \node (v1-2-3) at (0,0)  {$A^{a}/ B^{b}/ C^{c}$};
  \node (v123) at (0,4)  {$ABC^{a+b+c+1}$};
  \node (v12a-3) at (-7,2)  {$AB^{a+b}/ C^{c}$};
  \node (v13a-2) at (-4.5,2)  {$AC^{a+c}/ B^{b}$};
  \node (v1-23a) at (-1.5,2)  {$A^{a}/ BC^{b+c}$};   
  \node (v12b-3) at (1.5,2)  {$AB^{a+b+1}/ C^{c}$};
  \node (v13b-2) at (4.5,2)  {$AC^{a+c+1}/ B^{b}$};
  \node (v1-23b) at (7,2)  {$A^{a}/BC^ {b+c+1}$};
\path[thick](v1-2-3.north)                                            edge[red] node[sloped,yshift=0.5em,draw=none,fill=none,color=red]{{$\Large{(a,b)^0}$}}
(v12a-3.south)edge node[sloped,yshift=0.5em,draw=none,fill=none,color=blue]{{$\Large{(a,c)^0}$}}
(v13a-2.south)edge node[sloped,yshift=0.5em,draw=none,fill=none,color=blue]{{$\Large{(b,c)^0}$}}
(v1-23a.south)edge node[sloped,yshift=0.5em,draw=none,fill=none,color=blue]{{$\Large{(a,b)^1}$}}
(v12b-3.south)edge[red] node[sloped,yshift=0.5em,draw=none,fill=none,color=red]{{$\Large{(a,c)^1}$}}
(v13b-2.south)edge node[sloped,yshift=0.5em,draw=none,fill=none,color=blue]{{$\Large{(b,c)^1}$}}
(v1-23b.south);
\path[thick](v123.south)                                            edge[red] node[sloped,yshift=0.5em,draw=none,fill=none,color=red]{{$\Large{(a,c)^1}$}}
(v12a-3.north)edge node[sloped,yshift=0.5em,draw=none,fill=none,color=blue]{{$\Large{(a,b)^1}$}}
(v1-23a.north)edge node[sloped,yshift=0.5em,draw=none,fill=none,color=blue]{{$\Large{(a,c)^0}$}}
(v12b-3.north)edge node[sloped,yshift=0.5em,draw=none,fill=none,color=blue]{{$\Large{(a,b)^1}$}}
(v13a-2.north)edge[red] node[sloped,yshift=0.5em,draw=none,fill=none,color=red]{{$\Large{(a,b)^0}$}}
(v13b-2.north)edge node[sloped,yshift=0.5em,draw=none,fill=none,color=blue]{{$\Large{(a,b)^0}$}}
(v1-23b.north);

\end{scope}

\end{tikzpicture}

\label{fig:type1}
\caption[]{Rank two intervals in $\Pi_n^w$.   Here $A,B,C,D$ are the blocks that get merged in the interval and $a=\min(A)< b=\min(B)<c =\min(C)<d=\min(d)$.  The  blocks that are not changed in the interval are not depicted.   Edges given in red correspond to the rank two switching property. }\label{figure:ranktwointervals_weighted}
\end{figure}

To show that $\lambda_w$ is an EW-labeling, we are left to show that it 
satisfies the rank two switching property.
 As explained in~\cite{DleonWachs2016}, there are three types of rank two intervals in $\Pi_n^w$. These intervals are depicted in Figure~\ref{figure:ranktwointervals_weighted} together with their edge labels.   
 For each type, the reader can check that the rank two switching property holds. 

\begin{theorem}\label{theorem:lambdaE_EW}
 $\lambda_w$ is an EW-labeling of $\Pi_n^w$. Consequently, $\Pi_n^w$ has a Whitney dual.
\end{theorem}

We will give a combinatorial description of the corresponding Whitney dual in Section~\ref{section:WhitneyDualOfWeighted}.

\subsection{An EW-labeling of the pointed partition poset}\label{sec:EWLabOfPointPart}

A \emph{pointed set} is a pair $(A,p)$ where $A$ is a nonempty set and $p\in A$. In the following we will use the notation $A^p$ for   $(A,p)$.  A \textit{pointed partition} of $[n]$ is a collection $\bpi=\{B_1^{p_1},B_2^{p_2},\dots,B_m^{p_m}\}$ where $\pi=\{B_1,B_2,\dots,B_m\}$ is a partition of $[n]$, called its \emph{underlying partition}, and $B_i^{p_i}$ are pointed sets for all $i$.
We will also use the notation $B_1^{p_1}/B_2^{p_2}/\cdots/ B_m^{p_m}$ for $ \{B_1^{p_1},B_2^{p_2},\dots,B_m^{p_m}\}$. The \emph{poset of pointed partitions} $\Pi_{n}^{\bullet}$ is the partial order on the set of all pointed partitions of $[n]$ with cover order relation given by $\bpi=\{A_1^{q_1},A_2^{q_2},...,A_l^{q_l}\}\lessdot \bpi'=\{B_1^{p_1},B_2^{p_2},...,B_m^{p_m}\}$   whenever

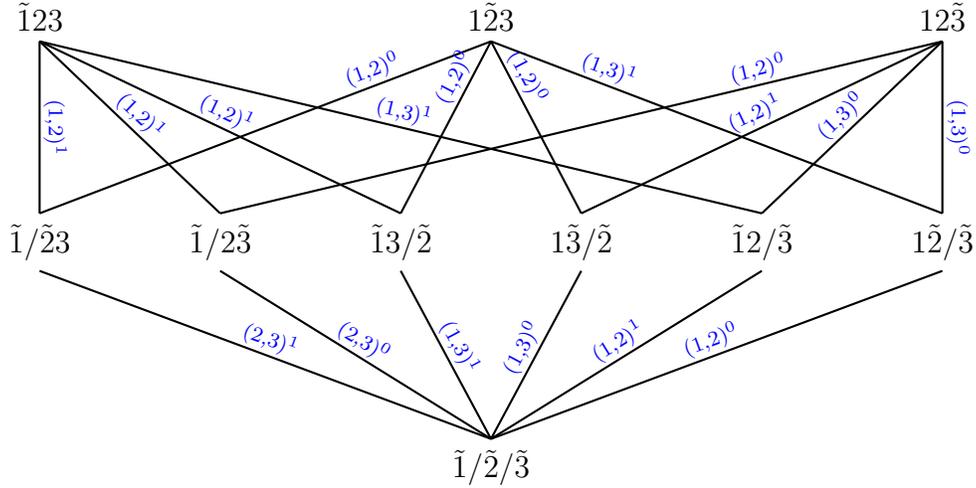
\begin{figure}
    \begin{center}
  \begin{tikzpicture}[scale=1.2]
\node(root){$\tilde{1}/\tilde{2}/\tilde{3}$};
\node at ($(root)+(-5,5)$)(one){$\tilde{1}23$};
\node at ($(root)+(0,5)$)(two){$1\tilde{2}3$};
\node at ($(root)+(5,5)$)(three){$12\tilde{3}$};
\node at ($(root)+(-5,2.5)$)(21){$\tilde{1}/\tilde{2}3$};
\node at ($(root)+(-3,2.5)$)(22){$\tilde{1}/2\tilde{3}$};
\node at ($(root)+(-1,2.5)$)(23){$\tilde{1}3/\tilde{2}$};
\node at ($(root)+(1,2.5)$)(24){$1\tilde{3}/\tilde{2}$};
\node at ($(root)+(3,2.5)$)(25){$\tilde{1}2/\tilde{3}$};
\node at ($(root)+(5,2.5)$)(26){$1\tilde{2}/\tilde{3}$};
\path[thick](root.north)                                            edge node[sloped,yshift=0.5em,draw=none,fill=none,color=blue]{{$\scriptstyle{(2,3)^1}$}}
(21.south)
edge node[sloped,yshift=0.5em,draw=none,fill=none,color=blue]{{$\scriptstyle{(2,3)^0}$}}
(22.south)
edge node[sloped,yshift=0.5em,draw=none,fill=none,color=blue]{{$\scriptstyle{(1,3)^1}$}}
(23.south)edge node[sloped,yshift=0.5em,draw=none,fill=none,color=blue]{{$\scriptstyle{(1,3)^0}$}}(24.south)edge node[sloped,yshift=0.5em,draw=none,fill=none,color=blue]{{$\scriptstyle{(1,2)^1}$}}(25.south)edge node[sloped,yshift=0.5em,draw=none,fill=none,color=blue]{{$\scriptstyle{(1,2)^0}$}}(26.south);
\path[thick](one.south)
edge node[sloped,yshift=0.5em,draw=none,fill=none,color=blue]{{$\scriptstyle{(1,2)^1}$}}
(21.north)
edge node[sloped,yshift=0.5em,draw=none,fill=none,color=blue]{{$\scriptstyle{(1,2)^1}$}}
(22.north)
edge node[sloped,yshift=0.5em,draw=none,fill=none,color=blue]{{$\scriptstyle{(1,2)^1}$}}
(23.north)edge node[sloped,yshift=0.5em,draw=none,fill=none,color=blue]{{$\scriptstyle{(1,3)^1}$}}(25.north);
\path[thick](two.south)edge node[near start,sloped,yshift=0.5em,draw=none,fill=none,color=blue]{{$\scriptstyle{(1,2)^0}$}}(21.north)edge node[near start,sloped,yshift=0.5em,draw=none,fill=none,color=blue]{{$\scriptstyle{(1,2)^0}$}}(23.north)edge node[near start,sloped,yshift=0.5em,draw=none,fill=none,color=blue]{{$\scriptstyle{(1,2)^0}$}}(24.north)edge node[near start,sloped,yshift=0.5em,draw=none,fill=none,color=blue]{{$\scriptstyle{(1,3)^1}$}}(26.north);
\path[thick](three.south)edge node[near start, sloped,yshift=0.5em,draw=none,fill=none,color=blue,color=blue]{{$\scriptstyle{(1,2)^0}$}}(22.north)edge node[sloped,yshift=0.5em,draw=none,fill=none,color=blue]{{$\scriptstyle{(1,2)^1}$}} (24.north)edge node[sloped,yshift=0.5em,draw=none,fill=none,color=blue]{{$\scriptstyle{(1,3)^0}$}}(25.north)edge node[sloped,yshift=0.5em,draw=none,fill=none,color=blue]{{$\scriptstyle{(1,3)^0}$}}(26.north);
 \end{tikzpicture}
 \end{center}
 \caption{$\Pi_{3}^{\bullet}$ with its edge labeling $\lambda_{\bullet}$.}
    \label{fig:pointed_3}
\end{figure}

\begin{itemize}
  \item $\pi\lessdot \pi'$ in $\Pi_n$.
  \item if $B_h=A_i\cup A_j$ then $p_h \in \{q_i,q_j\}$.
  \item if $B_h=A_i$ then $p_h=q_i$.
\end{itemize}
Thus to move up in a cover, exactly two blocks are merged and the pointed element of this new block is one of the pointed elements of the merged blocks.  We will represent the pointed element for each block by placing a tilde above the pointed element. For example, $\{1478\}^4$ will be denoted by $1\tilde{4}78$. The Hasse diagram of $\Pi_3^{\bullet}$ is illustrated in Figure \ref{fig:pointed_3}. As noted in the introduction,  $\Pi_n^\bullet$ is (isomorphic to) the poset of partitions for the operad  $\Perm$.

Suppose we are merging two blocks $A$ and $B$ with $\min A< \min B$. We say that this merge is a  \emph{$0$-merge} if the pointed element of $A\cup B$ is the pointed element of $B$. Similarly, we say the merge is a \emph{$1$-merge} if the pointed element of $A\cup B$ is the pointed element of $A$.  For example, if we merge the blocks $1\tilde{2}4$ with $3\tilde{5}$ to get $1234\tilde{5}$ we have done a $0$-merge.  On other hand, if we had obtained $1\tilde{2}345$, we would have done a $1$-merge.  From time to time, we will need to discuss merges where we do not know whether it is a $0$ or $1$-merge. In these cases, we will refer to it as an \emph{$u$-merge}, always bearing in mind that $u\in \{0,1\}$.

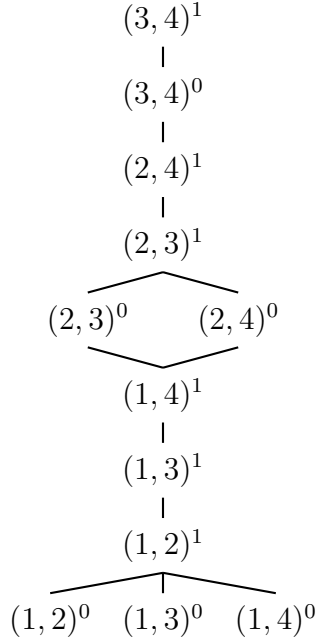
\begin{figure}
    \centering
\begin{align*}
    \begin{tikzpicture}
\node(one){$(1,3)^0$};
\node at ($(one)+(-1.5,0)$)(two){$(1,2)^0$};
\node at ($(one)+(1.5,0)$)(three){$(1,4)^0$};
\node at ($(one)+(0,1)$)(four){$(1,2)^1$};
\node at ($(one)+(0,2)$)(five){$(1,3)^1$};
\node at ($(one)+(0,3)$)(six){$(1,4)^1$};
\node at ($(one)+(-1,4)$)(seven){$(2,3)^0$};
\node at ($(one)+(1,4)$)(eight){$(2,4)^0$};
\node at ($(one)+(0,5)$)(nine){$(2,3)^1$};
\node at ($(one)+(0,6)$)(ten){$(2,4)^1$};
\node at ($(one)+(0,7)$)(eleven){$(3,4)^0$};
\node at ($(one)+(0,8)$)(twelve){$(3,4)^1$};
\path[thick](one.north)edge 
(four.south);
\path[thick](two.north)edge 
(four.south);
\path[thick](three.north)edge 
(four.south);
\path[thick](four.north)edge 
(five.south);
\path[thick](five.north)edge 
(six.south);
\path[thick](six.north)edge 
(seven.south)edge(eight.south);
\path[thick](seven.north)edge 
(nine.south);
\path[thick](eight.north)edge 
(nine.south);
\path[thick](nine.north)edge 
(ten.south);
\path[thick](ten.north)edge 
(eleven.south);
\path[thick](eleven.north)edge 
(twelve.south);
\end{tikzpicture}
\end{align*}
\caption{Poset of labels $\Lambda_4^{\bullet}$.}
    \label{fig:labels_pointed_3}
\end{figure}

We now give an edge labeling of $\Pi_n^\bullet$.  We first define the poset of labels. Let  $\Lambda_n^{\bullet}$  be the set $\{(a,b)^u\mid 1\le a < b\le n \mbox{ and } u\in \{0,1\}\}$. To define the order relation on  $\Lambda_n^{\bullet}$, let $A_a$ be the antichain on the set $A_a=\{(a, b)^0 \mid a<b\le n \}$ and let $C_a$ be the chain on the set  $\{(a,b)^1\mid a <b\le n\}$ where $(a, b)^1< ( a,c )^1$  whenever $b<c$. Then we define $\Lambda_n^{\bullet}$ as the ordinal sum $$\Lambda_n^{\bullet}:=A_1\oplus C_1\oplus A_2\oplus C_2 \oplus \cdots \oplus A_{n-1}\oplus C_{n-1}.$$
The Hasse diagram of $\Lambda_4^{\bullet}$ is given in  Figure \ref{fig:labels_pointed_3}.  Note that the underlying sets of $\Lambda^w_{n}$ and $\Lambda_n^\bullet$ are the same, but their partial orders are different.  Now suppose that in the cover relation $\bpi'\cover \bpi$, we $u$-merge blocks $A$ and $B$.  Then we define the labeling $\lambda_\bullet: \Pi_n^\bullet \rightarrow \Lambda_n^\bullet$ by 
\begin{align}
    \lambda_{\bullet}(\bpi \lessdot \bpi')= (\min A, \min B)^u.
    \label{equation:definition_lambda_pointed}
\end{align}
 In Figure \ref{fig:pointed_3} we illustrate the labeling $\lambda_{\bullet}$ of $\Pi_3^{\bullet}$.

  We now turn our attention to proving that $\lambda_\bullet$ is an EW-labeling.  First, let us note that a label $\lambda_{\bullet}(\bpi \lessdot \bpi')$ completely determines which two blocks of $\bpi$ merge to form a block of $\bpi'$ and which element in the resulting block is pointed. Hence, for every $\bpi\in \Pi_n^{\bullet}$ the cover relations over $\bpi$ have distinct labels. Thus starting at any element $\bpi \in \Pi_n^{\bullet}$, a sequence of valid labels completely determines a saturated chain starting at $\bpi$. Thus we obtain the following proposition.

\begin{proposition} \label{proposition:injectivity_pointed} The labeling $\lambda_{\bullet}$ of equation \eqref{equation:definition_lambda_pointed} is injective on maximal chains in any interval of $\Pi_n^{\bullet}$.
\end{proposition}

Next we show that $\lambda_{\bullet}$ is an ER-labeling. 
For a finite $A\subset \NN$ we denote $\Pi_A$ the poset of partitions of $A$ and $\Pi_A^{\bullet}$ the poset of pointed partitions supported on $A$. We also use $U(x)=\{y\in P \mid y\geq x\}$ to denote the \emph{(principal) upper filter} generated by an element $x$ in a poset $P$.  It turns out that the upper filter of any element of $\Pi_n^\bullet$ is isomorphic to another pointed partition poset and that this isomorphism preserves the labeling $\lambda_\bullet$.  We make this explicit next.

\begin{lemma}\label{lemma:isomorphism_pointed_upper_filters} Let $\balpha=\{B_1^{p_1},\dots, B_l^{p_l}\} \in \Pi_n^{\bullet}$ with $\min B_1< \cdots < \min B_l$. Let $$\Phi:U(\balpha)\rightarrow \Pi_{\{\min B_1, \dots, \min B_l\}}^{\bullet}$$ be the map defined as follows:
\begin{enumerate}
    \item For a pointed set $A^q$  with  $A=B_{j_1}\cup\cdots \cup B_{j_r}$ with $j_1<\dots<j_r$ and $q=p_{j_s}$ for some $s\in [r]$ we define $\Phi(A^q):=\{\min B_{j_1}\cup\dots \cup \min B_{j_r}\}^{\min B_{j_s}}$. 
    \item For any $\bpi \in U(\balpha)$ we define $\Phi(\bpi):=\{\Phi(A^q)\mid A^q \in \bpi\}$.
\end{enumerate}
Then the map $\Phi$ is an isomorphism  preserving the labeling $\lambda_{\bullet}$ defined in equation \eqref{equation:definition_lambda_pointed}, i.e., for any $\bpi \lessdot \bpi'$ in $U(\balpha)$ we have that 
$$\lambda_{\bullet}(\Phi(\bpi) \lessdot \Phi(\bpi'))=\lambda_{\bullet}(\bpi \lessdot \bpi').$$
\end{lemma}

Before we prove the lemma, let us provide a quick example of the map $\Phi$. Suppose that $\alpha = 14\tilde{5}6/ 2\tilde{7}9/3\tilde{8}$.  Then $1256\tilde{7}9/ 3\tilde{8}$ is in $U(\alpha)$ and $\Phi(1256\tilde{7}9/ 3\tilde{8}) = 1\tilde{2}/\tilde{3}$.  The pointed block $1\tilde{2}$ comes from the fact that we merged the blocks $14\tilde{5}6$ and $ 2\tilde{7}9$ and chose to keep $7$ pointed. As a result we point $2$ when we apply the map $\Phi$ since $2$ is the minimum element in the block containing $7$.  The block $3\tilde{8}$ does not get merged, but since we reduce to the minimum element of the block when applying $\Phi$, we get the pointed block $\tilde{3}$.

\begin{proof}
We will show first that the function $\Phi$ preserves the $u$-merging of two blocks $u\in \{0,1\}$. Let $A_1=B_{j_1}\cup\cdots \cup B_{j_r}$ with $j_1<j_2<\cdots <j_r$ and $q_1=p_{j_l}$ for some $l\in [r]$ and let $A_2=B_{k_1}\cup\cdots \cup B_{k_t}$ with $k_1<k_2<\cdots <k_t$ and $q_2=p_{k_m}$ for some $m\in [t]$. Without loss of generality we assume $j_1<k_1$ so $\min A_1 < \min A_2$. 
We denote $A_1^{q_1}\cup_{u} A_2^{q_2}=(A_1\cup A_2)^{q}$ the $u$-merging of the pointed blocks $A_1^{q_1}$ and $A_2^{q_2}$ where $q=q_1$ if $u=1$ and $q=q_2$ if $u=0$. 
\begin{align*}
     \Phi(A_1^{q_1}\cup_u A_2^{q_2})&=\Phi(\left\{ B_{j_1}\cup\cdots \cup  B_{j_r}\cup  B_{k_1}\cup\cdots \cup  B_{k_t}\right\}^{q})\\
     &=\left\{\min B_{j_1}\cup\cdots \cup \min B_{j_r}\cup \min B_{k_1}\cup\cdots \cup \min B_{k_t}\right\}^{\tilde q}\\
     &=\left\{\min B_{j_1}\cup\cdots \cup \min B_{j_r}\right\}^{\min B_{j_s}}\cup_u \left\{\min B_{k_1}\cup\cdots \cup \min B_{k_t}\right\}^{\min B_{k_u}}\\
     &=\Phi(A_1^{q_1})\cup_u \Phi(A_2^{q_2}), 
 \end{align*}
 where $\tilde q=\min B_{j_s}$ if $u=1$ and $\tilde q=\min B_{k_u}$ if $u=0$. Since the blocks of $\balpha$ are in bijection with the blocks of $\min B_1/\cdots/\min B_l$ and all elements of $U(\balpha)$ are obtained uniquely by a sequence of $u$-merges of blocks of $\balpha$ and the elements of $\Pi_{\{\min B_1, \dots, \min B_l\}}^{\bullet}$ are obtained uniquely by a sequence of $u$-merges of the blocks of $\min B_1/\cdots/\min B_l$, we conclude that $\Phi$ is a bijection.  Moreover, $\Phi$ and $\Phi^{-1}$ preserve cover relations and hence $\Phi$ is a poset isomorphism.
 
 Now, to see that the labeling  according to $\lambda_{\bullet}$ of equation \eqref{equation:definition_lambda_pointed} is preserved, note that in a cover relation where we $u$-merge the blocks $A_1^{q_1}$ and $A_2^{q_2}$ the label is
 $$(\min A_1, \min A_2)^u=(\min B_{j_1},\min B_{k_1})^u,$$
 which is the same obtained by $u$-merging the blocks $\Phi(A_1^{q_1})$ and $\Phi(A_2^{q_2})$.
\end{proof}

As we explain in the proof of the following proposition, Lemma \ref{lemma:isomorphism_pointed_upper_filters} essentially reduces the task of finding a unique increasing chain in each interval to finding an increasing chain in every maximal interval.

\begin{proposition}\label{proposition:ER_pointed}
 The labeling $\lambda_{\bullet}$ of equation \eqref{equation:definition_lambda_pointed} is an ER-labeling of $\Pi_{n}^{\bullet}$.
\end{proposition}

\begin{proof}
Let $\bpi,\bpi' \in \Pi_n^{\bullet}$ such that $\bpi \leq \bpi'$. We want to show that there is a unique increasing saturated chain in $[\bpi,\bpi']$. 

Assume first that $\bpi=\hat{0}$ and $\bpi'=[n]^p$, so $[\bpi,\bpi']=[\hat{0},[n]^p]$ is a maximal interval. We will construct an increasing saturated chain in $[\hat{0},[n]^p]$ and show that such chain is the only increasing saturated chain in  $[\hat{0},[n]^p]$.  Consider the chain $c_{[n]^p}$ whose label sequence is as follows.

\begin{equation}\label{equation:increasing_chain_pointed}
    \lambda_{\bullet}(c_{[n]^p})=
    \begin{cases}
     (1,2)^1 (1,3)^1\cdots (1,n-1)^1 (1,n)^1 & \text{ if } p =1,\\ 
     (1,p)^0(1,2)^1\cdots (1,p-1)^1(1,p+1)^1 \cdots (1,n)^1  & \text{ if } 1<p<n,\\
     (1,n)^0(1,2)^1\cdots (1,n-1)^1& \text{ if } p=n.
    \end{cases}
 \end{equation}

Because of Proposition \ref{proposition:injectivity_pointed}, there is at most one such chain with the above label sequence.  It is not hard to check that such a chain does in fact exist. In the case $p=1$, it is easy to see that the chain is increasing.  On the other hand, if $p\neq 1$, the chain is also increasing since $(1,p)^0$ is smaller than any  label of the form $(1,b)^1$ and the remaining values are increasing in $\Lambda_n^\bullet$. 

We now show that the chain $c_{[n]^p}$ is indeed the only increasing chain in $[\hat{0},[n]^p]$. We discuss the case when $p\neq 1$. The case when $p=1$ follows the same idea. Note that if $c'$ is any other chain in $[\hat{0},[n]^p]$ it must have as final label either $(1,a)^0$ or $(1,a)^1$ for some $a\neq 1$ since in the last step the block with minimal label $1$ always be involved.  It follows that for $c'$ to be increasing all the labels along the chain must be of the form $(1,b)^u$ for some $b$ and $u$.  Hence $c'$ has to be constructed by a step-by-step process of merging blocks with the block that contains the element $1$. Hence, the labels in the second component will form a permutation  of the elements $\{2,3,\dots, n\}$. Since $p$ has to be the pointed element, we will have a step where the label $(1,p)^0$ appears. Since  $(1,p)^0$ and $(1,a)^0$ are not comparable when $a\neq p$, we see that $c'$ cannot have the label $(1,a)^0$ where $a\neq p$ as $c'$ would not be increasing.  Hence, all other labels are of the form  $(1,a)^1$ and the only way to order them increasingly is as in equation \eqref{equation:increasing_chain_pointed}. By Proposition \ref{proposition:injectivity_pointed}, $\lambda_{\bullet}$ is injective and so $c'=c_{[n]^p}$.

Now, we consider an interval of the form $[\hat{0},\bpi]$ where $\bpi \in \Pi_n^{\bullet}$ and $\bpi$ has at least two blocks. Let $\bpi=\{B_1^{p_1},\dots,B_l^{p_l}\}$ where $\min B_1 <\cdots <\min B_l$. For each $i=1,\dots,l$,  let  $c_{B_i^{p_i}}$ be the unique increasing chain of $[\hat{0},B_i^{p_i}]$. To see why such  chains exist and are unique, apply the same idea from the previous paragraph to each of the intervals  $[\hat{0},B_i^{p_i}]$.  We will now consider the word of labels of  $c_{B_i^{p_i}}$,  $\lambda_{\bullet}(c_{B_i^{p_i}})$.  Note that this word will be empty if $|B_i|=1$. Now let $c_{\bpi}$ be the chain in $[\hat{0},\bpi]$ that first merges the elements with labels in $B_1$ as instructed in $c_{B_1^{p_1}}$, then merges the elements with labels in $B_2$ as instructed in $c_{B_2^{p_2}}$, and so on.  Then  $c_{\bpi}$ has the  word of labels obtained by the concatenation of words
$$\lambda_{\bullet}(c_{\bpi})=\lambda_{\bullet}(c_{B_1^{p_1}})\lambda_{\bullet}(c_{B_2^{p_2}})\cdots\lambda_{\bullet}(c_{B_l^{p_l}}).
$$
Note that this chain is increasing because $\min B_1 <\cdots <\min B_l$ and there is only one chain with this word of labels because of 
Proposition \ref{proposition:injectivity_pointed}.  In order to see that $\lambda_{\bullet}(c_{\bpi})$ is the unique increasing chain in $[\hat{0},\bpi]$, let $c'$ be any other increasing chain in this interval and for every $i=1,\dots,l$, and let 
$$w_i=\lambda_{\bullet}(c')_{i_1}\lambda_{\bullet}(c')_{i_2}\cdots \lambda_{\bullet}(c')_{i_{|B_i|}}$$ be the subword 
of $\lambda_{\bullet}(c')$ whose labels belong to the steps in $c'$
where blocks with elements in $B_i$ were merged. Since $w_i$ is a subword of an increasing word it must also be increasing. Then by the discussion in the paragraph above, we conclude that there is a unique way to apply the merges in order to get an increasing word and this word is $\lambda_{\bullet}(c_{B_i^{p_i}})$. Note that all the labels from all these words are comparable among each other since the $\min B_i$ are all  different. There is then a unique shuffle of the subwords  $\lambda_{\bullet}(c_{B_i^{p_i}})$ that leads to an increasing word  $\lambda_{\bullet}(c')$ which is $\lambda_{\bullet}(c_{B_1^{p_1}})\lambda_{\bullet}(c_{B_2^{p_2}})\cdots\lambda_{\bullet}(c_{B_l^{p_l}})$. So we have that $c'=c_{\bpi}$. 

Finally, consider an interval of the form $[\bpi, \bpi']$ in $\Pi_n^{\bullet}$ with $\bpi=\{B_1^{p_1},\dots,B_l^{p_l}\}$. We have by Lemma \ref{lemma:isomorphism_pointed_upper_filters} that $[\bpi, \bpi']$ is isomorphic to an interval $[\hat{0},\bpi'']$ in the poset $\Pi_{\{\min B_1, \dots, \min B_l\}}^{\bullet}$ through an isomorphism that preserves the labels of the maximal chains. Hence by the discussion in the paragraph before we have that there is a unique increasing chain in the interval  $[\hat{0},\bpi'']$ of the latter poset and hence in $[\bpi, \bpi']$, completing the proof.
\end{proof}

\begin{figure}
    \centering
    \begin{tikzpicture}[scale=1.33]
\node at (-8.65,2) {\textbf{Type I intervals}};
\node at (-8.5,-1.5) {\textbf{Type II intervals}};
\begin{scope}[scale=0.7,xshift=-60]
\tikzstyle{every node}=[scale=0.8]
\node(one){$A^{p_A}/B^{p_B}/C^{p_C}/D^{p_D}$};
\node at ($(one)+(-2,2)$)(two){$AB^{p_1}/C^{p_C}/D^{p_D}$};
\node at ($(one)+(2,2)$)(five){$A^{p_A}/B^{p_B}/CD^{p_2}$};
\node at ($(one)+(0,4)$)(six){$AB^{p_1}/CD^{p_2}$};
\path[thick](one.north)edge[red] node[sloped,yshift=0.5em,draw=none,fill=none]{$\Large{(a,b)^{u_1}}$}
(two.south)edge[red] node[sloped,yshift=0.5em,draw=none,fill=none]{$\Large{(c,d)^{u_2}}$}(five.south);
\path[thick](two.north)edge[red] node[sloped,yshift=0.5em,draw=none,fill=none]{$\Large{(c,d)^{u_2}}$}(six.south);
\path[thick](five.north)edge[red] node[sloped,yshift=0.5em,draw=none,fill=none]{$\Large{(a,b)^{u_1}}$}(six.south);

\end{scope}

\begin{scope}[scale=0.7,xshift=-60,yshift=-150]
\tikzstyle{every node}=[scale=0.8]

\node(one){$A^{p_A}/B^{p_B}/C^{p_C}$};
\node at ($(one)+(-4.5,2)$)(two){$A^{p_A}/BC^{p_B}$};
\node at ($(one)+(-1.5,2)$)(three){$AB^{p_A}/C^{p_C}$};
\node at ($(one)+(1.5,2)$)(four){$AC^{p_A}/B^{p_B}$};
\node at ($(one)+(4.5,2)$)(five){$A^{p_A}/BC^{p_C}$};
\node at ($(one)+(0,4)$)(six){$ABC^{p_A}$};

\path[thick](one.north)edge node[sloped,yshift=0.5em,draw=none,fill=none,color=blue]{$\Large{(b,c)^1}$}
(two.south)edge[red] node[sloped,yshift=0.5em,draw=none,fill=none]{$\Large{(a,b)^1}$}(three.south)edge[red] node[sloped,yshift=0.5em,draw=none,fill=none]{$\Large{(a,c)^1}$}
(four.south)edge node[sloped,yshift=0.5em,draw=none,fill=none,color=blue]{$\Large{(b,c)^{0}}$}(five.south);
\path[thick](two.north)edge node[sloped,yshift=0.5em,draw=none,fill=none,color=blue]{$\Large{(a,b)^1}$}(six.south);
\path[thick](three.north)edge[red] node[near start,sloped,yshift=0.5em,draw=none,fill=none]{$\Large{(a,c)^1}$}(six.south);
\path[thick](four.north)edge[red] node[sloped,yshift=0.5em,draw=none,fill=none]{$\Large{(a,b)^1}$}(six.south);
\path[thick](five.north)edge node[sloped,yshift=0.5em,draw=none,fill=none,color=blue]{$\Large{(a,b)^{1}}$}(six.south);
\end{scope}
\begin{scope}
[scale=0.7,xshift=-250,yshift=-250]
\tikzstyle{every node}=[scale=0.8]

\node(one){$A^{p_A}/B^{p_B}/C^{p_C}$};
\node at ($(one)+(-4.5,2)$)(two){$A^{p_A}/BC^{p_B}$};
\node at ($(one)+(-1.5,2)$)(three){$AB^{p_B}/C^{p_C}$};
\node at ($(one)+(1.5,2)$)(four){$AC^{p_A}/B^{p_B}$};
\node at ($(one)+(4.5,2)$)(five){$AC^{p_C}/B^{p_B}$};
\node at ($(one)+(0,4)$)(six){$ABC^{p_B}$};
\path[thick](one.north)edge node[sloped,yshift=0.5em,draw=none,fill=none,color=blue]{$\Large{(b,c)^1}$}
(two.south)edge[red] node[sloped,yshift=0.5em,draw=none,fill=none]{$\Large{(a,b)^{0}}$}(three.south)edge[red] node[sloped,yshift=0.5em,draw=none,fill=none]{$\Large{(a,c)^1}$}
(four.south)edge node[sloped,yshift=0.5em,draw=none,fill=none,color=blue]{$\Large{(a,c)^{0}}$}(five.south);
\path[thick](two.north)edge node[sloped,yshift=0.5em,draw=none,fill=none,color=blue]{$\Large{(a,b)^{0}}$}(six.south);
\path[thick](three.north)edge[red] node[near start,sloped,yshift=0.5em,draw=none,fill=none]{$\Large{(a,c)^1}$}(six.south);
\path[thick](four.north)edge[red] node[sloped,yshift=0.5em,draw=none,fill=none]{$\Large{(a,b)^{0}}$}(six.south);
\path[thick](five.north)edge node[sloped,yshift=0.5em,draw=none,fill=none,color=blue]{$\Large{(a,b)^{0}}$}(six.south);
\end{scope}

\begin{scope}[scale=0.7,xshift=-60,yshift=-350]
\tikzstyle{every node}=[scale=0.8]

\node(one){$A^{p_A}/B^{p_B}/C^{p_C}$};
\node at ($(one)+(-4.5,2)$)(two){$A^{p_A}/BC^{p_C}$};
\node at ($(one)+(-1.5,2)$)(three){$AC^{p_C}/B^{p_B}$};
\node at ($(one)+(1.5,2)$)(four){$AB^{p_A}/C^{p_C}$};
\node at ($(one)+(4.5,2)$)(five){$AB^{p_B}/C^{p_C}$};
\node at ($(one)+(0,4)$)(six){$ABC^{p_C}$};
\path[thick](one.north)edge node[sloped,yshift=0.5em,draw=none,fill=none,color=blue]{$\Large{(b,c)^0}$}
(two.south)edge[red] node[sloped,yshift=0.5em,draw=none,fill=none]{$\Large{(a,c)^{0}}$}(three.south)edge[red] node[sloped,yshift=0.5em,draw=none,fill=none]{$\Large{(a,b)^1}$}
(four.south)edge node[sloped,yshift=0.5em,draw=none,fill=none,color=blue]{$\Large{(a,b)^{0}}$}(five.south);
\path[thick](two.north)edge node[sloped,yshift=0.5em,draw=none,fill=none,color=blue]{$\Large{(a,b)^{0}}$}(six.south);
\path[thick](three.north)edge[red] node[near start,sloped,yshift=0.5em,draw=none,fill=none]{$\Large{(a,b)^1}$}(six.south);
\path[thick](four.north)edge[red] node[sloped,yshift=0.5em,draw=none,fill=none]{$\Large{(a,c)^{0}}$}(six.south);
\path[thick](five.north)edge node[sloped,yshift=0.5em,draw=none,fill=none,color=blue]{$\Large{(a,c)^{0}}$}(six.south);
\end{scope}

\end{tikzpicture}

\caption{Rank two intervals in $\Pi_3^{\bullet}$. Here $A,B,C,D$ are the blocks that get merged in the interval and $a=\min(A)< b=\min(B)<c =\min(C)<d=\min(d)$.  The  blocks that are not changed in the interval are not depicted.  Edges given in red correspond to the rank two switching property.}

    \label{fig:rank2_pointed}
\end{figure}
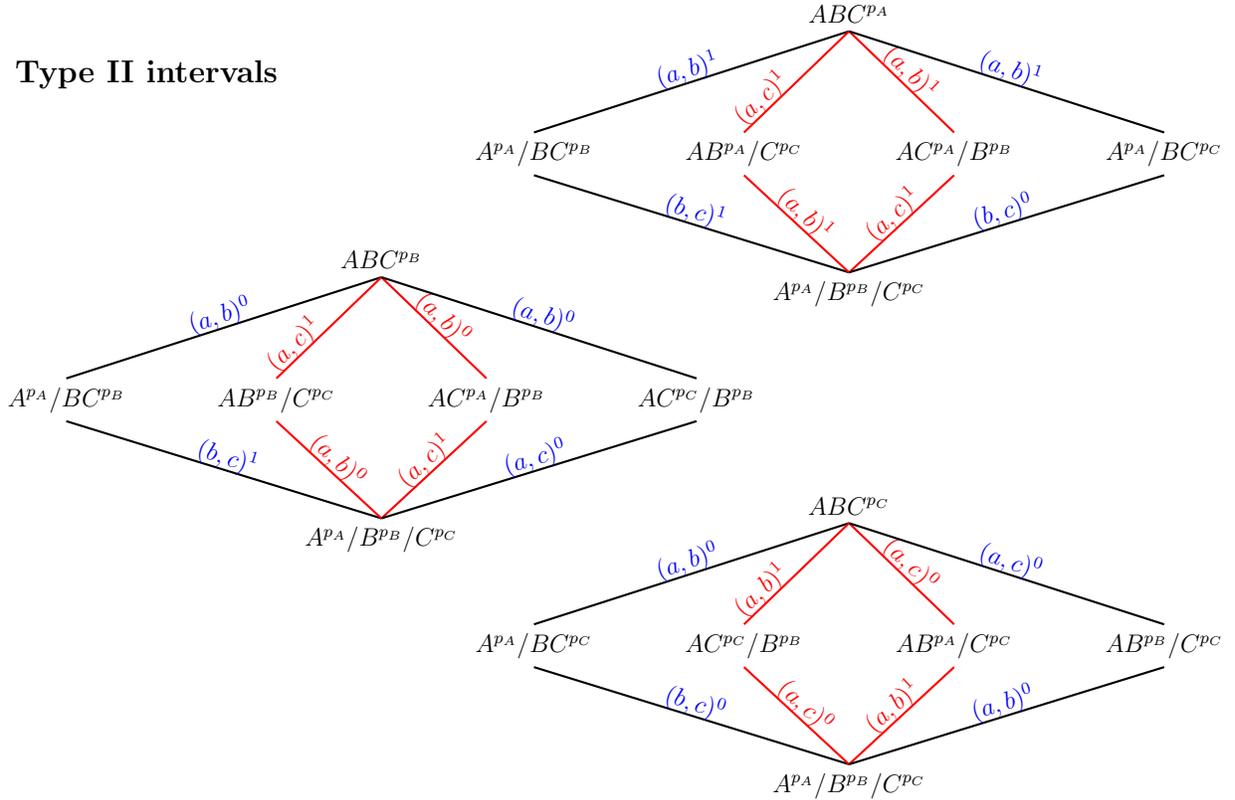

To finish showing that the labeling $\lambda_{\bullet}$ is an EW-labeling, we just need to show $\lambda_\bullet$ has the rank two switching property. There are two types of rank two intervals in $\Pi_n^\bullet$.   The type I is where there are two pairs of blocks that get merged independently of each other and the type II is where there are three blocks all of which get merged. These two types are shown in Figure~\ref{fig:rank2_pointed}. In type II, there are 3 possible choices for which initial block the pointed element at the top of the interval comes from.  From Figure~\ref{fig:rank2_pointed} one can readily see that we have the following proposition. 

\begin{proposition}\label{proposition:rank_two-pointed}
The labeling $\lambda_{\bullet}$ of equation \eqref{equation:definition_lambda_pointed} satisfies the rank two switching property.
\end{proposition}

By Propositions \ref{proposition:injectivity_pointed}, \ref{proposition:ER_pointed} and \ref{proposition:rank_two-pointed} we have that the labeling $\lambda_{\bullet}$ on $\Pi_n^{\bullet}$ satisfies Definition \ref{definition:EW-labeling} which proves the following theorem.

\begin{theorem}\label{theorem:EW-pointed}
The labeling $\lambda_{\bullet}$ is an EW-labeling of $\Pi_{n}^{\bullet}$. As a consequence, the poset $\Pi_{n}^{\bullet}$ has a Whitney dual.
\end{theorem}

We will give a combinatorial description of the corresponding Whitney dual in Section~\ref{subsection:whitney_dual_pointed}.

\subsection{EL-labelings for the pointed partition poset}\label{section:El-labeling_pointed_partition}

In~\cite[Theorem 1.11]{ChapotonVallete2006} Chapoton and Vallete show that $\Pi_n^\bullet$ has a CL-labeling and hence is Cohen-Macaulay.   In Remark 1.11 of that paper, they  leave open the question of if  maximal intervals of $\Pi_n^\bullet$  have  EL-labelings.    We give a positive answer to their question below by providing an EL-labeling for $\Pi_n^\bullet$ (which restricts to an EL-labeling in every maximal interval). Before we do, let us note that Bellier-Mill{\`e}s, Delcroix-Oger and Hoffbeck~\cite[Proposition 3.13]{BelliermillesDelcroixogerHoffbeck2021} propose an edge labeling for $\Pi_n^{\bullet}$ that the authors claim is an EL-labeling. However, we later argue that the proposed labeling does not satisfy the conditions to be an EL-labeling (see Remark \ref{remark:proposed_EL_labeling}). 

One might hope that our previous EW-labeling $\lambda_\bullet$ is an EL-labeling. Unfortunately, this is not the case. Indeed,  consider the rank two interval $[ \tilde{1}/\tilde{2}/\tilde{3},12\tilde{3}]$.  The unique increasing chain $\tilde{1}/\tilde{2}/\tilde{3}\cover 1\tilde{3}/\tilde{2}\cover 12\tilde{3}$ has word of labels $(1,3)^0(1,2)^1$.  However, this sequence is not lexicographically comparable with the word of labels $(1,2)^0(1,3)^0$ of the chain $\tilde{1}/\tilde{2}/\tilde{3}\cover 1\tilde{2}/\tilde{3}\cover 12\tilde{3}$ in the same interval (see Figure~\ref{fig:pointed_3}).  Although $\lambda_\bullet$  is not an EL-labeling, if we keep the same edge labels, but instead use the ordering of the labels that we used for the weighted partition poset, we do get an EL-labeling.  More specifically, we claim that the labeling  $\lambda_{\bullet_2}: \mathcal{E}(\Pi_n^\bullet) \rightarrow \Lambda_n^w$  where $\lambda_{\bullet_2}(\bpi \cover \bpi') := \lambda_{\bullet}(\bpi \cover \bpi')$ is an EL-labeling.

The following theorem has a very similar proof to the one in Proposition \ref{proposition:ER_pointed} and \cite[Theorem 3.2]{DleonWachs2016}. To avoid a lengthy discussion we just provide the relevant steps in the proof, which can be verified by the reader.
\begin{theorem}\label{theorem:EL-labeling_pointed}
 The labeling $\lambda_{\bullet_2}$ is an EL-labeling of $\Pi_n^{\bullet}$.  Consequently, $\Pi_n^{\bullet}$ is EL-shellable and its maximal intervals are Cohen-Macaulay.
\end{theorem}
\begin{proof}[Proof idea:] We need to show that in each interval $[\bpi,\bpi']$ of $\Pi_n^{\bullet}$ there is a unique increasing maximal chain and that this chain is lexicographically first.

First we consider an interval of the form $[\hat{0},[n]^p]$. For this type of interval the reader can verify that there is an  increasing maximal chain  that has word of labels

\begin{equation}\label{equation:increasing_chain_pointed_2}
    \lambda_{\bullet_2}(c' _{[n]^p})=
    \begin{cases}
      (1,2)^{1}\cdots(1,n)^1 & \text{ if } p =1,\\ (1,2)^{0}\cdots(1,p)^0(1,p+1)^1\cdots(1,n)^1 & \text{ if } 1<p<n,\\
       (1,2)^{0}\cdots(1,n)^0, & \text{ if } p=n,
    \end{cases}
 \end{equation}   
 and which is of the form (where we represent each pointed partition by its unique nonsingleton block):
\begin{align*}
    c'_{[n]^p}=(\hat{0}\lessdot [2]^2\lessdot \cdots \lessdot [p]^{p} \lessdot [p+1]^{p} \lessdot \cdots \lessdot [n]^{p}).
\end{align*}
We note that a similar argument given in the proof of Proposition~\ref{proposition:ER_pointed} shows that this is the only increasing maximal chain in  $[\hat{0},[n]^p]$. To show that this chain is lexicographically smallest, suppose this was not the case.   Then there is some other maximal chain, $d$, whose words of labels  is not lexicographically larger.  Let $d_1d_2\cdots d_{n-1}$ be the word of labels of $d$ and assume that the first time it disagrees with the increasing chain at $d_{i-1}$. Note that we may assume that $i>2$ since the first label along the increasing chain is the smallest possible label.

First, suppose that $i-1<p$ and that  $(1,i)^0 \not < d_{i-1}$.  Then based on $\Lambda_n^w$, $d_{i-1}$ must be of the form $(1,b)^1$ where $b<i$.  But this is impossible since by the time $d$ adds the label $d_{i-1}$, $b$ was already in the same block as $1$ as $d$ agrees with $c'_{[n]^p}$ up to this step.  If $i-1\geq p$ and $(1,i)^1 \not < d_{i-1}$. This would imply that $d_{i-1}$ is of one of the following forms: $(1,b)^0$ with $b>i$ or $(1,b)^u$ with $b< i$ and $u\in \{0,1\}$.   By this point along $d$, $p$ is the pointed element in its block and this block contains $1$.  So, all the labels at this point must have an exponent of $1$ and thus it cannot be of the form $(1,b)^0$.  It also cannot be  of the form $(1,b)^1$ with $b<i$ since $b$ is already in the block with $1$ by this point along $d$. We conclude that the unique maximal chain is lexicographically smallest. 

For an interval of the form $[\hat{0},\bpi]$ where $\bpi$ is of the form $\bpi=\{B_1^{p_1},\dots,B_l^{p_l}\}$ with $\min B_1 <\cdots <\min B_l$ and $l\ge 2$, we consider the unique increasing word of labels $c'_{B_i^{p_i}}$ in $[\hat{0},B_i^{p_i}]$  and then the unique maximal chain  $c'_{\bpi}$ in  $[\hat{0},\bpi]$ with word of labels
$$\lambda_{\bullet_2}(c'_{\bpi})=\lambda_{\bullet}(c'_{B_1^{p_1}})\lambda_{\bullet}(c'_{B_2^{p_2}})\cdots\lambda_{\bullet}(c'_{B_l^{p_l}})$$
is the unique increasing chain and is lexicographically first among maximal chains  in $[\hat{0},\bpi]$.

Finally, for an interval of the form $[\bpi, \bpi']$ in $\Pi_n^{\bullet}$, we use Lemma \ref{lemma:isomorphism_pointed_upper_filters} to reduce to any of the two cases before. Note that the lemma still applies in this case since the functions $\lambda_{\bullet_2}$ and  $\lambda_{\bullet}$    only differ in the  order structure on the poset of labels.
\end{proof}

At this point, the reader may be wondering if  $\lambda_{\bullet_2}$ is an EW-labeling. By looking at the last occurrence of rank $2$ intervals of type II in Figure \ref{fig:rank2_pointed}, we can conclude that the unique increasing chain has a word of labels $(a,b)^0(a,c)^0$, but there is no chain with a word of labels $(a,c)^0(a,b)^0$. Hence, the EL-labeling $\lambda_{\bullet_2}$ of Theorem \ref{theorem:EL-labeling_pointed} fails the rank two switching property and is not an EW-labeling.  

As we mentioned earlier, $\lambda_\bullet$ is an EW-labeling, but is not an EL-labeling. Nevertheless, if we take the order dual of $\Pi_n^\bullet$ and reverse the ordering on the labels for $\lambda_{\bullet}$, we do get an EL-labeling of the order dual.

Given a poset $P$, let $P^*$ be the order dual of $P$.  Moreover, given a labeling $\lambda:\mathcal{E}(P)\rightarrow \Lambda$ of a poset $P$ with label poset $\Lambda$, we define the \emph{dual labeling} $\lambda^*:\mathcal{E}(P^*)\rightarrow \Lambda^*$ of the order dual poset $P^*$ to be given by
$$
\lambda^*(y\cover_{P^*} x)  = \lambda(x\cover_{P} y).
$$
In other words, the edge labels do not change when passing from $P$ to its order dual $P^*$, just the ordering on the labels.

\begin{theorem}\label{theorem:pointed_dual_EL_labeling}
The labeling  $\lambda^*_{\bullet}$  is an EL-labeling  of ${\Pi_n^{\bullet}}^*$.  Consequently the maximal intervals of the order dual are EL-shellable.
\end{theorem}
\begin{proof}

First note that since we reverse the order of the labels from  $\lambda_{\bullet}$  to get $\lambda^*_{\bullet}$,  an  increasing chain in an interval $[\balpha, \bpi]$ of ${\Pi_n^{\bullet}}^*$ is exactly the order dual of an  increasing chain in the interval $[\bpi, \balpha]$ of  ${\Pi_n^{\bullet}}$.  It follows that since  $\lambda_\bullet$ is an ER-labeling,   $\lambda^*_{\bullet}$ is also an ER-labeling.  So to finish the proof, we need only show that in every interval of ${\Pi_n^{\bullet}}^*$, the increasing chain with respect to $\lambda^*_\bullet$ is lexicographically smallest. Note that when we restrict the (unique) increasing chain on an interval to a smaller subinterval, that restriction is again the (unique) increasing chain in the said subinterval. Now, since the order on the labels are reversed, it is enough to show that in any interval  of ${\Pi_n^{\bullet}}$ the last label along the increasing chain is strictly larger than the other possible last labels of other chains in that interval. The rest of the argument will follow by induction on the smaller subinterval that is obtained by removing the last step on the unique increasing chain.  By appealing to Lemma~\ref{lemma:isomorphism_pointed_upper_filters}, we need only need to check this   condition for increasing chains   in intervals of the form $[\hat{0},\bpi]$ in ${\Pi_n^{\bullet}}$. This is what we do next.

Consider the interval $[\hat{0}, \bpi]$ and let $c_{\bpi}$  be the increasing chain.  Suppose that the last cover relation on $c_{\bpi}$  is $\balpha\cover \bpi$.  We will show that the label $\lambda_\bullet(\balpha\cover \bpi)$ is strictly larger than any other label of the form $\lambda_\bullet(\balpha'\cover \bpi)$.    Suppose that $\bpi$   is of the form 
$$
\bpi=\{B_1^{p_1},\dots,B_s^{p_s},B_{s+1}^{p_{s+1}},\dots,B_l^{p_l}\}
$$
and where $B_1^{p_1},\dots,B_s^{p_s}$ are the nonsingleton blocks and $\min B_1 <\cdots <\min B_s$.

Note that for any $\balpha'\lessdot \bpi$, we have labels of the form $\lambda_{\bullet}(\balpha'\lessdot \bpi)=(\min B_i,a)^u$ with $a\in B_i$, $i=1,\dots,s$ and $u\in\{0,1\}$. Hence,  the largest possible label that appears along the edges of $[\hat{0},\bpi]$ is of the form $(\min B_{s},a)^u$.  Moreover, since $\min B_s \geq \min B_i$ for all $1\leq i\leq s$, we know that $(\min B_{s},a)^u$ is larger than any label in $[\hat{0}, \bpi]$ of the form $(\min B_i, b)^u$ with $i\neq s$.   Since the elements of $B_s$ must be merged together when going from $\hat{0}$ to $\bpi$,  a label of this form must occur on every maximal chain in $[\hat{0},\bpi]$.  It follows that the increasing chain must end with a label of the form $(\min B_{s},a)^u$.  

Now we distinguish between two cases whether $|B_{s}|=2$ or $|B_{s}|>2$. If $|B_{s}|=2$ there is only one label of the form $(\min B_{s},a)^u$ where $a$ is the unique element in $B_{s}\setminus \{\min B_{s}\}$ and $u\in\{0,1\}$ is uniquely determined to be $u=0$ if $p_s=a$ or $u=1$ otherwise. In either case, all other labels in $[\hat{0},\bpi]$ are of the form $(\min B_{i},b)^u$, which are strictly smaller than $(\min B_{s},a)^u$.

Now suppose that $|B_{s}|>2$ and let $b=\max B_s$. If $p_s\neq b$ then the label in the uppermost cover relation of $c_{\bpi}$ is $(\min B_{s},b)^1$ (see the proof of Proposition \ref{proposition:ER_pointed}) which is the largest label among all labels that appear in the interval $[\hat{0},\bpi]$.  Since $b$ is the largest value in its block, there is only one $\balpha'$ with this property. That is $\balpha'=\{B_1^{p_1},\dots,B_s\setminus\{b\}^{p_s},\{b\}^{b},B_{s+1}^{p_{s+1}},\dots,B_l^{p_l}\}$ which is the second to last element of $c_{\bpi}$. 

If $p_s= b$ then  the label in the uppermost cover relation of $c_{\bpi}$ is $(\min B_{s},c)^1$ where $c=\max (B_s \setminus \{b\})$. Note that the only label of the form $(\min B_{s},a)^u$ larger than   $(\min B_{s},c)^1$ is $(\min B_{s},b)^1$.   
But this label cannot actually appear among the cover relations  $\balpha'\lessdot \bpi$ since this would indicate that $p_s=b$ which is not the case. Thus $(\min B_{s},c)^1$ is the largest possible label that can occur  in $[\hat{0}, \bpi]$  and this also determines uniquely $\balpha'=\{B_1^{p_1},\dots,B_s\setminus\{c\}^{b},\{c\}^{c},B_{s+1}^{p_{s+1}},\dots,B_l^{p_l}\}$ which is the second to last element of $c_{\bpi}$. We conclude that $\lambda^*_\bullet$ is an EL-labeling.
\end{proof}

We finish this section with a remark regarding a proposed  edge labeling of $\Pi_n^\bullet$ given by Bellier-Mill{\`e}s, Delcroix-Oger, and Hoffbeck in~\cite {BelliermillesDelcroixogerHoffbeck2021} \footnote{In the ArXiv v2 version of this paper, the authors have not included this labeling anymore. We include here an explanation of why the proposed labeling is not an EL-labeling.}.

\begin{remark}\label{remark:proposed_EL_labeling} The authors in \cite{BelliermillesDelcroixogerHoffbeck2021} define a labeling $\tilde \lambda:\mathcal{E}(\Pi_{n}^{\bullet})\rightarrow \NN\times \NN$ where  $\NN\times \NN$ has the lexicographic order. To describe their labeling, let $\bpi \lessdot \bpi'$ be such that the two pointed blocks in $\bpi$ which were $u$-merged to get $\bpi'$ are $A_i^{q_i}$ and $A_j^{q_j}$, with $a=\min(A_i)< b=\min(A_j)$. Then define 
$$\tilde \lambda(\bpi \lessdot \bpi')=  \begin{cases}
    (b,a+n-|\bpi|)       & \quad \text{if }u=0,\\
    (b,b+n-|\bpi|)  & \quad \text{if } u=1.
  \end{cases}$$
 
     With this labeling in the interval $[\hat{0},12\tilde{3}]$ of $\Pi_3^{\bullet}$, we see that the chains $\hat{0}\lessdot  \tilde{1}2/\tilde {3}< 12\tilde{3}$ and  $\hat{0}\lessdot 1\tilde{2}/ \tilde {3}<12\tilde{3}$ have words of labels $(2,2)(3,2)$   and $(2,1)(3,2)$ respectively, which are both increasing in the lexicographic order on $\NN\times \NN$. This shows that $\tilde \lambda$ already fails to satisfy the requirement of the uniqueness of the increasing chain in the interval $[\hat{0},[3]^3]$.   Note that this issue does not arise in the interval $[\hat{0},[3]^1]$. So $\tilde \lambda$ is an EL-labeling of at least one maximal interval of $\Pi^{\bullet}_{3}$. Moreover since  $[\hat{0},[3]^1] \cong  [\hat{0},[3]^3]$,  their labeling shows all maximal intervals of $\Pi_{3}^\bullet$   have an EL-labeling.
     However, by extending this idea one can show that if $n\geq 6$, $\tilde \lambda$ is not an EL-labeling for any maximal interval of $\Pi_n^{\bullet}$.  To see why, note that if $n\geq 6$, $[\hat{0}, [n]^i]$ always contains an interval of the form $[\hat{0}, [6]^j/\tilde{7}/\cdots/ \tilde{n}]$ where $j\in [6]$.  If $j=4,5,6$ then the interval $[\hat{0}, 12\tilde{3}/\tilde{4}/\cdots/\tilde{n}]$   is in $[\hat{0}, [6]^j/\tilde{7}/\cdots/ \tilde{n}]$ and so we have the same problem as before. If $j=1,2,3$, then we claim the interval $[123^{j}/\tilde{4}/\dots /\tilde{n}, 123^j/45\tilde{6}/\tilde{7}\cdots/\tilde{n}]$ has two increasing chains.   The chain 
     $$
     123^{j}/\tilde{4}/\tilde{5}/\tilde{6}/\cdots /\tilde{n} \cover 123^{j}/\tilde{4}5/\tilde{6}/\cdots /\tilde{n} \cover 123^{j}/45\tilde{6}/\cdots /\tilde{n}
     $$
    has label sequence $(5,7)(6,7)$ and the chain
      $$
     123^{j}/\tilde{4}/\tilde{5}/\tilde{6}/\cdots /\tilde{n} \cover 123^{j}/4\tilde{5}/\tilde{6}/\cdots /\tilde{n} \cover 123^{j}/45\tilde{6}/\cdots /\tilde{n}
     $$
     has label sequence $(5,6)(6,7)$, both of which are increasing.  It follows that $\tilde \lambda$ is not an EL-labeling in general.
\end{remark}

\section{Combinatorial Descriptions of  Whitney duals} \label{sec:CombDescriptOfWhitneyDual}

In this section we give combinatorial descriptions of the Whitney duals of $\Pi_n^w$ and $\Pi_n^\bullet$ that come from the EW-labelings we discussed in Section \ref{sec:EWLab}.

\subsection{Constructing Whitney duals}\label{sec:constructingWD}
  We start with a quick review of how to construct Whitney duals from EW-labelings. The full details can be found in~\cite{GonzalezHallam2021}. Note that in~\cite{GonzalezHallam2021}, the authors introduce two (isomorphic) Whitney duals that can be constructed using an EW-labeling.  The first construction, which is denoted $Q_{\lambda}(P)$ in \cite{GonzalezHallam2021}, is obtained by taking a quotient of the poset of saturated chains containing $\hat{0}$. This quotient is based on a quadratic relation on these chains. The second construction, denoted $R_{\lambda}(P)$, is a poset on ascent-free saturated chains containing $\hat{0}$ and whose order relation is given by sorting a label into an existing ascent-free chain.  Here  we will   only use the sorting description.

Let $\Lambda$ be a poset of labels and consider the following sorting algorithm on words over $\Lambda$.  Given the word  $w=w_1w_2\cdots w_iw_{i+1}\cdots w_n$, let $i$ be the smallest index such that  $w_i<w_{i+1}$.  If no such $i$ exists, the algorithm terminates and returns $w$.  If $i$ does exists, swap $w_i$ and $w_{i+1}$ to get the word $w'=w_1w_2\cdots w_{i+1}w_i\cdots w_n$.  Next, apply this procedure to $w'$ and continue until the algorithm terminates. Once the algorithm terminates one obtains a unique ascent-free word with the same underlying multiset of labels as the original word. We denote this unique word as $\sort(w)$.

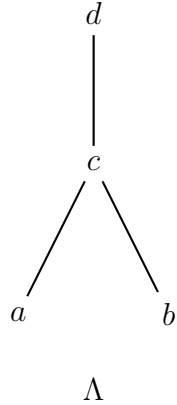
\begin{figure}
  \begin{tikzpicture}
    \node (a) at (-1,0) {$a$};
    \node (b) at (1,0) {$b$};
  \node (c) at (0,2) {$c$};
    \node (d) at (0,4) {$d$};

\draw[thick] (a)--(c)--(d);
 \draw[thick] (b)--(c);

\node at (0,-1) {$\Lambda$};

  \end{tikzpicture}
  \caption{A poset of labels}\label{fig:posetOfLabels}
\end{figure}

As an example of the sorting algorithm described above, consider the poset of labels  $\Lambda$ shown in Figure~\ref{fig:posetOfLabels}.  Applying the sorting algorithm to $adbca$ gives the following.
$$
abdca \rightarrow adbca \rightarrow dabca \rightarrow dacba \rightarrow dcaba
$$
Thus, $\sort(adbca) = dcaba$.

\begin{definition}[Definition 4.2 \cite{GonzalezHallam2021}]\label{definition:R_lambda}
Suppose that $P$ is a poset with an EW-labeling $\lambda$.  Let $R_\lambda(P)$ be the set of pairs $(x,w)$ where $x \in P$ and $w$ is the sequence of labels along a $\hat{0}-x$ ascent-free  chain.  Order the elements of $R_\lambda(P)$ by $(x,w) \cover (y,u)$ if and only  if  $x\cover y$ in $P$ and $u =\sort(w\lambda(x\cover y))$, where $w\lambda(x\cover y)$ denotes the concatenation of the words $w$ and $\lambda(x\cover y)$ and sort is done with respect to the ordering of the labels of $\lambda$.   
\end{definition}

\begin{theorem}[Theorem 4.4 \cite{GonzalezHallam2021}]\label{thm:WhitneyDualRLambda}
  Suppose that $P$ is a poset with an EW-labeling $\lambda$.  Then $P$ and $R_\lambda(P)$ are Whitney duals.
\end{theorem}

\begin{figure}
    \begin{tikzpicture}[]

\begin{scope}[yshift=150,xshift=130,scale=1.4]
\tikzstyle{every node}=[inner sep=0pt, scale=0.7*1.2, minimum width=4pt]
\node (n1232) at (3,4) {$123^{2}$};
  \node (n13020) at (-3,2) {$13^{ 0}/ 2^{ 0}$};
  \node (n102030) at (0,0)  {$1^{0}/ 2^{0}/ 3^{0}$};
  \node (n1231) at (0,4) {$123^{1}$};
  \node (n12030) at (-5,2) {$12^{0}/ 3^{0}$};
  \node (n13120) at (3,2)  {$13^{1}/ 2^{0}$};
  \node (n1230) at (-3,4){$123^ {0}$};
  \node (n10230) at (-1,2)  {$1^{0}/ 23^{0}$};
  \node (n12130) at (1,2)  {$12^{ 1}/ 3^{0}$};
  \node (n10231) at (5,2) {$23^ {1}/ 1^{0}$};

  \draw (n1231) -- (n10230)  node [midway,sloped, above] {$\textcolor{blue}{(1,2)^1}$};
  \draw [] (n13020) -- (n102030)   node [midway,sloped, above] {$\textcolor{blue}{(1,3)^0}$};
  \draw [] (n1232) -- (n13120)  node [midway,sloped, above] {$\textcolor{blue}{(1,2)^1}$};
  \draw [] (n1231)-- (n13020)node [midway,sloped, above] {$\textcolor{blue}{(1,2)^1}$};
  \draw [] (n10230)--(n102030) node [midway,sloped, above] {$\textcolor{blue}{(2,3)^0}$};
  \draw [] (n1230) -- (n10230) node [near start, sloped, above] {$\textcolor{blue}{(1,2)^0}$};
  \draw [] (n1231) -- (n13120) node [midway,sloped, above] {$\textcolor{blue}{(1,2)^0}$};
  \draw [] (n12030)-- (n102030) node [midway,sloped, above] {$\textcolor{blue}{(1,2)^0}$};
  \draw [] (n1231) --(n12130)node [midway,sloped, above] {$\textcolor{blue}{(1,3)^0}$};
  \draw [] (n1232) -- (n12130) node [near start,sloped, above] {$\textcolor{blue}{(1,3)^1}$};
  \draw [] (n13120) --(n102030) node [midway,sloped, above] {$\textcolor{blue}{(1,3)^1}$};
  \draw [] (n1231) --(n10231) node [near start,sloped, above] {$\textcolor{blue}{(1,2)^0}$};
  \draw [] (n1230) -- (n13020) node [near start,sloped, above] {$\textcolor{blue}{(1,2)^0}$};
  \draw [] (n1230)  -- (n12030) node [midway,sloped, above] {$\textcolor{blue}{(1,3)^0}$};
  \draw [] (n12130) --  (n102030) node [midway,sloped, above] {$\textcolor{blue}{(1,2)^1}$};
  \draw [] (n1232)  --  (n10231) node [midway,sloped, above] {$\textcolor{blue}{(1,2)^1}$};
  \draw [] (n10231)  --  (n102030)node [midway,sloped, above] {$\textcolor{blue}{(2,3)^1}$};
  \draw [] (n1231) -- (n12030) node [near start,sloped, above] {$\textcolor{blue}{(1,3)^1}$};	
  
   \tikzstyle{every node}= [scale=1]
    \node at (0,-.5) {$\Pi_3^w$};

  \end{scope}

\begin{scope}[scale=.35]
\tikzstyle{every node}=[inner sep=0pt,minimum size=10,scale=0.45*1.2]
\node (p1) at (13,0.25) {\Large$\emptyset$};
\node (q1) at (-4.5,5.25) {\Large$(2,3)^0$};
\node (q2) at (5.5,5.25) {\Large$(1,3)^0$};
\node (q3) at (10.5,5.25) {\Large$(1,2)^0$};
\node (q4) at (15.5,5.25) {\Large$(1,3)^1$};
\node (q5) at (20.5,5.25) {\Large$(1,2)^1$};
\node (q6) at (30.5,5.25) {\Large$(2,3)^1$};

\node (r1) at (-7,10.25) {\Large$(2,3)^0(1,2)^0$};
\node (r2) at (-2,10.25) {\Large$(2,3)^0(1,2)^1$};
\node (r3) at (3,10.25) {\Large$(1,3)^0(1,2)^1$};
\node (r4) at (8,10.25) {\Large$(1,3)^0(1,2)^0$};
\node (r5) at (13,10.25) {\Large$(1,3)^1(1,2)^0$};
\node (r6) at (18,10.25) {\Large$(1,3)^1(1,2)^1$};
\node (r7) at (23,10.25) {\Large$(1,2)^1(1,3)^0$};
\node (r8) at (28,10.25) {\Large$(2,3)^1(1,2)^0$};
\node (r9) at (33,10.25) {\Large$(2,3)^1(1,2)^1$};

\tikzstyle{every path}=[]
 \draw (p1) -- (q1) ;
 \draw (p1) -- (q2) ;
 \draw (p1) -- (q3) ;
 \draw (p1) -- (q4) ;
 \draw (p1) -- (q5) ;
 \draw (p1) -- (q6) ;

 \draw (q1) -- (r1) ;
 \draw (q1) -- (r2) ;
 \draw (q2) -- (r3) ;
 \draw (q2) -- (r4) ;
 \draw (q3) -- (r4) ;
 \draw (q3) -- (r5) ;
 \draw (q4) -- (r5) ;
 \draw (q4) -- (r6) ;
 \draw (q5) -- (r6) ;
 \draw (q5) -- (r7) ;
 \draw (q6) -- (r8) ;
 \draw (q6) -- (r9) ;

 \tikzstyle{every node}= [scale=1]

\node at (13,-1) {$R_{\lambda_w}(\Pi_3^w)$};
\end{scope}

\end{tikzpicture}
\caption{The weighted partition poset an with an EW-labeling and  the corresponding Whitney dual.}
\label{fig:Pi3wAndWhitneyDual}
\end{figure}
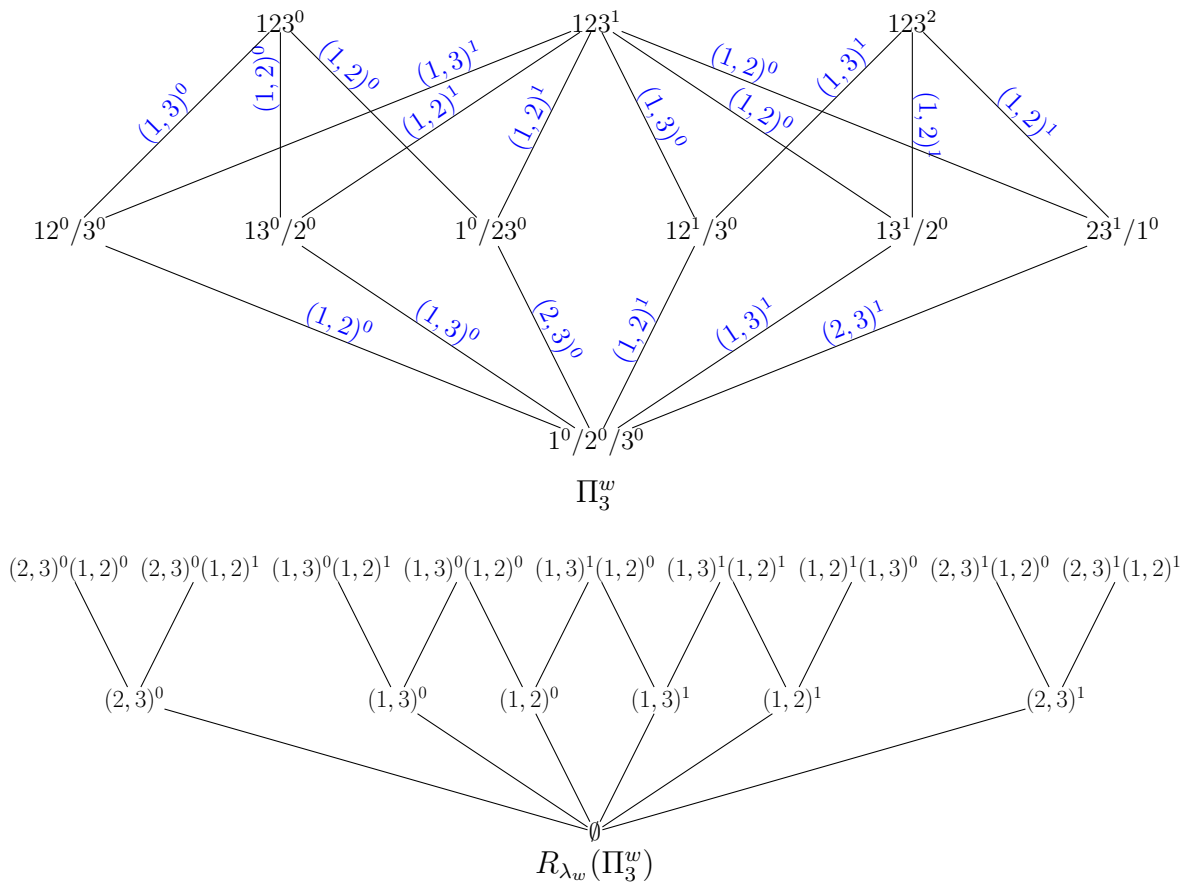
The reader can verify that the poset labeled $Q$ in Figure~\ref{fig:WhitneyDualExample} is $R_\lambda(P)$ where $\lambda$ is the EW-labeling of $P$ given in the figure. In addition to being EW-labelings, where ascent-free chains have a unique word of labels in their own interval, $\lambda_w$ and $\lambda_\bullet$ have the additional property that the sequence of labels along a saturated chain starting at $\hat{0}$ completely determines the elements on that chain. So the word of labels identifies the saturated chain uniquely in the poset and as a result, when we consider $R_\lambda(P)$ the ``$x$" in the pair $(x,w)$ is redundant. In other words, we need only to consider ascent-free chains as elements of $R_\lambda(P)$.  
Figure~\ref{fig:Pi3wAndWhitneyDual} depicts $\Pi_3^w$ and the Whitney dual corresponding to the EW-labeling described in the previous section. See Figure~\ref{fig:FLyn_3_weighted} for an isomorphic version of $R_{\lambda_w}(\Pi_n^w)$ whose elements are described in terms of a family of forests.

\subsection{Combinatorial families for the Whitney duals} \label{section:combinatorial_families}
We now turn our attention to giving combinatorial descriptions of the Whitney duals of $\Pi_n^\bullet$ and $\Pi_n^w$. First, we need to describe the combinatorial objects on which the Whitney duals will be defined. 

A \emph{tree} is an undirected graph in which any two vertices are connected by exactly one path. We say that a tree is \emph{rooted} if there is a distinguished vertex that we call the \emph{root}.  If, in order to travel through the unique path from a vertex $b$ to the root we need to pass through a vertex $a$, we say that $a$ is an \emph{ancestor} of $b$. In the specific cases that $\{a,b\}$ is an edge, we say that $a$ is the \emph{parent} of $b$ (or equivalently, $b$ is a \emph{child} of $a$). Every vertex in a rooted tree $T$ which has at least one child is considered an \emph{internal vertex}. If it has no child we say that the vertex is a \emph{leaf}. A \emph{planar tree} is a rooted tree in which the set of children of each internal vertex comes equipped with a total order (which we represent by placing the vertices from left to right in this order). A \emph{binary tree} is a rooted planar tree in which every internal vertex has two children, a \emph{left child} and a \emph{right child}.  All the trees we consider from now on are both rooted and planar, so we will be referring to them (informally) as ``binary trees" when the context makes it clear.

We say a binary tree is a \emph{bicolored binary tree} if there is a function $\clr$ that assigns to each internal vertex $x$ a number $\clr(x)\in \{\tbl{0},\tre{1}\}$ (a \emph{color}).  Note that in all of our figures, we represent the color \tbl{0} with blue and the color \tre{1} with red.

\begin{figure}
    \centering
    
     \begin{subfigure}[t]{0.45\textwidth}
         \centering
         \begin{tikzpicture}[scale=0.8]
     
     \tikzstyle{every node}=[draw,scale=0.8]
     \node[fill=blue,scale=1.2](one){};
     \node at ($(one)+(-1.5,-1.5)$)[circle,fill=red](two){};
     \node at ($(two)+(-1,-1)$)[circle,fill=red](three){};
     \node at ($(three)+(-0.75,-1)$)[draw=none](3son1){1};
     \node at ($(three)+(0.75,-1)$)[draw=none](3son2){9};
     \node at ($(two)+(1,-1)$)[fill=blue,scale=1.2](four){};
     \node at ($(four)+(-0.75,-1)$)[draw=none](4son){2};
     \node at ($(four)+(0.75,-1)$)[circle,fill=red](4sonn){};
     \node at ($(one)+(1.5,-1.5)$)[fill=blue,scale =1.2](five){};
     \node at ($(five)+(0.75,-1)$)[circle,fill=red](six){};
     \node at ($(six)+(-0.75,-1)$)[draw=none](6son){4};
     \node at ($(4sonn)+(-0.75,-1)$)[draw=none](4sonn1){5};
     \node at ($(4sonn)+(0.75,-1)$)[draw=none](4sonn2){8};
     \node at ($(five)+(-0.75,-1)$)[draw=none](5son1){3};
     \node at ($(six)+(0.75,-1)$)[circle,fill=red](6sonn){};
     \node at ($(6sonn)+(-0.75,-1)$)[draw=none](last1){6};
     \node at ($(6sonn)+(0.75,-1)$)[draw=none](last2){7};
     \path[thick](six)edge(6son)edge(6sonn);
     \path[thick](6sonn)edge(last1)edge(last2);
     \path[thick](five)edge(5son1);
     \path[thick](4sonn)edge(4sonn1)edge(4sonn2);
     \path[thick](one)edge(two)edge(five);
     \path[thick](two)edge(three)edge(four);
     \path[thick](five)edge(six);
     \path[thick](three)edge(3son1)edge(3son2);
     \path[thick](four)edge(4son)edge(4sonn);

     \node[fill=blue,scale=1.2,below of=two,align=center,anchor=center,xshift=0.25cm,yshift=-2cm](label1){};
     \node at ($(label1)+(1,0)$)[draw=none](label11){0};
     \node[circle,fill=red,right of=label11,align=center,anchor=center](label2){};
     \node at ($(label2)+(1,0)$)[draw=none](label22){1};

     \end{tikzpicture}
         \caption{Bicolored tree}
         \label{fig:bicolored_tree}
     \end{subfigure}
     \begin{subfigure}[t]{0.45\textwidth}
         \centering
              \begin{tikzpicture}[scale=0.8]
	     \tikzstyle{every node}=[draw,scale=0.8]
	     \node[	blue,scale=1.2](one){8};
	     \node at ($(one)+(-1.5,-1.5)$)[circle,red](two){7};
	     \node at ($(two)+(-1,-1)$)[circle,red](three){6};
	     \node at ($(three)+(-0.75,-1)$)[draw=none](3son1){1};
	     \node at ($(three)+(0.75,-1)$)[draw=none](3son2){9};
	     \node at ($(two)+(1,-1)$)[blue,scale=1.2](four){5};
	     \node at ($(four)+(-0.75,-1)$)[draw=none](4son){2};
	     \node at ($(four)+(0.75,-1)$)[circle,red](4sonn){2};
	     \node at ($(one)+(1.5,-1.5)$)[blue,scale =1.2](five){4};
	     \node at ($(five)+(0.75,-1)$)[circle,red](six){3};
	     \node at ($(six)+(-0.75,-1)$)[draw=none](6son){4};
	     \node at ($(4sonn)+(-0.75,-1)$)[draw=none](4sonn1){5};
	     \node at ($(4sonn)+(0.75,-1)$)[draw=none](4sonn2){8};
	     \node at ($(five)+(-0.75,-1)$)[draw=none](5son1){3};
	     \node at ($(six)+(0.75,-1)$)[circle,red](6sonn){1};
	     \node at ($(6sonn)+(-0.75,-1)$)[draw=none](last1){6};
	     \node at ($(6sonn)+(0.75,-1)$)[draw=none](last2){7};
	     \path[thick](six)edge(6son)edge(6sonn);
	     \path[thick](6sonn)edge(last1)edge(last2);
	     \path[thick](five)edge(5son1);
	     \path[thick](4sonn)edge(4sonn1)edge(4sonn2);
	     \path[thick](one)edge(two)edge(five);
	     \path[thick](two)edge(three)edge(four);
	     \path[thick](five)edge(six);
	     \path[thick](three)edge(3son1)edge(3son2);
	     \path[thick](four)edge(4son)edge(4sonn);
	
	     \node[fill=blue,scale=1.2,below of=two,align=center,anchor=center,xshift=0.25cm,yshift=-2cm](label1){};
	     \node at ($(label1)+(1,0)$)[draw=none](label11){0};
	     \node[circle,fill=red,right of=label11,align=center,anchor=center](label2){};
	     \node at ($(label2)+(1,0)$)[draw=none](label22){1};

	     \end{tikzpicture}
         \caption{The internal vertices have been labeled in the reverse-minimal linear extension of the tree to the left}
         \label{fig:bicolored_tree_reverse_minimal_linear_extension}
     \end{subfigure}

    \caption{Example of a bicolored binary tree and its reverse-minimal linear extension}
    \label{fig:bicolored_tree_example}
\end{figure}
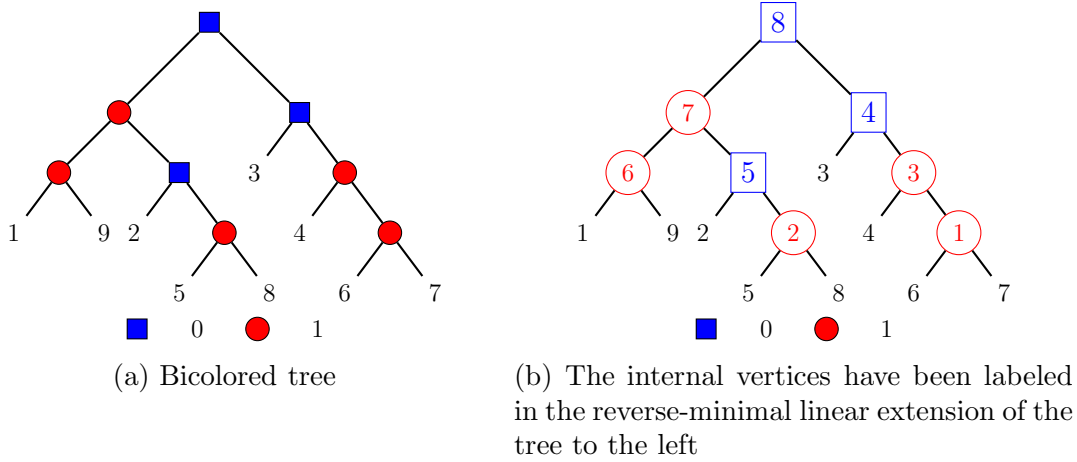

A \emph{linear extension} of a binary tree $T$ is a listing $v_1,v_2, \dots,v_{n-1}$ of the internal vertices of $T$ such that each vertex precedes its parent. 
Let $T$ be a bicolored binary tree and $v$ a vertex of $T$. We define the valency of $v$, $\nu(v)$,  to be the smallest leaf label of the subtree rooted at $v$.  Note that, by this definition, if $w$ is an ancestor of $v$ we have that $\nu(v)\ge \nu(w)$. Hence, since the leaves are labeled by the totally ordered set $[n]$, there is a unique linear extension $v_1,v_2,\dots,v_{n-1}$  of $T$ such that
\begin{align}\label{equation:reverse_minimal_le}
    \nu(v_1)\ge \nu(v_2) \ge \cdots \ge \nu(v_{n-1}).
\end{align}
We will call this linear extension the \emph{reverse-minimal} linear extension of the internal vertices of $T$. Figure \ref{fig:bicolored_tree_reverse_minimal_linear_extension} depicts the reverse-minimal linear extension of the tree in Figure \ref{fig:bicolored_tree}.     Note that we can extend the notion of reverse linear extension to forests with leaf set $[n]$.   For example, see the forest in Figure~\ref{fig:pointed_Lyndon_Tree_and_chain}. As we will see later, the reverse-minimal linear extension gives us a recipe to build  a forest  step-by-step in a way that corresponds to the ascent-free chains of the weighted and pointed partition posets.

Let $v$ be an internal vertex of a bicolored binary tree $T$. We denote as $L(v)$ the left child of $v$ and as $R(v)$ the right child of $v$. 
We say that $T$ is \emph{normalized} if for every internal vertex $v$ we have that
\begin{equation}\label{equation:normalized}\tag{N}
    \nu(v)=\nu(L(v)).
\end{equation}
In other words, a tree is normalized if  the smallest leaf label always appears in a leaf to its left.  One can check that the tree depicted in Figure~\ref{fig:bicolored_tree_example} is normalized. We say a forest is normalized if all the trees in the forest are normalized. Whenever $T$ is normalized, we say that an internal vertex $v$ is \emph{Lyndon} if $L(v)$ is a leaf or else if $L(v)$ is not a leaf then we have that
 \begin{equation}\label{eq:1}\tag{L}
     \nu(R(L(v)))>\nu(R(v)).
 \end{equation}
 Returning to our example in Figure~\ref{fig:bicolored_tree_example}, we see that the internal vertices that are labeled by $1,2,3,4,5,6$ are all Lyndon since for each,  their left child is a leaf. The internal vertex labeled as $7$ is also Lyndon but $8$ is not. To see why $7$ is Lyndon, note that $R(L(7))$ is the leaf labeled $9$ (and hence has valency 9) and $R(7)=5$ which has valency $2$. Thus, $ \nu(R(L(7)))>\nu(R(7))$. To see that $8$ is not Lyndon, note that  $\nu(R(L(8)))=2\not> 3= \nu(R(8))$, i.e.~inequality \ref{eq:1} is not satisfied.

The Whitney duals for the weighted partition poset and the pointed partition poset can both   be described using special types of normalized bicolored binary trees together with a sliding procedure used to merge   such trees.   We now discuss these special types of trees.

\subsubsection{Pointed Lyndon forests}

Let us define first the objects used to describe the Whitney dual of the pointed partition poset.

\begin{definition}
A normalized bicolored binary tree $T$ is said to be a \emph{pointed Lyndon tree} if for each internal vertex $v$ of $T$ such that $L(v)$ is not a leaf it must be that:
\begin{equation}\tag{PL1}\label{condition:pointed_lyndon_1}
\clr(L(v))\geq \clr(v)
\end{equation}
\begin{center} and \end{center}
\begin{equation}\label{condition:pointed_lyndon_2}
     \tag{PL2} \text{If } \clr(L(v))=\clr(v)=\tre{1}\text{, then $v$ is a Lyndon node.}
\end{equation}
 We say a forest $F$ is a \emph{pointed Lyndon forest} if  all the connected components of $F$ are pointed Lyndon trees.  We denote as   $\Flyn_{n}^{\bullet}$ the set of all pointed Lyndon forests whose leaf label set is $[n]$.   
\end{definition}

The tree in Figure~\ref{fig:bicolored_tree_example} is a pointed Lyndon tree. Indeed, we need only check the conditions for the internal vertices labeled by 7 and 8 since the left children of the other internal vertices are leaves.  For 7, both it and its left child 6 are colored by \tre{1} (red) and 7 is Lyndon. For 8, it is colored \tbl{0} (blue) so we do not need to check anything else.

\begin{figure}
     \begin{tikzpicture}[every node/.style={draw=black},every 
node/.append style={transform shape},scale=.8]
     \node[fill=red,circle](one){};
     \node at ($(one)+(0.75,-1)$)[draw=none](two){2};
     \node at ($(one)+(-0.75,-1)$)[scale=1.2,fill=blue](three){};
     \node at ($(three)+(-0.75,-1)$)[draw=none](four){1};
     \node at ($(three)+(0.75,-1)$)[draw=none](five){3};
     \path[thick](one)edge(two)edge(three);
     \path[thick](three)edge(four)edge(five);
     \node at (0,-3) [draw=none] {\Large $T_1$};
     \begin{scope}[shift={(5,0)}, every node/.style={draw=black},every 
node/.append style={transform shape}]
     \node[fill=blue,scale=1.2](one){};
     \node at ($(one)+(0.75,-1)$)[draw=none](two){3};
     \node at ($(one)+(-0.75,-1)$)[scale=1.2,fill=blue](three){};
     \node at ($(three)+(-0.75,-1)$)[draw=none](four){1};
     \node at ($(three)+(0.75,-1)$)[draw=none](five){$2$};
     \path[thick](one)edge(two)edge(three);
     \path[thick](three)edge(four)edge(five);
          \node at (0,-3) [draw=none] {\Large $T_2$};

    \end{scope}

       \node[fill=blue,scale=1.2,below of=two,align=center,anchor=center] at (1,-3)(label1){};
     \node at ($(label1)+(1,0)$)[draw=none](label11){0};
     \node[circle,fill=red,right of=label11,align=center,anchor=center](label2){};
     \node at ($(label2)+(1,0)$)[draw=none](label22){1};

        \end{tikzpicture}
        \caption{A bicolored Lyndon tree on the left (which is not a pointed Lyndon tree) and a pointed Lyndon tree on the right (which is not a bicolored Lyndon tree).}. \label{fig:bicoloredAndPointed}
\end{figure}
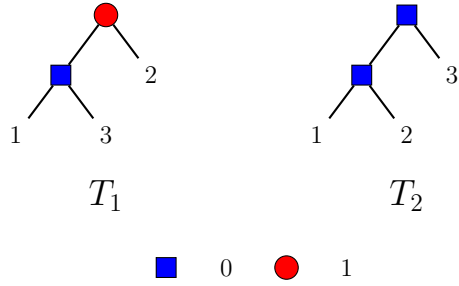

\subsubsection{Bicolored Lyndon forests}

Let us now define the trees used for the Whitney dual of the weighted partition poset.

\begin{definition}
Let $T$ be a normalized bicolored binary tree, we say that $T$ is a \textit{bicolored Lyndon tree} if, for each internal vertex $v$ of $T$ whose left child is an internal vertex, either $v$ is a Lyndon vertex or it must be that:
\begin{equation}\label{cleq:2}
 \tag{CL}    \clr(L(v))>\clr(v).
\end{equation}
We say a forest $F$ is a \emph{bicolored Lyndon forest} if  all the connected components of $F$ are bicolored Lyndon trees.   We denote as   $\Flyn_{n}^{w}$ the set of all bicolored Lyndon forests whose leaf label set is $[n]$.   
\end{definition}

The tree in Figure~\ref{fig:bicolored_tree_example} is a bicolored Lyndon tree.  As mentioned earlier, the only internal vertex that is not a Lyndon vertex is the one labeled by 8.   Its color is \tbl{0} (blue) and its left child is colored by \tre{1} (red), so it satisfies condition~\eqref{cleq:2}.

At this point, the reader may be wondering if pointed Lyndon trees and bicolored Lyndon trees are the same.  To see that this is not the case, consider the trees in Figure~\ref{fig:bicoloredAndPointed}.  We claim that the tree $T_1$ is bicolored, but not pointed.  All the internal vertices of $T_1$ are Lyndon, so it is automatically a bicolored Lyndon tree.  However,  it is not a pointed Lyndon tree because the color of the root is larger than its left child, a violation of condition~\eqref{condition:pointed_lyndon_1}.  On the other hand, we claim that $T_2$ is a pointed Lyndon tree, but not a bicolored Lyndon tree.  It is pointed  since both vertices are colored by \tbl{0} (blue).   However, it is not a bicolored Lyndon tree since the root is not a Lyndon vertex and the root's color is not strictly larger than its left child, a violation of condition~\eqref{cleq:2}.

\begin{remark} Note that condition ~\eqref{cleq:2} implies that the family of  classical Lyndon trees coincide with the subfamily of bicolored Lyndon trees that either have only internal vertices of color \tbl{0} (blue) or that have only internal vertices of color \tre{1} (red).  On the other hand condition ~\eqref{condition:pointed_lyndon_2} implies that the family of classical Lyndon trees coincide with the subfamily of pointed Lyndon trees which have all internal nodes colored \tre{1} (red).
\end{remark}

\subsection{A Whitney dual for \texorpdfstring{$\Pi_n^\bullet$}-} \label{subsection:whitney_dual_pointed}

 Let $F$ be a pointed Lyndon forest.  We explain how to associate an ascent-free saturated chain $c(F)$ starting at $\hat{0}$ in $\Pi_n^\bullet$.  Recall that $F$ has a unique reverse-minimal linear extension order on the internal vertices.  To construct an ascent-free chain from $F$, we start with the the bottom element, $\tilde{1}/\tilde{2}/\cdots /\tilde{n}$.  In the first step, we merge together the two blocks which are the in the left subtree and right subtree rooted at the first internal vertex.   We keep the point on the element of the left tree if the first internal vertex is colored \tre{1} (red) and we keep the point on the element of the right tree if it is  colored is \tbl{0} (blue).   We continue doing this so that at the $i^{th}$ step we merge together the elements of the left and right subtrees rooted at the $i^{th}$ internal vertex and keep the pointed element of the left subtree if it is colored \tre{1} (red) or the right subtree if it is colored \tbl{0} (blue). 

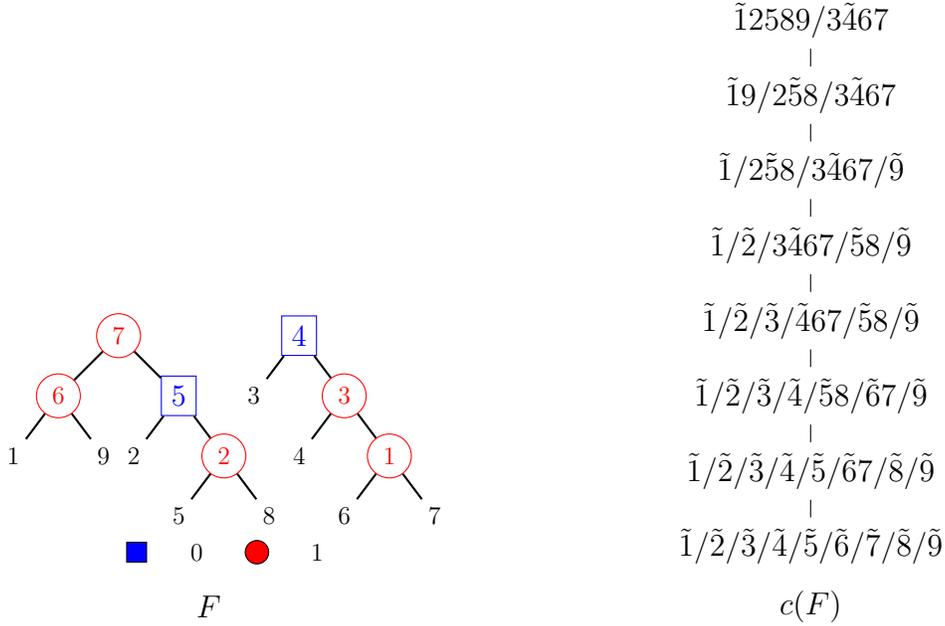
\begin{figure}
                 \begin{tikzpicture}[scale=0.8]
	    
	     \node[](one){};
	     	    \node[] at (0,-6) {$F$};

	     \node at ($(one)+(-1.5,-1.5)$)(two){};
	     	     \tikzstyle{every node}=[draw,scale=0.8]
	     \node at ($(one)+(-1.5,-1.5)$)[circle,red](two){7};

	     \node at ($(two)+(-1,-1)$)[circle,red](three){6};
	     \node at ($(three)+(-0.75,-1)$)[draw=none](3son1){1};
	     \node at ($(three)+(0.75,-1)$)[draw=none](3son2){9};
	     \node at ($(two)+(1,-1)$)[blue,scale=1.2](four){5};
	     \node at ($(four)+(-0.75,-1)$)[draw=none](4son){2};
	     \node at ($(four)+(0.75,-1)$)[circle,red](4sonn){2};
	     \node at ($(one)+(1.5,-1.5)$)[blue,scale =1.2](five){4};
	     \node at ($(five)+(0.75,-1)$)[circle,red](six){3};
	     \node at ($(six)+(-0.75,-1)$)[draw=none](6son){4};
	     \node at ($(4sonn)+(-0.75,-1)$)[draw=none](4sonn1){5};
	     \node at ($(4sonn)+(0.75,-1)$)[draw=none](4sonn2){8};
	     \node at ($(five)+(-0.75,-1)$)[draw=none](5son1){3};
	     \node at ($(six)+(0.75,-1)$)[circle,red](6sonn){1};
	     \node at ($(6sonn)+(-0.75,-1)$)[draw=none](last1){6};
	     \node at ($(6sonn)+(0.75,-1)$)[draw=none](last2){7};
	     \path[thick](six)edge(6son)edge(6sonn);
	     \path[thick](6sonn)edge(last1)edge(last2);
	     \path[thick](five)edge(5son1);
	     \path[thick](4sonn)edge(4sonn1)edge(4sonn2);

	     \path[thick](two)edge(three)edge(four);
	     \path[thick](five)edge(six);
	     \path[thick](three)edge(3son1)edge(3son2);
	     \path[thick](four)edge(4son)edge(4sonn);
	
	     \node[fill=blue,scale=1.2,below of=two,align=center,anchor=center,xshift=0.25cm,yshift=-2cm](label1){};
	     \node at ($(label1)+(1,0)$)[draw=none](label11){0};
	     \node[circle,fill=red,right of=label11,align=center,anchor=center](label2){};
	     \node at ($(label2)+(1,0)$)[draw=none](label22){1};

	    \begin{scope} [shift = {(10, -5)}]
	    	     \tikzstyle{every node}=[]
         \node[] at (0,-1) {$c(F)$};
	       \node (0) at (0,0) {$\tilde{1}/\tilde{2}/\tilde{3}/\tilde{4}/\tilde{5}/\tilde{6}/\tilde{7}/\tilde{8}/\tilde{9}$};
	       
	       \node (1) at (0,1.25) {$\tilde{1}/\tilde{2}/\tilde{3}/\tilde{4}/\tilde{5}/\tilde{6} 7/\tilde{8}/\tilde{9}$};
	       
	       \draw (0)--(1);
	       
	         \node (2) at (0,2.5) {$\tilde{1}/\tilde{2}/\tilde{3}/\tilde{4}/\tilde{5}8/\tilde{6}7/\tilde{9}$};
	       
	       \draw (1)--(2);
	       
	       \node (3) at (0,3.75) {$\tilde{1}/\tilde{2}/\tilde{3}/\tilde{4}67 /\tilde{5}8/\tilde{9}$};
	       
	       \draw (2)--(3);
	       
	       \node (4) at (0,5) 
	       {$\tilde{1}/\tilde{2}/3\tilde{4}67 /\tilde{5}8/\tilde{9}$};
	       
	        \draw (3)--(4);
	        
	        \node (5) at (0,6.25) {$\tilde{1}/2\tilde{5}8/3\tilde{4}67/\tilde{9}$};
	        
	        \draw (4)--(5);
	        
	         \node (6) at (0,7.5) {$ \tilde{1}9/2\tilde{5}8/3\tilde{4}67 $};
	        
	        \draw (5)--(6);
	        
	         \node (7) at (0,8.75) {$ \tilde{1}2589/3\tilde{4}67 $};
	         
	         \draw (6)--(7);

	    \end{scope}
	   
	     \end{tikzpicture}
	     
	     \caption{Pointed Lyndon forest $F$ and corresponding ascent-free chain $c(F)$.}\label{fig:pointed_Lyndon_Tree_and_chain}
         
\end{figure}

 Figure~\ref{fig:pointed_Lyndon_Tree_and_chain} has  a depiction of a pointed Lyndon forest $F$ and its corresponding chain $c(F)$.  In the first step of the chain, we merge the blocks containing $6$ and $7$ since they are the leafs in the left and right tree rooted at the first internal vertex. We keep $6$ pointed because this first internal vertex is colored by \tre{1} (red). Continuing this process gives the saturated chain seen in the figure.

Let $F$ be a pointed Lyndon forest and let $a_i$ be the valency of left child of the $i^{th}$ internal vertex  and $b_i$ the valency of right child of $i^{th}$  internal vertex.  That is $a_i = \nu(L(i))$ and $b_i =\nu(R(i))$. Note that $a_i$ is also the valency of $i$ since our trees are normalized (i.e.~$\nu(i) =\nu(L(i)))$.  Then, it is straightforward to see that the sequence of labels along the chain $c(F)$ is $(a_1,b_1)^{u_1}(a_2,b_2)^{u_2}\cdots (a_{k},b_{k})^{u_{k}}$ where $u_i$ is the color of the $i^{th}$ internal vertex.

Let us return to the forest in Figure~\ref{fig:pointed_Lyndon_Tree_and_chain}. Since the left child of the first internal vertex is the leaf labeled 6, the right child is a leaf labeled 7, and the first internal vertex is colored \tre{1} (red), we have that $(a_1,b_1)^{u_1} = (6,7)^1$. Moving to the second internal vertex we see its left child is the leaf $5$ and the right child is the leaf $8$. Since it is colored by \tre{1} (red) the next label is $(5,8)^{\tre{1}}$.  Next, moving to the third internal vertex, wee see that the left child is a leaf $4$ and the right child is the first internal vertex  whose valency is 6. Since the third internal vertex is colored \tre{1} (red), the corresponding label is $(4,6)^{\tre{1}}$.   Continuing this, we see the sequence we get from the forest is 
$$
(6,7)^{\tre{1}} (5,8)^{\tre{1}} (4,6)^{\tre{1}}  (3,4)^{\tbl{0}} (2,5)^{\tbl{0}} (1,9)^{\tre{1}} (1,2)^{\tre{1}}
$$
 Note that the sequence (and hence its chain) is ascent-free with respect to the ordering of the labels for $\lambda_\bullet$.   It turns out that this not a coincidence as we show next.

\begin{theorem}\label{thm:bijectionBetweenForestAndChain}
The map sending $F$ to $c(F)$ is a bijection between pointed Lyndon forest whose leaf label set is $[n]$ and ascent-free chains starting at $\hat{0}$ of $\Pi_n^\bullet$, where the  ascent-free condition is defined with respect to $\lambda_\bullet$.
\end{theorem}
\begin{proof}
First, we show that the map  is well-defined.   That is, $c(F)$ is in fact ascent-free for all pointed Lyndon forest $F$.  Assume that the internal vertices are  $1,2,\dots, k$ which is also the reverse-minimal ordering.   Let $(a_1,b_1)^{u_1}(a_2,b_2)^{u_2}\cdots (a_{k},b_{k})^{u_{k}}$ be the sequence of labels along $c(F)$.   As mentioned earlier, $a_i$ is the valency of $i$. Since we are using the reverse-minimal order, this gives us $a_1\geq a_2\geq \cdots \geq a_k$.  If $a_i> a_{i+1}$, then $(a_i,b_i)^{u_i} > (a_{i+1},b_{i+1})^{u_{i+1}}$ in the ordering of labels of $\lambda_\bullet$.  On the other hand, if $a_i=a_{i+1}$,  then  $i+1$ is an ancestor of $i$ and $i$ must be in the left tree rooted at $i+1$.  Since the reverse-minimal ordering is a linear extension of the internal vertices, it must be the case that $i$ is the left child of $i+1$.  Since $F$ is a  pointed Lyndon forest, condition~\eqref{condition:pointed_lyndon_1} implies that $u_i = \clr(i) \geq \clr(i+1) =u_{i+1}$.  If $i+1$ is not Lyndon, then condition~\eqref{condition:pointed_lyndon_2} implies that $u_i>u_{i+1}$. Since $u_i,u_{i+1}\in \{\tbl{0},\tre{1}\}$, this implies that $u_i=\tre{1}$ and $u_{i+1}=\tbl{0}$. So, if $i+1$ is not Lyndon, $(a_i,b_i)^{u_i} > (a_{i+1},b_{i+1})^{u_{i+1}}$.  On the other hand, if $i+1$ is Lyndon, $b_{i+1}=\nu(R(i+1)) <\nu(R(L(i+1))) =b_{i}$.  Then either $(a_i,b_i)^{u_i} > (a_{i+1},b_{i+1})^{u_{i+1}}$ (if $u_i>u_{i+1}$  or $u_i=1=u_{i+1}$) or   $(a_i,b_i)^{u_i}$ and $(a_{i+1},b_{i+1})^{u_{i+1}}$  are incomparable (if $u_i=\tbl{0}=u_{i+1})$.   It follows that the map is well-defined.
  
  We claim that the map sending $F$ to $c(F)$ is invertible.  Suppose we have an ascent-free chain with sequence 
  $$
  (a_1,b_1)^{u_1}(a_2,b_2)^{u_2}\cdots (a_{k},b_{k})^{u_{k}}.
  $$
  Build a forest recursively by first  placing $n$ isolated vertices (which will be the leafs) labeled by $[n]$. Now assume that you are at the $i^{th}$ step of this process.    Let $T_1$ be the connected component of the forest with minimal leaf label $a_i$ and let $T_2$ be the connected component of the forest with minimal leaf label $b_i$.  Add a vertex colored $u_i$  and add edges from this vertex to the roots of $T_1$ and $T_2$.   Repeat this process until each pair $(a_j,b_j)^{u_j}$ has been used and call the resulting forest $F$. Since $  (a_1,b_1)^{u_1}(a_2,b_2)^{u_2}\cdots (a_{k},b_{k})^{u_{k}}$ is ascent-free, it must be the case that $a_1\geq a_2\geq \cdots \geq a_k$.  It then follows that the reverse minimal ordering on $F$ is exactly the order that the internal vertices were added in the process. Upon observing this, it is clear that $c(F)$ is  the chain with label  sequence $(a_1,b_1)^{u_1}(a_2,b_2)^{u_2}\cdots (a_{k},b_{k})^{u_{k}}$. So, if we can show this map is well-defined, we will have that the map sending $F$ to $c(F)$ is invertible and thus a bijection.  We do this next.

  Let $F$ be a forest obtained by using the inverse procedure described in the previous paragraph. We need to show that $F$ is  a pointed Lyndon forest. First, it is clear from the construction that $F$ is a bicolored binary tree.  Since $a_i<b_i$ for all $i$, we also have that $F$ is  normalized.  Now consider an internal vertex $v$  and suppose that $v$ is the $i^{th}$ internal vertex in the reverse minimal order.  If $L(v)$ is not  a leaf,  then $L(v)$ must immediately precede $v$ in the  reverse minimal ordering.  That is, $L(v)$ is the $(i-1)^{th}$ internal vertex.  Since $F$ is normalized, $v$ and $L(v)$ have the same valency and so $a_{i-1} = a_{i}$.  Since $(a_1,b_1)^{u_1}(a_2,b_2)^{u_2}\cdots (a_{k},b_{k})^{u_{k}}$ is ascent-free, this means that $u_{i-1} \geq u_{i}$.  Thus, $\clr(L(v)) =u_{i-1} \geq u_i=\clr(v)$ and so condition~\eqref{condition:pointed_lyndon_1} is satisfied.  If $u_{i-1} =\tre{1} = u_{i}$, then the fact that $a_{i-1} = a_{i}$ implies that $b_{i-1}>b_i$ (otherwise either a label would be repeated in the sequence or the sequence would have an ascent).  Thus, 
  $\nu(R(L(v))= b_{i-1} > b_i = \nu(R(v))$ and so $v$ is a Lyndon vertex (equation~\eqref{eq:1}).  It follows that condition~\eqref{condition:pointed_lyndon_2} is also satisfied. We conclude that $F$ is a pointed Lyndon forest.  Thus, the inverse map is well-defined, completing the proof.
  \end{proof}

 By the previous theorem and Theorem~\ref{thm:WhitneyDualRLambda}, we can describe the Whitney dual of the pointed partition poset using pointed Lyndon forests.  To do this, we will need to describe how to merge trees in a Lyndon forest. Suppose that $T_1$ and $T_2$ are trees in a Lyndon forest with roots $r_1$ and $r_2$ where  the minimum leaf label of $T_1$ is less than the minimum leaf label of $T_2$.  Let $u\in \{\tbl{0},\tre{1}\}$.   To \emph{$u$-merge} $T_1$ with $T_2$, we first create a new vertex $r$ and color it so $\clr(r)=u$. Then we add edges from $r$ to $r_1$ and $r_2$.   If the resulting tree is a pointed Lyndon tree, we stop.  If it is not, we  \emph{slide} the new internal vertex $r$ together with its right subtree  past its left child and check if the result is a pointed Lyndon tree.  We continue this process until we obtain a pointed Lyndon tree.

 \begin{figure}
     \centering
     
\begin{tikzpicture}[scale=0.8]

\tikzstyle{every node} = [draw=black]
     \node[fill=blue,scale=1.2](root){};
     \node at ($(root)+(1.1,-1)$)[inner sep=3pt,fill=red,circle](one){};
     \node at ($(one)+(-0.75,-1)$)[draw=none](two){2};
     \node at ($(one)+(0.75,-1)$)[draw=none](three){3};
     \path[thick](one)edge(two)edge(three);
     \node at ($(root)+(-1.1,-1)$)[inner sep=3pt,fill=red,circle](one1){};
     \node at ($(one1)+(-0.75,-1)$)[draw=none](two1){1};
     \node at ($(one1)+(0.75,-1)$)[draw=none](three1){4};
     \path[thick](one1)edge(two1)edge(three1);
     \path[thick](root)edge(one1)edge(one);
     \node at ($(root)+(0,-2.6)$)[draw=none](arrow1){$T_1$};
     \node at ($(root)+(4,0)$)[scale=1.2,fill=blue](onee){};
     \node at ($(onee)+(0.75,-1)$)[draw=none](twoo){6};
     \node at ($(onee)+(-0.75,-1)$)[inner sep=3pt,fill=red,circle](threee){};
     \node at ($(threee)+(-0.75,-1)$)[draw=none](fourr){5};
     \node at ($(threee)+(0.75,-1)$)[draw=none](fivee){7};
     \path[thick](onee)edge(twoo)edge(threee);
     \path[thick](threee)edge(fourr)edge(fivee);
     \node at ($(onee)+(0,-2.6)$)[draw=none](arrow1){$T_2$};
     \node at ($(one)+(0.7,-2.4)$)[draw=none](o){};
     \node at ($(o)+(5,2)$)[draw=none](arrow1){$\longrightarrow$};
     \node at ($(o)+(10,3.5)$)[inner sep=3pt,red,circle](rootp){r};
     \node at ($(o)+(8,2.5)$)[fill=blue,scale=1.2](root1){};
     \node at ($(root1)+(1.1,-1)$)[inner sep=3pt,fill=red,circle](one1){};
     \node at ($(one1)+(-0.75,-1)$)[draw=none](two1){2};
     \node at ($(one1)+(0.75,-1)$)[draw=none](three1){3};
     \path[thick](one1)edge(two1)edge(three1);
     \node at ($(root1)+(-1.1,-1)$)[inner sep=3pt,fill=red,circle](one11){};
     \node at ($(one11)+(-0.75,-1)$)[draw=none](two11){1};
     \node at ($(one11)+(0.75,-1)$)[draw=none](three11){4};
     \path[thick](one11)edge(two11)edge(three11);
     \path[thick](root1)edge(one11);
     \path[densely dotted](root1)edge(one1);
     \node at ($(root1)+(4,0)$)[scale=1.2,fill=blue](onee1){};
     \node at ($(onee1)+(0.75,-1)$)[draw=none](twoo1){6};
     \node at ($(onee1)+(-0.75,-1)$)[inner sep=3pt,fill=red,circle](threee1){};
     \node at ($(threee1)+(-0.75,-1)$)[draw=none](fourr1){5};
     \node at ($(threee1)+(0.75,-1)$)[draw=none](fivee1){7};
     \path[thick](onee1)edge(twoo1)edge(threee1);
     \path[thick](threee1)edge(fourr1)edge(fivee1);
     \path[densely dotted](rootp)edge(onee1);
     \path[thick](rootp)edge(root1);
     \node at ($(o)+(10,-1)$)[draw=none](arrow1){$\downarrow$};
     \node at ($(arrow1)+(0,-1)$)[fill=blue,scale=1.2](rootp2){};
     \node at ($(arrow1)+(-2,-2.5)$)[inner sep=3pt,red,circle](root12){r};
     \node at ($(root12)+(4,0)$)[inner sep=3pt,fill=red,circle](one12){};
     \node at ($(one12)+(-0.75,-1)$)[draw=none](two12){2};
     \node at ($(one12)+(0.75,-1)$)[draw=none](three12){3};
     \path[thick](one12)edge(two12)edge(three12);
     \node at ($(root12)+(-1.1,-1)$)[inner sep=3pt,fill=red,circle](one112){};
     \node at ($(one112)+(-0.75,-1)$)[draw=none](two112){1};
     \node at ($(one112)+(0.75,-1)$)[draw=none](three112){4};
     \path[thick](one112)edge(two112);
     \path[densely dotted](one112)edge(three112);
     \node at ($(root12)+(1.1,-1)$)[scale=1.2,fill=blue](onee12){};
     \node at ($(onee12)+(0.75,-1)$)[draw=none](twoo12){6};
     \node at ($(onee12)+(-0.75,-1)$)[inner sep=3pt,fill=red,circle](threee12){};
      \path[densely dotted](root12)edge(onee12);
     \path[thick](root12)edge(one112);
     \node at ($(threee12)+(-0.75,-1)$)[draw=none](fourr12){5};
     \node at ($(threee12)+(0.75,-1)$)[draw=none](fivee12){7};
     \path[thick](onee12)edge(twoo12)edge(threee12);
     \path[thick](threee12)edge(fourr12)edge(fivee12);
     \path[thick](rootp2)edge(one12);
     \path[thick](rootp2)edge(root12);
    \filldraw[black] (16,-8)  node (5.6,-14.15)[anchor=west][draw=none,below](e12) {};
     \node at ($(e12)+(0,-1)$)[draw=none](arrow1){};
     \node at ($(o)+(5,-3)$)[draw=none](arrow1){$\longleftarrow$};
     \node at ($(o)+(0,-2)$)[fill=blue,scale=1.2](froot){};
     \node at ($(o)+(-2,-3.5)$)[inner sep=3pt,fill=red,circle](froot1){};
     \node at ($(froot1)+(4,0)$)[inner sep=3pt,fill=red,circle](fone1){};
     \node at ($(fone1)+(-0.75,-1)$)[draw=none](ftwo1){2};
     \node at ($(fone1)+(0.75,-1)$)[draw=none](fthree1){3};
     \path[thick](fone1)edge(ftwo1)edge(fthree1);
     \node at ($(froot1)+(1.1,-1)$)[draw=none](fone11){4};
     \node at ($(froot1)+(-1.1,-1)$)[inner sep=3pt,red,circle](f1){r};
     \node at ($(f1)+(-.75,-1)$)[draw=none](f11){1};
     \node at ($(f1)+(.75,-1)$)[scale=1.2,fill=blue](f12){};
     \node at ($(f12)+(-.75,-1)$)[inner sep=3pt,fill=red,circle](f21){};
     \node at ($(f21)+(-.75,-1)$)[draw=none](f31){5};
     \node at ($(f21)+(.75,-1)$)[draw=none](f32){7};
     \node at ($(f12)+(.75,-1)$)[draw=none](f22){6};
     \path[thick](froot)edge(froot1)edge(fone1);
     \path[thick](froot1)edge(fone11)edge(f1);
     \path[thick](f1)edge(f11)edge(f12);
     \path[thick](f12)edge(f21)edge(f22);
     \path[thick](f21)edge(f31)edge(f32);
     
     \end{tikzpicture}
     \caption{Pointed Lyndon tree obtained by a $1$-merge of two pointed Lyndon trees $T_1$ and $T_2$.}
     \label{fig:my_2-merge-pointed}
 \end{figure}
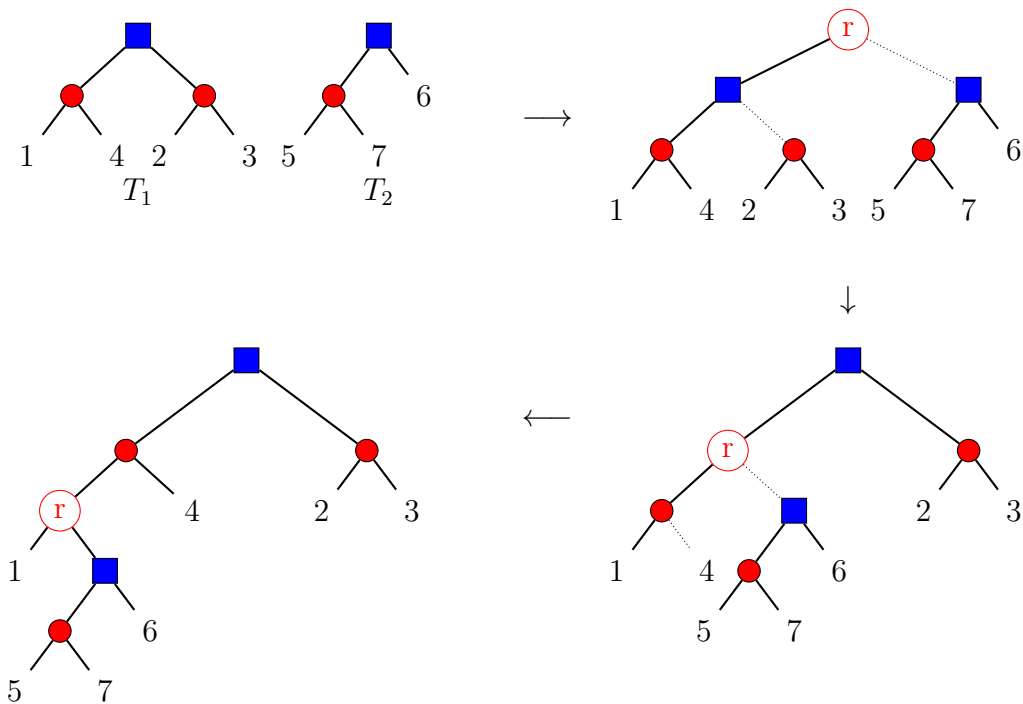
 
  An example of a $1$-merge of two pointed Lyndon trees is illustrated in Figure \ref{fig:my_2-merge-pointed}. First,  we add a vertex $r$ colored by \tre{1} (red) and add edges from $r$ to the roots of $T_1$ and $T_2$.  We then need to check whether or not this construction results into a pointed Lyndon tree. Since $\clr(L(r))=\tbl{0}<\tre{1}=\clr(r)$, we  have a violation of condition \eqref{condition:pointed_lyndon_1} at $r$ and hence, this is not yet a pointed Lyndon tree. We then slide $r$ together with its right subtree, to its left, interchanging $r$ and $L(r)$. After sliding $r$ to its left we have that $\clr(L(r))=\clr(r)=\tre{1}$. However,  $4=v(R(L(r)))<v(R(r))=5$ resulting in a violation of condition \eqref{condition:pointed_lyndon_2}, so we slide $r$ once more to its left to finally get a valid pointed Lyndon tree.

Let us remark here that this process always terminates in  a valid pointed Lyndon tree.  This is the case since if we keep sliding until we cannot anymore, the left child of the root will be a leaf. We are  now in a position to put an ordering on the pointed Lyndon forests.

 \begin{definition}   The \emph{poset of pointed Lyndon forests} is the set $\Flyn_{n}^{\bullet}$ together with the  cover relation $F\lessdot F'$   whenever $F'$ is obtained from $F$ when exactly two trees of $F$ are $u$-merged for some $u\in\left\{\tbl{0},\tre{1}\right\}$.
\end{definition}

    \begin{figure}[htbp]
    \begin{center}

    \resizebox{\columnwidth}{!}{%
\begin{tikzpicture}[scale=0.8,
rdot/.style = {fill=red, circle, draw, 
               minimum size=4mm, inner sep=0pt, outer sep=0pt,
               node contents={}},
sdot/.style  ={draw, fill=blue, scale=1.2,
               minimum size=2.7mm, inner sep=0pt, outer sep=0pt,
               node contents={}},
sibling distance = 7mm,scale=0.8
                    ]
\node (one) [rdot] 
    child{ node {$\tilde{1}$}} 
    child{ node {} edge from parent[draw=none]}
    child{ node [fill=red, circle,draw]{}
    child{node{2}}
    child{node (one south) {} edge from parent[draw=none]}
    child{node{3}}};
\begin{scope}[xshift=33mm]
\node (two) [sdot] 
    child{ node {1}} 
    child{ node {} edge from parent[draw=none]}
    child{ node [fill=red, circle,draw]{}
    child{node{$\tilde{2}$}}
    child{node(two south){} edge from parent[draw=none]}
    child{node{3}}};
\end{scope}
\begin{scope}[xshift=66mm]
\node (three) [sdot] 
    child{ node [draw,fill=blue,scale=1.2]{}
    child{node{1}}
    child{node(three south){} edge from parent[draw=none]}
    child{node{3}}}
    child{node{} edge from parent[draw=none]}
    child{node{$\tilde{2}$}};
\end{scope}
\begin{scope}[xshift=99mm]
\node (four) [sdot] 
    child{ node [fill=red, circle,draw]{}
    child{node{1}}
    child{node(four south){} edge from parent[draw=none]}
    child{node{2}}}
    child{node{} edge from parent[draw=none]}
    child{node{$\tilde{3}$}};
\end{scope}
\begin{scope}[xshift=132mm]
\node (five) [rdot] 
    child{ node [fill=red, circle,draw]{}
    child{node{$\tilde{1}$}}
    child{node(five south){} edge from parent[draw=none]}
    child{node{3}}}
    child{node{} edge from parent[draw=none]}
    child{node{2}};
\end{scope}
\begin{scope}[xshift=165mm]
\node (six) [sdot] 
    child{ node [fill=red, circle,draw]{}
    child{node{1}}
    child{node(six south){} edge from parent[draw=none]}
    child{node{3}}}
    child{node{} edge from parent[draw=none]}
    child{node{$\tilde{2}$}};
\end{scope}
\begin{scope}[xshift=198mm]
\node (seven) [sdot] 
    child{ node [draw,fill=blue,scale=1.2]{}
    child{node{1}}
    child{node(seven south){} edge from parent[draw=none]}
    child{node{2}}}
    child{node{} edge from parent[draw=none]}
    child{node{$\tilde{3}$}};
\end{scope}
\begin{scope}[xshift=231mm]
\node (eight) [rdot] 
    child{ node {$\tilde{1}$}} 
    child{ node {} edge from parent[draw=none]}
    child{ node [draw,fill=blue,scale=1.2]{}
    child{node{2}}
    child{node(eight south){} edge from parent[draw=none]}
    child{node{3}}};
\end{scope}
\begin{scope}[xshift=264mm]
\node (nine) [sdot] 
    child{ node {1}} 
    child{ node {} edge from parent[draw=none]}
    child{ node [draw,fill=blue,scale=1.2]{}
    child{node{2}}
    child{node(nine south){} edge from parent[draw=none]}
    child{node{$\tilde{3}$}}};
\end{scope}
\begin{scope}[shift={(49.5mm,-85mm)}]
\node (one1) [rdot] 
    child{ node {$\tilde{2}$}} 
    child{ node (one1 south){} edge from parent[draw=none]}
    child{ node {3}};
    \filldraw[black] (1.2,0) circle (1pt) node[anchor=west][draw=none,below] {$\tilde{1}$};
\end{scope}
\begin{scope}[shift={(82.5mm,-85mm)}]
\node (two1) [sdot] 
    child{ node {1}} 
    child{ node (two1 south){} edge from parent[draw=none]}
    child{ node {$\tilde{3}$}};
    \filldraw[black] (1.2,0) circle (1pt) node[anchor=west][draw=none,below] {$\tilde{2}$};
\end{scope}
\begin{scope}[shift={(115.5mm,-85mm)}]
\node (three1) [rdot] 
    child{ node {$\tilde{1}$}} 
    child{ node (three1 south){} edge from parent[draw=none]}
    child{ node {2}};
    \filldraw[black] (1.2,0) circle (1pt) node[anchor=west][draw=none,below] {$\tilde{3}$};
\end{scope}
\begin{scope}[shift={(148.5mm,-85mm)}]
\node (four1) [rdot] 
    child{ node {$\tilde{1}$}} 
    child{ node (four1 south){} edge from parent[draw=none]}
    child{ node {3}};
    \filldraw[black] (1.2,0) circle (1pt) node[anchor=west][draw=none,below] {$\tilde{2}$};
\end{scope}
\begin{scope}[shift={(181.5mm,-85mm)}]
\node (five1) [sdot] 
    child{ node {1}} 
    child{ node (five1 south) {} edge from parent[draw=none]}
    child{ node {$\tilde{2}$}};
    \filldraw[black] (1.2,0) circle (1pt) node[anchor=west][draw=none,below] {$\tilde{3}$};
\end{scope}
\begin{scope}[shift={(214.5mm,-85mm)}]
\node (six1) [sdot] 
    child{ node {2}} 
    child{ node(six1 south) {} edge from parent[draw=none]}
    child{ node {$\tilde{3}$}};
    \filldraw[black] (1.2,0) circle (1pt) node[anchor=west][draw=none,below] {$\tilde{1}$};
\end{scope}
\begin{scope}[shift={(132mm,-160mm)}]
\node(et) {$\dot{\tilde{1}}\; \dot{\tilde{2}}\; \dot{\tilde{3}}$};
\end{scope}

\draw[very thick]   ([yshift=3mm] one1.north) -- (one south.south)
([yshift=3mm] one1.north) -- (two south.south)
([yshift=3mm] two1.north) -- (three south.south)
([yshift=3mm] two1.north) -- (four south.south)
([yshift=3mm] three1.north) -- (four south.south)
([yshift=3mm] three1.north) -- (five south.south)
([yshift=3mm] four1.north) -- (five south.south)
([yshift=3mm] four1.north) -- (six south.south)
([yshift=3mm] five1.north) -- (six south.south)
([yshift=3mm] five1.north) -- (seven south.south)
([yshift=3mm] six1.north) -- (eight south.south)
([yshift=3mm] six1.north) -- (nine south.south)
([yshift=3mm] et.north) -- (one1 south.south)
([yshift=3mm] et.north) -- (two1 south.south)
([yshift=3mm] et.north) -- (three1 south.south)
([yshift=3mm] et.north) -- (four1 south.south)
([yshift=3mm] et.north) -- (five1 south.south)
([yshift=3mm] et.north) -- (six1 south.south);
    \end{tikzpicture}
    }
    \end{center}
\caption{$\Flyn_{3}^{\bullet}$.}
\label{fig:FLyn_3-pointed}

    \end{figure}

We illustrate $\Flyn_{3}^{\bullet}$ in Figure \ref{fig:FLyn_3-pointed}.    The sliding procedure described to merge two pointed Lyndon trees is just a way to explain the sorting procedure used to define the Whitney dual $R_\lambda$ in Definition~~\ref{definition:R_lambda}.   Note that in the definition of the map sending $F$ to $c(F)$ we did not need $F$ to be a pointed Lyndon forest nor do we need to use the reverse-minimal ordering on the internal vertices. Indeed as long as $F$ is a normalized bicolored binary forest and we use some linear extension of the internal vetices, the map still produces a saturated chain in $\Pi_n^\bullet$ starting at $\hat{0}$. In the case that the corresponding chain is not necessarily ascent-free, swapping the labels in an ascent either corresponds to reordering the internal vertices to get the reverse-minimal order or corresponds to the sliding procedure.  For example, say we have the following label sequence in $\Pi_4^\bullet$
$$
(1,2)^{\tre{1}}(3,4)^{\tre{1}}
$$
This corresponds to a pointed Lyndon forest with two components where the internal vertices are ordered so that the internal vertex above leafs 1 and 2 comes first.  When we swap the labels in the sequence to get the ascent-free sequence  
$$
(3,4)^{\tre{1}}(1,2)^{\tre{1}}
$$
we are just reordering the internal vertices so the one above leafs 3 and 4  comes first. 

For an example where sliding occurs, consider the sliding procedure shown in Figure~\ref{fig:my_2-merge-pointed}.  The sequence of labels of the pointed Lyndon  forest in the upper left corner of the figure is 
$$
(5,7)^{\tre{1}}(5,6)^{\tbl{0}}(2,3)^{\tre{1}}(1,4)^{\tre{1}}(1,2)^{\tbl{0}}
$$
Adding the red vertex $r$ is corresponds to adding the the label $(1,5)^{\tre{1}}$ at the end of the sequence since it merges together the components with minimal leaf label $1$ and $5$.  So we would have the sequence
$$
(5,7)^{\tre{1}}(5,6)^{\tbl{0}}(2,3)^{\tre{1}}(1,4)^{\tre{1}}(1,2)^{\tbl{0}}(1,5)^{\tre{1}}.
$$
This is the label sequence for a saturated chain starting at $\hat{0}$ of $\Pi_n^\bullet$.  However, it is not ascent-free since $(1,2)^{\tbl{0}}(1,5)^{\tre{1}}$ is an ascent. Moreover, the corresponding tree in the upper right corner of Figure~\ref{fig:my_2-merge-pointed} is not a pointed Lyndon tree. As $\lambda_\bullet$ has the rank two switching property, we can swap these to labels to get the sequence
$$
(5,7)^{\tre{1}}(5,6)^{\tbl{0}}(2,3)^{\tre{1}}(1,4)^{\tre{1}}(1,5)^{\tre{1}}(1,2)^{\tbl{0}}.
$$
This swap corresponds to sliding the root $r$ to the left of its left child giving the tree in the bottom right corner of Figure~\ref{fig:my_2-merge-pointed} whose label sequence is the one given above.  This sequence has an ascent at $(1,4)^{\tre{1}}(1,5)^{\tre{1}}$ and the corresponding tree is not a pointed Lyndon tree. Swapping labels, we  get
$$
(5,7)^{\tre{1}}(5,6)^{\tbl{0}}(2,3)^{\tre{1}}(1,5)^{\tre{1}}(1,4)^{\tre{1}}(1,2)^{\tbl{0}}
$$
which is ascent-free. Again, this swap corresponds to sliding $r$ once again past its left child to get the tree in the bottom left corner to Figure~\ref{fig:my_2-merge-pointed}.  Note that this sequence is ascent-free and the tree we finish with is the corresponding pointed Lyndon tree for this sequence.

Because the sliding procedure corresponds to the sorting procedure in the Whitney dual, we get the following.

\begin{theorem} \label{theorem:whitney_dual_pointed} \label{theorem:flynpointed_dual}$\Flyn_{n}^{\bullet}$ is a Whitney dual to $\Pi_n^\bullet$. In particular, $\Flyn_{n}^{\bullet}\cong R_{\lambda_\bullet}(\Pi_n^\bullet)$.
\end{theorem}

We omit the full details of this proof here since is rather technical and a case by case analysis of the ways one can have ascents in the chains.     For all the details, the reader can consult ~\cite[Section 3.2.2]{QuicenoDuran2020}

\subsection{A Whitney dual for \texorpdfstring{$\Pi_n^w$}-}\label{section:WhitneyDualOfWeighted}  Here we give a combinatorial description for the Whitney dual of the weighted partition poset.  The method closely follows what we did in the previous subsection. As such, we do not provide as many detailed examples.

In~\cite[Theorem 5.7]{DleonWachs2016}, it was shown that the maximal ascent-free saturated chains of $\Pi_n^w$   with respect to $\lambda_w$ are in bijection with bicolored Lyndon trees.   This bijection can be modified  to give a bijection between ascent-free chains of $\Pi_n^w$ starting at $\hat{0}$ and bicolored Lyndon forests.  It follows that the elements of the Whitney dual $R_{\lambda_w}(\Pi_n^w)$ can be described using bicolored Lyndon forests.    We briefly explain this bijection.

Let $F$ be a bicolored Lyndon forest. As mentioned in Section \ref{section:combinatorial_families},  $F$ has a unique reverse-minimal linear extension on the internal vertices.  To construct our ascent-free chain $c(F)$, we follow a similar procedure as in the case for pointed Lyndon forests (see Figure ~\ref{fig:pointed_Lyndon_Tree_and_chain}). We start with the the bottom element, $1^0/2^0/\cdots /n^0$.  At the $i^{th}$ step we merge together the elements of the left and right subtrees rooted at the $i^{th}$ internal vertex. We add together the weights of these blocks and add 1 to this weight if this internal vertex is colored by \tre{1}. If the internal vertex is colored by \tbl{0}, we do not add 1 to this sum. See \cite[Figure 4]{DleonWachs2016} for an example.

 \begin{figure}
     \centering
     
\begin{tikzpicture}[scale=0.8]

\tikzstyle{every node} = [draw=black]
     \node[inner sep=3pt,fill=red,circle](root){};
     \node at ($(root)+(1.1,-1)$)[fill=blue,scale=1.2](one){};
     \node at ($(one)+(-0.75,-1)$)[draw=none](two){2};
     \node at ($(one)+(0.75,-1)$)[draw=none](three){3};
     \path[thick](one)edge(two)edge(three);
     \node at ($(root)+(-1.1,-1)$)[fill,blue,scale=1.2](one1){};
     \node at ($(one1)+(-0.75,-1)$)[draw=none](two1){1};
     \node at ($(one1)+(0.75,-1)$)[draw=none](three1){4};
     \path[thick](one1)edge(two1)edge(three1);
     \path[thick](root)edge(one1)edge(one);
     \node at ($(root)+(0,-2.6)$)[draw=none](arrow1){$T_1$};
     \node at ($(root)+(4,0)$)[scale=1.2,fill=blue](onee){};
     \node at ($(onee)+(0.75,-1)$)[draw=none](twoo){6};
     \node at ($(onee)+(-0.75,-1)$)[inner sep=3pt,fill=red,circle](threee){};
     \node at ($(threee)+(-0.75,-1)$)[draw=none](fourr){5};
     \node at ($(threee)+(0.75,-1)$)[draw=none](fivee){7};
     \path[thick](onee)edge(twoo)edge(threee);
     \path[thick](threee)edge(fourr)edge(fivee);
     \node at ($(onee)+(0,-2.6)$)[draw=none](arrow1){$T_2$};
     \node at ($(one)+(0.7,-2.4)$)[draw=none](o){};
     \node at ($(o)+(5,2)$)[draw=none](arrow1){$\longrightarrow$};
     \node at ($(o)+(10,3.5)$)[inner sep=3pt,red,circle](rootp){r};
     \node at ($(o)+(8,2.5)$)[inner sep=3pt,fill=red,circle](root1){};
     \node at ($(root1)+(1.1,-1)$)[fill,blue,scale=1.2](one1){};
     \node at ($(one1)+(-0.75,-1)$)[draw=none](two1){2};
     \node at ($(one1)+(0.75,-1)$)[draw=none](three1){3};
     \path[thick](one1)edge(two1)edge(three1);
     \node at ($(root1)+(-1.1,-1)$)[scale=1.2,fill=blue](one11){};
     \node at ($(one11)+(-0.75,-1)$)[draw=none](two11){1};
     \node at ($(one11)+(0.75,-1)$)[draw=none](three11){4};
     \path[thick](one11)edge(two11)edge(three11);
     \path[thick](root1)edge(one11);
     \path[densely dotted](root1)edge(one1);
     \node at ($(root1)+(4,0)$)[scale=1.2,fill=blue](onee1){};
     \node at ($(onee1)+(0.75,-1)$)[draw=none](twoo1){6};
     \node at ($(onee1)+(-0.75,-1)$)[inner sep=3pt,fill=red,circle](threee1){};
     \node at ($(threee1)+(-0.75,-1)$)[draw=none](fourr1){5};
     \node at ($(threee1)+(0.75,-1)$)[draw=none](fivee1){7};
     \path[thick](onee1)edge(twoo1)edge(threee1);
     \path[thick](threee1)edge(fourr1)edge(fivee1);
     \path[densely dotted](rootp)edge(onee1);
     \path[thick](rootp)edge(root1);
     \node at ($(o)+(10,-1)$)[draw=none](arrow1){$\downarrow$};
     \node at ($(arrow1)+(0,-1)$)[inner sep=3pt,fill=red,circle](rootp2){};
     \node at ($(arrow1)+(-2,-2.5)$)[inner sep=3pt,red,circle](root12){r};
     \node at ($(root12)+(4,0)$)[scale=1.2,fill=blue](one12){};
     \node at ($(one12)+(-0.75,-1)$)[draw=none](two12){2};
     \node at ($(one12)+(0.75,-1)$)[draw=none](three12){3};
     \path[thick](one12)edge(two12)edge(three12);
     \node at ($(root12)+(-1.1,-1)$)[scale=1.2,fill=blue](one112){};
     \node at ($(one112)+(-0.75,-1)$)[draw=none](two112){1};
     \node at ($(one112)+(0.75,-1)$)[draw=none](three112){4};
     \path[thick](one112)edge(two112);
     \path[densely dotted](one112)edge(three112);
     \node at ($(root12)+(1.1,-1)$)[scale=1.2,fill=blue](onee12){};
     \node at ($(onee12)+(0.75,-1)$)[draw=none](twoo12){6};
     \node at ($(onee12)+(-0.75,-1)$)[inner sep=3pt,fill=red,circle](threee12){};
      \path[densely dotted](root12)edge(onee12);
     \path[thick](root12)edge(one112);
     \node at ($(threee12)+(-0.75,-1)$)[draw=none](fourr12){5};
     \node at ($(threee12)+(0.75,-1)$)[draw=none](fivee12){7};
     \path[thick](onee12)edge(twoo12)edge(threee12);
     \path[thick](threee12)edge(fourr12)edge(fivee12);
     \path[thick](rootp2)edge(one12);
     \path[thick](rootp2)edge(root12);
    \filldraw[black] (16,-8)  node (5.6,-14.15)[anchor=west][draw=none,below](e12) {};
     \node at ($(e12)+(0,-1)$)[draw=none](arrow1){};

     \node at ($(o)+(5,-3)$)[draw=none](arrow1){$\longleftarrow$};
     \node at ($(o)+(0,-2)$)[inner sep=3pt,fill=red,circle](froot){};
     \node at ($(o)+(-2,-3.5)$)[scale=1.2,fill=blue](froot1){};
     \node at ($(froot1)+(4,0)$)[scale=1.2,fill=blue](fone1){};
     \node at ($(fone1)+(-0.75,-1)$)[draw=none](ftwo1){2};
     \node at ($(fone1)+(0.75,-1)$)[draw=none](fthree1){3};
     \path[thick](fone1)edge(ftwo1)edge(fthree1);
     \node at ($(froot1)+(1.1,-1)$)[draw=none](fone11){4};
     \node at ($(froot1)+(-1.1,-1)$)[inner sep=3pt,red,circle](f1){r};
     \node at ($(f1)+(-.75,-1)$)[draw=none](f11){1};
     \node at ($(f1)+(.75,-1)$)[scale=1.2,fill=blue](f12){};
     \node at ($(f12)+(-.75,-1)$)[inner sep=3pt,fill=red,circle](f21){};
     \node at ($(f21)+(-.75,-1)$)[draw=none](f31){5};
     \node at ($(f21)+(.75,-1)$)[draw=none](f32){7};
     \node at ($(f12)+(.75,-1)$)[draw=none](f22){6};
     \path[thick](froot)edge(froot1)edge(fone1);
     \path[thick](froot1)edge(fone11)edge(f1);
     \path[thick](f1)edge(f11)edge(f12);
     \path[thick](f12)edge(f21)edge(f22);
     \path[thick](f21)edge(f31)edge(f32);

     \end{tikzpicture}
     \caption{Bicolored Lyndon tree obtained by a $1$-merge of two bicolored Lyndon trees $T_1$ and $T_2$.}
     \label{fig:my_2-merge-bicolored}
 \end{figure}
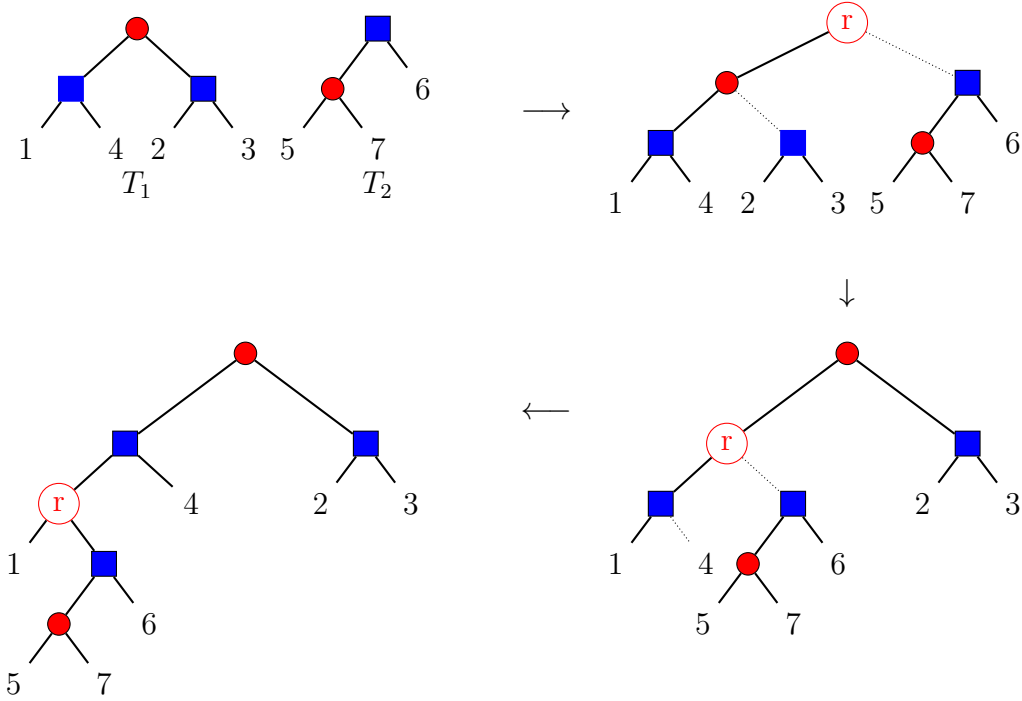
 
 To describe the covering relation in $R_\lambda(\Pi_n^w)$, we need to describe a method to merge bicolored Lyndon trees.  Just like the case for the pointed partition poset, the covering relation is defined using a sliding process.  Let $T_1$ and $T_2$ be bicolored Lyndon trees with the minimum leaf label of $T_1$ less than the minimum leaf label of $T_2$. Let $u\in \{\tbl{0},\tre{1}\}$. To \emph{$u$-merge} $T_1$ and $T_2$, we first create a new vertex $r$ and color it so that $\clr{(r)}=u$.  Then we add edges from $r$ to the roots of $T_1$ and $T_2$.  If the resulting tree is a bicolored Lyndon tree, we stop. If not, we slide the new internal vertex $r$ together with its left subtree past its left child. We continue this procedure until we get a bicolored Lyndon tree.  See Figure~\ref{fig:my_2-merge-bicolored} for an example of this procedure.

 \begin{definition}   The \emph{poset of bicolored Lyndon forests} is the set of bicolored Lyndon forests on $[n]$, $\Flyn_{n}^{w}$, together with the  cover relation $F\lessdot F'$   whenever $F'$ is obtained from $F$ when exactly two trees of $F$ are $u$-merged for some $u\in\left\{\tbl{0},\tre{1}\right\}$.
\end{definition}

 \begin{figure}[htbp]
    \begin{center}
          \resizebox{\columnwidth}{!}{%

\begin{tikzpicture}[scale=0.8,
rdot/.style = {circle, draw, fill=red, 
               minimum size=4mm, inner sep=0pt, outer sep=0pt,
               node contents={}},
sdot/.style  ={draw, fill=blue, scale=1.2,
               minimum size=2.7mm, inner sep=0pt, outer sep=0pt,
               node contents={}},
sibling distance = 7mm,scale=0.8
                    ]
\node (one) [rdot] 
    child{ node {1}} 
    child{ node {} edge from parent[draw=none]}
    child{ node [circle,draw,fill=red]{}
    child{node{2}}
    child{node (one south) {} edge from parent[draw=none]}
    child{node{3}}};
\begin{scope}[xshift=33mm]
\node (two) [sdot] 
    child{ node {1}} 
    child{ node {} edge from parent[draw=none]}
    child{ node [circle,draw,fill=red]{}
    child{node{2}}
    child{node(two south){} edge from parent[draw=none]}
    child{node{3}}};
\end{scope}
\begin{scope}[xshift=66mm]
\node (three) [sdot] 
    child{ node [circle,draw,fill=red]{}
    child{node{1}}
    child{node(three south){} edge from parent[draw=none]}
    child{node{2}}}
    child{node{} edge from parent[draw=none]}
    child{node{3}};
\end{scope}
\begin{scope}[xshift=99mm]
\node (four) [rdot] 
    child{ node [circle,draw,fill=red]{}
    child{node{1}}
    child{node(four south){} edge from parent[draw=none]}
    child{node{3}}}
    child{node{} edge from parent[draw=none]}
    child{node{2}};
\end{scope}
\begin{scope}[xshift=132mm]
\node (five) [sdot] 
    child{ node [rdot]{}
    child{node{1}}
    child{node(five south){} edge from parent[draw=none]}
    child{node{3}}}
    child{node{} edge from parent[draw=none]}
    child{node{2}};
\end{scope}
\begin{scope}[xshift=165mm]
\node (six) [sdot] 
    child{ node [draw,fill=blue,scale=1.2]{}
    child{node{1}}
    child{node(six south){} edge from parent[draw=none]}
    child{node{3}}}
    child{node{} edge from parent[draw=none]}
    child{node{2}};
\end{scope}
\begin{scope}[xshift=198mm]
\node (seven) [rdot] 
    child{ node [draw,fill=blue,scale=1.2]{}
    child{node{1}}
    child{node(seven south){} edge from parent[draw=none]}
    child{node{3}}}
    child{node{} edge from parent[draw=none]}
    child{node{2}};
\end{scope}
\begin{scope}[xshift=231mm]
\node (eight) [rdot] 
    child{ node {1}} 
    child{ node {} edge from parent[draw=none]}
    child{ node [draw,fill=blue,scale=1.2]{}
    child{node{2}}
    child{node(eight south){} edge from parent[draw=none]}
    child{node{3}}};
\end{scope}
\begin{scope}[xshift=264mm]
\node (nine) [sdot] 
    child{ node {1}} 
    child{ node {} edge from parent[draw=none]}
    child{ node [draw,fill=blue,scale=1.2]{}
    child{node{2}}
    child{node(nine south){} edge from parent[draw=none]}
    child{node{3}}};
\end{scope}
\begin{scope}[shift={(49.5mm,-85mm)}]
\node (one1) [rdot] 
    child{ node {2}} 
    child{ node (one1 south){} edge from parent[draw=none]}
    child{ node {3}};
    \filldraw[black] (1.2,0) circle (1pt) node[anchor=west][draw=none,below] {1};
\end{scope}
\begin{scope}[shift={(82.5mm,-85mm)}]
\node (two1) [rdot] 
    child{ node {1}} 
    child{ node (two1 south){} edge from parent[draw=none]}
    child{ node {2}};
    \filldraw[black] (1.2,0) circle (1pt) node[anchor=west][draw=none,below] {3};
\end{scope}
\begin{scope}[shift={(115.5mm,-85mm)}]
\node (three1) [rdot] 
    child{ node {1}} 
    child{ node (three1 south){} edge from parent[draw=none]}
    child{ node {3}};
    \filldraw[black] (1.2,0) circle (1pt) node[anchor=west][draw=none,below] {2};
\end{scope}
\begin{scope}[shift={(148.5mm,-85mm)}]
\node (four1) [sdot] 
    child{ node {1}} 
    child{ node (four1 south){} edge from parent[draw=none]}
    child{ node {2}};
    \filldraw[black] (1.2,0) circle (1pt) node[anchor=west][draw=none,below] {3};
\end{scope}
\begin{scope}[shift={(181.5mm,-85mm)}]
\node (five1) [sdot] 
    child{ node {1}} 
    child{ node (five1 south) {} edge from parent[draw=none]}
    child{ node {3}};
    \filldraw[black] (1.2,0) circle (1pt) node[anchor=west][draw=none,below] {2};
\end{scope}
\begin{scope}[shift={(214.5mm,-85mm)}]
\node (six1) [sdot] 
    child{ node {2}} 
    child{ node(six1 south) {} edge from parent[draw=none]}
    child{ node {3}};
    \filldraw[black] (1.2,0) circle (1pt) node[anchor=west][draw=none,below] {1};
\end{scope}
\begin{scope}[shift={(132mm,-160mm)}]
\node(et) {$\dot{1}\; \dot{2}\; \dot{3}$};
\end{scope}

\draw[very thick]   ([yshift=3mm] one1.north) -- (one south.south)
([yshift=3mm] one1.north) -- (two south.south)
([yshift=3mm] two1.north) -- (three south.south)
([yshift=3mm] two1.north) -- (four south.south)
([yshift=3mm] three1.north) -- (four south.south)
([yshift=3mm] three1.north) -- (five south.south)
([yshift=3mm] four1.north) -- (five south.south)
([yshift=3mm] four1.north) -- (six south.south)
([yshift=3mm] five1.north) -- (six south.south)
([yshift=3mm] five1.north) -- (seven south.south)
([yshift=3mm] six1.north) -- (eight south.south)
([yshift=3mm] six1.north) -- (nine south.south)
([yshift=3mm] et.north) -- (one1 south.south)
([yshift=3mm] et.north) -- (two1 south.south)
([yshift=3mm] et.north) -- (three1 south.south)
([yshift=3mm] et.north) -- (four1 south.south)
([yshift=3mm] et.north) -- (five1 south.south)
([yshift=3mm] et.north) -- (six1 south.south);
    \end{tikzpicture}
}
\caption{$\Flyn^w_{3}$.}
\label{fig:FLyn_3_weighted}

\end{center}
    \end{figure}

Figure~\ref{fig:FLyn_3_weighted} depicts $\Flyn_{3}^{w}$.

\begin{theorem} \label{theorem:flynweighted_dual}  $\Flyn_{n}^{w}$ is a Whitney dual of $\Pi_n^w$. In particular, $\Flyn_{n}^{w}\cong R_{\lambda_w}(\Pi_n^w)$.
\end{theorem}

As is the case with the pointed partition poset, the proof is a case by case analysis of the ways one might have an ascent on an ascent-free chain when we add a new label. See~\cite[Section 3.1.2]{QuicenoDuran2020} for all the details.

\begin{remark}
In $\Pi_n^w$, the interval $[\hat{0}, [n]^0]$ is isomorphic to the partition lattice $\Pi_n$. The labeling $\lambda_w$ restricted to this interval is also an EW-labeling, and hence a subposet of $\Flyn_{n}^{w}$ is also a Whitney dual for the partition lattice.  In $[\hat{0}, [n]^0]$ all labels are of the form $(a,b)^0$. Thus, all the ascent-free chains correspond to  bicolored Lyndon forests where every vertex must be Lyndon.  Moreover, all the internal vertices are colored $\tbl{0}$ (blue), so we can just assume that the internal vertices are uncolored. We then get a Whitney dual of $\Pi_n$ where the elements are normalized binary forests whose internal vertices are Lyndon and where the covering relation is given by merging together trees to get normalized binary trees whose vertices are Lyndon. By restricting to forests with internal vertices colored $\tbl{0}$ (blue) in Figure~\ref{fig:FLyn_3_weighted}, we get a depiction of this Whitney dual of $\Pi_3$.  We should note that restricting $\lambda_w$ to $[\hat{0}, [n]^0]$ yields  Stanley's labeling ~\cite{Stanley1974} for $\Pi_n$.  Stanley's labeling  was used in~\cite[Corollary 5.7]{GonzalezHallam2021}  to construct a Whitney dual to $\Pi_n$ isomorphic to the poset of increasing spanning forests $\ISF_n$.  It follows then that the subposet of $\Flyn_{n}^{w}$ formed by Lyndon forests (with blue internal nodes) is isomorphic to $\ISF_n$.

\end{remark}

\section{Algebraic consequences}\label{sec:algConseq}

\subsection{Homological consequences of the EL-labelings}
     In Section~\ref{sec:EWLabOfPointPart} we showed that $\lambda_{\bullet}$  is an EL-labeling of the order dual of $\Pi_n^\bullet$ and $\lambda_{\bullet_2}$ is an EL-labeling of $\Pi_n^\bullet$. A poset and its order dual have the same order complex and hence also  have the same cohomology. As a result, Theorem \ref{theorem:bjorner_wachs} implies that the two EL-labelings give bases for the cohomology of maximal intervals of $\Pi_n^\bullet$.  These bases are indexed by the ascent-free chains in the two labelings.

     In Section~\ref{subsection:whitney_dual_pointed} we showed that the maximal ascent-free chains of   $\Pi_n^\bullet$ with respect to $\lambda_{\bullet}$ are indexed by pointed Lyndon trees. To each pointed Lyndon tree $T$, we gave a bijection mapping $T$ to the ascent-free chain $c(T)$. 
     
     Let  $\mathcal{T}Lyn_{n,p}^{\bullet}$ be the set of pointed Lyndon trees such that along the path from the leaf labeled $p$ to the root, if the path moves to the left, the internal vertex label is $0$ and if it is to the right, the label is given by $1$. In Figure \ref{fig:FLyn_3-pointed} the reader can easily observe all the trees in $\mathcal{T}Lyn_{3,p}^{\bullet}$ for  $p=1,2,3$ by selecting the subset of maximal elements in $\Flyn_3^{\bullet}$ with the leaf $p$ decorated with a tilde ($\tilde{p}$). Note that $T\in  \mathcal{T}Lyn_{n,p}^{\bullet}$, if and only if the top element of $c(T)$ is $[n]^p$. We have the following consequences of Theorem \ref{theorem:pointed_dual_EL_labeling} and the characterization of the ascent-free chains of $\lambda_{\bullet}$ presented in Section~\ref{subsection:whitney_dual_pointed}.

\begin{theorem}\label{theorem:consequence_first_pointed_labeling} 
  For every $p\in [n]$ we have that
 \begin{enumerate}
   \item   The order complex $\Delta((\hat{0},[n]^p))$ is shellable and has the homotopy type of a wedge of $| \mathcal{T}Lyn_{n,p}^{\bullet}|$ many spheres of dimension $n-3$. As a consequence, the interval $[\hat{0},[n]^p]$ is Cohen-Macaulay.

\item The set
$ \{c(T) \mid T\in \mathcal{T}Lyn_{n,p}^{\bullet}\}$ forms a basis for $\widetilde H^{n-3}(\hat{0},[n]^p)$.
 \end{enumerate}
\end{theorem}

\begin{remark}
    Note that as a consequence of Theorem \ref{theorem:consequence_first_pointed_labeling} and the fact that the intervals $[\hat{0},[n]^p]$ are  isomorphic,   all the sets  $\mathcal{T}Lyn_{n,p}^{\bullet}$ are  equinumerous for any given $p\in [n]$.
\end{remark}

Now let us turn our attention to the consequences of the EL-labelings of $\Pi_n^\bullet$ for the operad $\Prelie$. In the theory of nonassociative algebras, a $\Prelie$-algebra is a vector space $V$ that comes equipped with a binary operation $\circ$ which satisfies for every $v,w,z \in V$ the relation
 \begin{align*}
     (v \circ w) \circ z - v \circ (w \circ z) =(v \circ z) \circ w - v \circ (z \circ w).
 \end{align*}
 Let $\Prelie(n)$ be the multilinear component of the free $\Prelie$ algebra on $n$ generators. In~\cite{Vallette2007}, Vallette 
 proved the following theorem (in terms of homology, which we reinterpret here in terms of cohomology).

 \begin{theorem}[Theorem 13 \cite{Vallette2007}]\label{theorem:vallette_prelie}
 We have the following $\sym_n$-module isomorphism
$$\Prelie(n)\cong_{\sym_n}\bigoplus_{p\in[n]}\widetilde H^{n-3}(\hat{0},[n]^p)\otimes \sgn_{n},$$
 where $\sgn_n$ is the sign representation of $\sym_n$.
 \end{theorem}
 
 Under Theorem \ref{theorem:vallette_prelie} we obtain a corresponding basis for $\Prelie(n)$, which we describe now. Let  $T= T_L \wedge^u T_R$ denote a normalized bicolored binary tree where $T_L$ and $T_R$ are respectively the left and right subtrees from the root and $u$ is the color of the root.
 Define $\Theta(T)$ to be the element in $\Prelie(n)$ defined recursively by $\Theta(T)=a$ when $T=a$ is the one-leaved tree with leaf-label $a$, and if $T=T_L \wedge^u T_R$ then
 \begin{align*}
     \Theta(T)=\begin{cases}
         \Theta(T_L)\circ \Theta(T_R) \text{ if } u=1,\\
         \Theta(T_R)\circ \Theta(T_L) \text{ if } u=0.
     \end{cases}
 \end{align*}

     As an example of this definition, let $T$ be the  pointed Lyndon tree in the bottom left of Figure~\ref{fig:my_2-merge-pointed}.  One can check that the associated monomial $\Theta(T)$ is $(2\circ 3)\circ ((1\circ (6 \circ (5\circ 7))) \circ 4)$. Theorems \ref{theorem:consequence_first_pointed_labeling} and \ref{theorem:vallette_prelie} imply the following theorem.

    \begin{theorem}\label{theorem:pointed_lyndon_basis_prelie}
The set
$\{\Theta(T) \mid T\in \mathcal{T}Lyn_{n}^{\bullet}\}$ forms a basis for $\mathcal{P}re\mathcal{L}ie(n)$.
    \end{theorem}

    In \cite{DleonWachs2016}, the authors proved the analogous theorem to Theorem \ref{theorem:consequence_first_pointed_labeling} providing a basis for the reduced cohomology $\widetilde H^{n-3}(\hat{0},[n]^i)$ of the maximal intervals of $\Pi_n^w$ for $i=0,\dots, n-1$, and for the multilinear component $\Lie^2(n)$ of the free bibracketed Lie algebra in $n$ generators. Those bases are indexed in terms of the bicolored Lyndon trees, $\mathcal{T}Lyn^{w}_{n}$,  since they index the ascent-free chains of $\lambda_{w}$. We will show that the same set of trees index a basis for $\widetilde H^{n-3}(\hat{0},[n]^p)$ and  $\Prelie(n)$.

    In \cite{DleonWachs2016} the authors prove that there is a rank-preserving bijection between $\Pi_n^{\bullet}$ and $\Pi_n^{w}$. We prove here the following further statement about their sets of saturated chains from $\hat{0}$.

    \begin{theorem}\label{Theorem:label-preserving_bijection}
        There is a label-preserving bijection between saturated chains from $\hat{0}$ in $(\Pi_n^{\bullet},\lambda_{\bullet})$ (or $(\Pi_n^{\bullet},\lambda_{\bullet_2})$) and in $(\Pi_n^{w},\lambda_{w})$.
    \end{theorem}
    \begin{proof}
    First note that between $\lambda_{\bullet}$ and $\lambda_{\bullet_2}$, the edge labels are the same, only the ordering of the labels is different. So we can find the bijection using $\lambda_{\bullet}$. 
        In $(\Pi_n^{\bullet},\lambda_{\bullet})$ at every step on a saturated chain from $\hat{0}$ we  $u$-merge two blocks $(A,p)$ and $(B,q)$ such that $\min A < \min B$ and assign the label $$\lambda_{\bullet_2}(\bpi\lessdot \bpi')=(\min A, \min B)^u.$$ 
        In $(\Pi_n^{w},\lambda_{w})$ at every step we  merge two blocks $(A,i)$ and $(B,j)$ to obtain the block $(A\cup B, i+j+u)$ and assign the label $$\lambda_{w}(\bpi\lessdot \bpi')=(\min A, \min B)^u.$$
         Note that in both cases, at every merging step from bottom to top we are free to choose between $u=0$ or $u=1$ and hence the sets of words of labels for saturated chains from $\hat{0}$ are equal. Since the saturated chains are uniquely determined by their words of labels in both labelings, the words of labels induce a bijection among saturated chains.
    \end{proof}

     Theorem \ref{theorem:EL-labeling_pointed} gives us analogous results to  Theorem \ref{theorem:consequence_first_pointed_labeling} and Theorem \ref{theorem:pointed_lyndon_basis_prelie}, but this time in terms of bicolored Lyndon trees.
     Let  $\mathcal{T}Lyn_{n,p}^{w}$ be the set of bicolored Lyndon trees such that along the path from the leaf labeled $p$ to the root, if the path moves to the left, the internal vertex label is $0$ and if it is to the right, the label is given by $1$.

     \begin{remark}
         Note that the definition of $\mathcal{T}Lyn_{n,p}^{w}$ amounts to selecting the bicolored Lyndon trees whose associated maximal chains belong to the interval $[\hat{0},[n]^p]$.
     \end{remark}

\begin{theorem}\label{theorem:consequence_second_pointed_labeling} 
  For every $p\in [n]$ we have that
 \begin{enumerate}
   \item   The order complex $\Delta((\hat{0},[n]^p))$ has the homotopy type of a wedge of $|\mathcal{T}Lyn_{n,p}^{w}|$  many spheres of dimension $n-3$. Hence the interval $[\hat{0},[n]^p]$ is Cohen-Macaulay.

\item The set
$ \{c(T) \mid T\in \mathcal{T}Lyn_{n,p}^{w}\}$ forms a basis for $\widetilde H^{n-3}(\hat{0},[n]^p)$, where $c(T)$ gives the corresponding maximal chain associated to $T$ in $\Pi_n^{\bullet}$.

\item The set
$\{\Theta(T) \mid T\in \mathcal{T}Lyn_{n}^{w}\}$ forms a basis for $\mathcal{P}re\mathcal{L}ie(n)$.
 \end{enumerate}
\end{theorem}

\begin{proof}
Theorem \ref{theorem:EL-labeling_pointed} says that  $\lambda_{\bullet_2}$ is an EL-labeling of $\Pi_n^{\bullet}$. Theorem \ref{Theorem:label-preserving_bijection} implies that the ascent-free words of labels according to $\lambda_{\bullet_2}$ and $\lambda_{w}$ are the same. The ascent-free words of labels of $\lambda_{w}$ are indexed by bicolored Lyndon trees by \cite[Theorem 5.7]{DleonWachs2016}. The reader can check that the bicolored Lyndon trees in $\mathcal{T}Lyn_{n,p}^{w}$ are precisely the ones who index the ascent-free chains in $[\hat{0},[n]^p]$ according to $\lambda_{\bullet_2}$.
\end{proof}

Vallette also concludes in \cite[Theorem 9]{Vallette2007} that a criterion to show that a basic-set quadratic operad $\P$ and its Koszul dual $\P^{\antiexcl}$ have the property of being Koszul is to show that all maximal intervals of its associated operadic partition poset $\Pi^{\P}$ are Cohen-Macaulay. Theorems \ref{theorem:consequence_first_pointed_labeling} and \ref{theorem:consequence_second_pointed_labeling} give then new proofs of the following theorem.

\begin{theorem}[Theorem 1.13 \cite{ChapotonVallete2006}]
    The operads $\Perm$ and $\Prelie$ are Koszul operads.   
\end{theorem}

   \subsection{CL-labelings compatible with isomorphisms and PBW bases}
   In~\cite{BelliermillesDelcroixogerHoffbeck2021} the authors introduce a new compatibility condition on CL-labelings of operadic posets which gives rise to a Poincar\'e–Birkhoff–Witt
  (PBW) basis for the corresponding operad. This PBW basis comes from the {\bf increasing chains} as opposed to the ascent-free chains that is used to give a basis for the cohomology of the poset, and hence for the Koszul dual of the operad. 

  The defined property in \cite{BelliermillesDelcroixogerHoffbeck2021} which a CL-labeling can have is called being \emph{compatible with isomorphism of subposets}. We refer the reader to such article for the complete context and proper definitions which we mostly omit here. Informally, this property requires for operadic posets (in our case of binary generators) that in intervals of the collection $\{\Pi_n^{\P}\}_{n\ge 1}$ that are isomorphic due to the action of the operad $\P$, but perhaps on different sets of inputs, there is also a consistent map between the words of labels of saturated chains on those ``$\P$-isomorphic'' intervals. The requirement is that increasing chains map to increasing chains, ascent-free chains map to ascent-free chains and the lexicographic partial order on chains is preserved on $\P$-isomorphic intervals.

 Both of the labelings $(\Pi_n^w,\lambda_w)$ and $(\Pi_n^{\bullet},\lambda_{\bullet_2})$ depend only on the minimal elements of the blocks that are being merged at each step and the generator of the corresponding operad that is being used to merge the blocks. Because the $\min$ function is preserved under the unique order isomorphism between two totally ordered sets of the same cardinality, we follow a very similar argument as the one in \cite[Proposition 3.11]{BelliermillesDelcroixogerHoffbeck2021} to conclude the following theorem.
 \begin{theorem}\label{theorem:pointed_weighted_isomorphism_subposets}
     The EL-labelings  $(\Pi_n^w,\lambda_w)$ and $(\Pi_n^{\bullet},\lambda_{\bullet_2})$ are compatible with isomorphisms of subposets.
 \end{theorem}

 The following two theorems then highlight the relevance of the notion of CL-labelings compatible with isomorphisms in the context of operad theory.

 \begin{theorem}[Theorem 3.9 \cite{BelliermillesDelcroixogerHoffbeck2021}]\label{theorem:BDH_CL_PBW} A quadratic basic-set operad $\P$ whose operadic poset $\Pi^{\P}_n$ admits a CL-labeling compatible with isomorphisms of subposets admits a partially ordered PBW basis given by the increasing maximal chains of the CL-labeling where the order is given by the lexicographic order on saturated chains.
 \end{theorem}

  \begin{theorem}[Theorem 1.6 \cite{BelliermillesDelcroixogerHoffbeck2021}]\label{theorem:BDH_PBW_Koszul} An operad equipped with a partially ordered PBW basis is Koszul.
 \end{theorem}

 We obtain as a corollary of Theorems \ref{theorem:pointed_weighted_isomorphism_subposets}, \ref{theorem:BDH_CL_PBW}, and \ref{theorem:BDH_PBW_Koszul} a new proof of the fact that the operads $\Com^2$, $\Perm$, and their Koszul duals $\Lie^2$ and $\Prelie$ are all Koszul operads.

\begin{figure}
     \centering
   \begin{tikzpicture}[every node/.style={draw=black},every 
node/.append style={transform shape},scale=0.8]
     \node at ($(one)+(-.75,-.75)$)[circle, inner sep=0](two){\tiny $c_{n-1}$};
     \node at ($(two)+(-.75,-.75)$)[circle, inner sep=0](three){\tiny $c_{n-2}$};
     \node at ($(two)+(.75,-.75)$)[circle,draw=none](2son){$n$};
     \node at ($(three)+(-1.5,-1.5)$)[circle, inner sep=4](four){\tiny $c_{1}$};
      \node at ($(three)+(.75,-.75)$)[circle,draw=none](3son){$n-1$};
     \node at ($(four)+(-0.75,-.75)$)[draw=none](4son1){$1$};
     \node at ($(four)+(0.75,-.75)$)[draw=none](4son2){$2$};
     \path[thick](two)edge(three)edge(2son);
     \path[thick,dashed](three)edge(four);
     \path[thick](three)edge(3son);
     \path[thick](four)edge(4son1)edge(4son2);
\end{tikzpicture}
     \caption{A generic left comb with bicolored internal nodes.}
     \label{fig:left_comb}
 \end{figure}
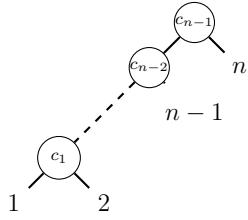

To determine the corresponding PBW bases predicted by Theorem \ref{theorem:BDH_CL_PBW} we use the increasing chains both of the EL-labelings $\lambda_{w}$  (described in \cite[Theorem 3.2]{DleonWachs2016}) and of  $\lambda_{\bullet_2}$  (described in the proof of Theorem \ref{theorem:EL-labeling_pointed}). Note that from Theorem \ref{Theorem:label-preserving_bijection} it follows that the increasing chains in both $(\Pi_n^w,\lambda_w)$ and $(\Pi_n^{\bullet},\lambda_{\bullet_2})$ have the same words of labels. These increasing words of labels are indexed by the following family of trees. Let $\mathrm{lcomb}_{n}^{w}$ be the set of left-combs of the form
 $$((1\wedge^{c_1} 2)\wedge^{c_2} 3)\cdots ) \wedge^{c_{n-1}} n,$$
 where for some $i\in [n]$ we have that $c_1=\cdots =c_{i-1}=0$  and $c_i=\cdots =c_{n-1}=1$ (See Figure \ref{fig:left_comb}).

 \begin{theorem} We have that
 
 \begin{enumerate}
     \item The EL-labeling $(\Pi_n^{\bullet},\lambda_{\bullet_2})$ determines the PBW basis for $\Perm$ formed by the identity and tree-monomials of the form
  $\{\Theta(T) \mid T\in \mathrm{lcomb}_{n}^{w}\}_{n\ge 1}$.
  \item The EL-labeling $(\Pi_n^{w},\lambda_{w})$ determines a PBW basis for $\Com^2$ formed by the identity and tree-monomials of the form
  $$((1\circ_{c_1} 2)\circ_{c_2} 3)\cdots ) \circ_{c_{n-1}} n,$$
   where for some $i\in [n]$ we have that $c_1=\cdots =c_{i-1}=0$  and $c_i=\cdots =c_{n-1}=1$.
 \end{enumerate}
 \end{theorem}

\section{Whitney twins and non-uniqueness of Whitney duals}\label{sec:twins}

\subsection{Whitney twins}
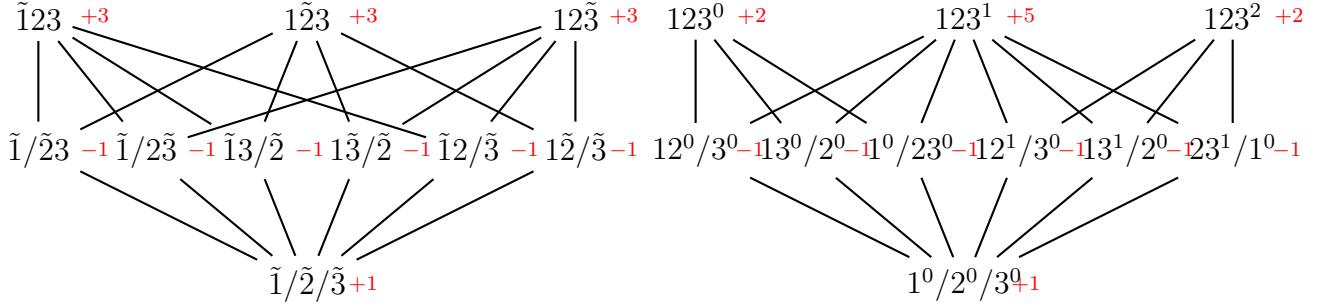
\begin{figure}
    \begin{center}
  \begin{tikzpicture}[scale=.68]
  \begin{scope}[scale=1.05]
\node(root){$\tilde{1}/\tilde{2}/\tilde{3} $};
\node at ($(root)+(-5,5)$)(one){$\tilde{1}23$};
\node at ($(root)+(0,5)$)(two){$1\tilde{2}3$};
\node at ($(root)+(5,5)$)(three){$12\tilde{3}$};
\node at ($(root)+(-5,2.5)$)(21){$\tilde{1}/\tilde{2}3$};
\node at ($(root)+(-3,2.5)$)(22){$\tilde{1}/2\tilde{3}$};
\node at ($(root)+(-1,2.5)$)(23){$\tilde{1}3/\tilde{2}$};
\node at ($(root)+(1,2.5)$)(24){$1\tilde{3}/\tilde{2}$};
\node at ($(root)+(3,2.5)$)(25){$\tilde{1}2/\tilde{3}$};
\node at ($(root)+(5,2.5)$)(26){$1\tilde{2}/\tilde{3}$};

\draw[thick] (root)--(21);
\draw[thick] (root)--(22);
\draw[thick] (root)--(23);
\draw[thick] (root)--(24);
\draw[thick] (root)--(25);
\draw[thick] (root)--(26);

\draw[thick] (21)--(one);
\draw[thick] (22)--(one);
\draw[thick] (23)--(one);
\draw[thick] (25)--(one);

\draw[thick] (21)--(two);
\draw[thick] (23)--(two);
\draw[thick] (24)--(two);
\draw[thick] (26)--(two);

\draw[thick] (22)--(three);
\draw[thick] (24)--(three);
\draw[thick] (25)--(three);
\draw[thick] (26)--(three);

\node at ($(root)+(1.05,0)$){$\textcolor{red}{\scriptstyle +1}$};
\node at ($(root)+(-5,5)+(1.05,0)$) {$\textcolor{red}{\scriptstyle +3}$};
\node at ($(root)+(0,5)+(1.05,0)$){$\textcolor{red}{\scriptstyle +3}$};
\node at ($(root)+(5,5)+(.9,0)$){$\textcolor{red}{\scriptstyle  +3}$};
\node at ($(root)+(-5,2.5)+(1.05,0)$){$\textcolor{red}{\scriptstyle -1}$};
\node at ($(root)+(-3,2.5)+(1.05,0)$){$\textcolor{red}{\scriptstyle -1}$};
\node at ($(root)+(-1,2.5)+(1.05,0)$){$\textcolor{red}{\scriptstyle  -1}$};
\node at ($(root)+(1,2.5)+(1.05,0)$){$\textcolor{red}{\scriptstyle -1}$};
\node at ($(root)+(3,2.5)+(1.05,0)$){$\textcolor{red}{\scriptstyle -1}$};
\node at ($(root)+(5,2.5)+(.9,0)$){$\textcolor{red}{\scriptstyle -1}$};
  \end{scope}

\begin{scope}[shift = {(12.85,0)},scale=1.05]

\node (n1232) at (5,5) {$123^{2}$};

 \node (n13020) at (-3,2.5) {$13^{ 0}/ 2^{ 0}$};
  \node (n102030) at (0,0)  {$1^{0}/ 2^{0}/ 3^{0}$};
  \node (n1231) at (0,5) {$123^{1}$};
  \node (n12030) at (-5,2.5) {$12^{0}/ 3^{0}$};
  \node (n13120) at (3,2.5)  {$13^{1}/ 2^{0}$};
  \node (n1230) at (-5,5){$123^ {0}$};
  \node (n10230) at (-1,2.5)  {$1^{0}/ 23^{0}$};
  \node (n12130) at (1,2.5)  {$12^{ 1}/ 3^{0}$};
  \node (n10231) at (5,2.5) {$23^ {1}/ 1^{0}$};

  \draw[thick] (n1231) -- (n10230) ;
  \draw [thick] (n13020) -- (n102030);
  \draw [thick] (n1232) -- (n13120);
  \draw [thick] (n1231)-- (n13020);
  \draw [thick] (n10230)--(n102030);
  \draw [thick] (n1230) -- (n10230);
  \draw [thick] (n1231) -- (n13120);
  \draw [thick] (n12030)-- (n102030);
  \draw [thick] (n1231) --(n12130);
  \draw [thick] (n1232) -- (n12130);
  \draw [thick] (n13120) --(n102030);
  \draw [thick] (n1231) --(n10231);
  \draw [thick] (n1230) -- (n13020);
  \draw [thick] (n1230)  -- (n12030);
  \draw [thick] (n12130) --  (n102030);
  \draw [thick] (n1232)  --  (n10231);
  \draw [thick] (n10231)  --  (n102030);
  \draw [thick] (n1231) -- (n12030);


\node  at ($(5,5)+(1.05,0)$){$\textcolor{red}{\scriptstyle +2}$};
\node  at ($(-5,5)+(1.05,0)$){$\textcolor{red}{\scriptstyle +2}$};
\node  at ($(0,5) +(1.05,0)$){$\textcolor{red}{\scriptstyle +5}$};
\node  at ($ (-3,2.5) +(1.01,0)$){$\textcolor{red}{\scriptstyle -1}$};
\node  at ($(-5,2.5) +(1.01,0)$){$\textcolor{red}{\scriptstyle -1}$};
\node  at ($(3,2.5)  +(1.01,0)$){$\textcolor{red}{\scriptstyle -1}$};
\node  at ($(-1,2.5)  +(1.01,0)$){$\textcolor{red}{\scriptstyle -1}$};
\node  at ($(1,2.5)   +(1.01,0)$){$\textcolor{red}{\scriptstyle -1}$};
\node  at ($ (5,2.5)+(1.01,0)$){$\textcolor{red}{\scriptstyle -1}$};
\node  at ($ (0,0) +(1.15,0)$){$\textcolor{red}{\scriptstyle +1}$};

\end{scope}

 \end{tikzpicture}
 \end{center}
 \caption{$\Pi_{3}^{\bullet}$ and $\Pi_{3}^{w}$. M\"obius values in red. } 
    \label{fig:PointedAndWeighted}
\end{figure}

The reader might have noticed at this point that the pointed and the weighted partition posets are closely related.  From Figure~\ref{fig:PointedAndWeighted}, we can see that already the posets $\Pi_3^{\bullet}$ and $\Pi_3^{w}$ are not isomorphic. This can be easily shown in general for $\Pi_n^{\bullet}$ and $\Pi_n^{w}$ since all maximal intervals in $\Pi_n^{\bullet}$  are isomorphic but this is not the case in $\Pi_n^{w}$. In particular, for the latter poset the intervals $[\hat{0},[n]^0]$ and $[\hat{0},[n]^{n-1}]$ are isomorphic to $\Pi_n$ which is not the case for any maximal interval in $\Pi_n^{\bullet}$ for $n\ge 3$.

The reader can verify from Figure~\ref{fig:PointedAndWeighted} that the Whitney numbers of the first and second kind are the same for $\Pi_3^{\bullet}$ and $ \Pi_{3}^{w}$. In \cite{DleonWachs2016} the authors prove that this is true for any $n\ge 1$. Indeed there is a rank preserving bijection $\Pi_n^{w}\rightarrow \Pi_{n}^{\bullet}$ induced by transforming, in a weighted partition, every weighted set $A^w$ into the pointed set $A^{p_{w}}$ where  $A=\{p_0< p_1<\dots <p_{|A|-1}\}$. The authors then use the fact that the two posets are uniform to conclude that their Whitney numbers of the first and second kind are the same. This gives an example of the next definition.

\begin{definition}\label{definition:whitney_twins}
Two graded posets $P$ and $Q$  are said to be \emph{Whitney twins} if their Whitney numbers of the first and second kind are the same, i.e., they satisfy
$$w_k(P)=w_k(Q)\text{ and } W_k(P)=W_k(Q)$$
for all $k$.
\end{definition}

Thus in our new terminology, the results in \cite{DleonWachs2016} can be recast into the following proposition.

\begin{proposition} [{\cite[Section 2.4]{DleonWachs2016}}] \label{proposition:whitney_twins_weighted_pointed} For all $n\ge 1$, the posets
$ \Pi_n^{\bullet}$ and $ \Pi_{n}^{w}$ are Whitney twins. 
 \end{proposition}
 
Note that if $P_1$ and $P_2$ are Whitney twins and $Q_1$ and $Q_2$ are Whitney duals of $P_1$ and $P_2$ respectively, then $Q_1$ and $Q_2$ are Whitney twins.  Thus we also have the following immediate corollary from Proposition \ref{proposition:whitney_twins_weighted_pointed} and Theorems \ref{theorem:flynpointed_dual} and \ref{theorem:flynweighted_dual}.

\begin{corollary}\label{cor:pointedLynAndBicoloredAreWhitneyTwins}
For all $n\geq 1$, the posets $\Flyn_{n}^{w}$ and $\Flyn_{n}^{\bullet}$ are Whitney twins.
\end{corollary}
We should note that if $P$ and $Q$ are isomorphic, they are  Whitney twins.  Thus, at this point, it could be that   $\Flyn_{n}^{w}$ and $\Flyn_{n}^{\bullet}$ are Whitney twins merely because they are isomorphic. We will show in Theorem~\ref{thm:nonIsomorphicWDual} that this is only true for $n\le 3$ and is not the case for $n\ge 4$.

\subsection{Non-uniqueness of Whitney duals}

   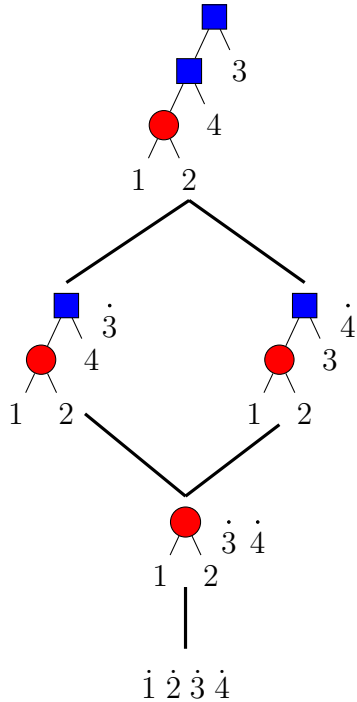
\begin{figure}[htbp]
    \hspace*{-2cm}
    \begin{center}
    \begin{tikzpicture}[scale=0.6,
rdot/.style = {circle, draw, fill=red,
               minimum size=4mm, inner sep=0pt, outer sep=0pt,
               node contents={}},
sdot/.style  ={draw, fill=blue, scale=1.2,
               minimum size=2.7mm, inner sep=0pt, outer sep=0pt,
               node contents={}},
sibling distance = 7mm,scale=0.8
                    ]
\begin{scope}[shift={(142mm,-140mm)}]
\node (three1) [rdot] 
    child{ node {$1$}} 
    child{ node (three1 south){} edge from parent[draw=none]}
    child{ node {2}};
    \filldraw[black] (1.2,0) circle (1pt) node[anchor=west][draw=none,below] {$3$};
    \filldraw[black] (2,0) circle (1pt) node[anchor=west][draw=none,below] {$4$};
\end{scope}
\begin{scope}[xshift=150mm]
\node (top) [sdot] 
    child{ node [scale=1.2,fill=blue,draw]{}
    child{node[circle,draw,fill=red]{}
    child{node{1}}
    child{node{} edge from parent[draw=none]}
    child{node(four s){2}}}
    child{node{} edge from parent[draw=none]}
    child{node{4}}}
    child{node{} edge from parent[draw=none]}
    child{node{3}};
\end{scope}

\begin{scope}[shift={(109mm,-80mm)}]
\node (four) [sdot] 
    child{ node [circle,fill=red,draw]{}
    child{node{1}}
    child{node{} edge from parent[draw=none]}
    child{node(four south){2}}}
    child{node{} edge from parent[draw=none]}
    child{node{$4$}};
    \filldraw[black] (1.2,0) circle (1pt) node[anchor=west][draw=none,below] {$3$};
\end{scope}
\begin{scope}[shift={(175mm,-80mm)}]
\node (five) [sdot] 
    child{ node [circle,draw,fill=red]{}
    child{node{$1$}}
    child{node(five south){} edge from parent[draw=none]}
    child{node{2}}}
    child{node{} edge from parent[draw=none]}
    child{node{3}};
    \filldraw[black] (1.2,0) circle (1pt) node[anchor=west][draw=none,below] {$4$};
\end{scope}

\begin{scope}[shift={(142mm,-185mm)}]
\node(et) {$\dot{1}\; \dot{2}\; \dot{3}\; \dot{4}$};
\end{scope}

\draw[very thick]
([yshift=3mm] three1.north) -- (four south.east)
([yshift=3mm] three1.north) -- (five south.south)
([yshift=3mm] four.north) -- (four s.south)
([yshift=3mm] five.north) -- (four s.south)
([yshift=3mm] et.north) -- (three1 south.south);
    \end{tikzpicture}
\end{center}
\caption{An interval of rank $3$ in $\Flyn_{4}^{w}$.}
\label{fig:FLyn_4-interval}

    \end{figure}

As mentioned in Corollary~\ref{cor:pointedLynAndBicoloredAreWhitneyTwins}, $\Flyn_{n}^{w}$ and  $\Flyn_{n}^{\bullet}$ are Whitney twins.  Here we explain why they are not isomorphic in general.  This in turn will show that a poset can have multiple (non-isomorphic) Whitney duals. We also argue that another poset $\SF_n$ already studied by Reiner \cite{Reiner1978} and Sagan \cite{Sagan1983} is a third non-isomorphic Whitney dual of $\Pi_n^w$ and $\Pi_n^{\bullet}$.

\begin{theorem}\label{thm:nonIsomorphicWDual}
  For $n\geq4$, $\Flyn_{n}^{w}$ and  $\Flyn_{n}^{\bullet}$ are not isomorphic. Consequently, $\Pi_n^w$ and $\Pi_n^\bullet$ have multiple Whitney duals.
\end{theorem}
   \begin{proof}
   Consider the maximal interval of $\Flyn_{4}^{w}$ depicted in Figure~\ref{fig:FLyn_4-interval}. This interval occurs in $\Flyn_{n}^{w}$ for all $n\geq 4$ since adding isolated vertices to the forests of the interval does not change the interval's structure.  We claim that  there are no   intervals in $\Flyn_{n}^{\bullet}$ (for $n\geq 4$) that start at $\hat{0}$ and are isomorphic to the interval in Figure~\ref{fig:FLyn_4-interval}.  Note that if we can verify this claim we will be done.  
   
   Suppose that such an  interval in $\Flyn_{n}^{\bullet}$  exists and let $I$ be this interval. Note that the cover relation on $\Flyn_{n}^{\bullet}$ only depends on the relative order of the leaf labels and not the actual leaf labels themselves. So $I$ must be isomorphic to an interval starting at $\hat{0}$ in $\Flyn_{4}^{\bullet}$. A simple check (see \cite{QuicenoDuran2020} for a complete argument) of the intervals of $\Flyn_{4}^{\bullet}$ shows that no intervals starting at $\hat{0}$ are isomorphic to $I$, completing the proof.
  \end{proof}
  
 Reiner~\cite{Reiner1978}  introduced a family of posets of rooted spanning forests $\SF_n$ and  Sagan~\cite{Sagan1983} computed the Whitney numbers of these posets. The poset  $\SF_n$  is formed by rooted spanning forests where cover relations happen when two rooted trees are merged by their roots selecting the new roots from the two that have been merged (see Figure \ref{fig:SF3} for an example, there the square (red) nodes represent the roots of the trees). 
 As mentioned in~\cite{DleonWachs2016}, the Whitney numbers $\Pi_n^w$ and $\Pi_n^{\bullet}$ are switched as compared to $\SF_n$, which implies that $\SF_n$ is also a Whitney dual to both posets.  From Figures \ref{fig:FLyn_3-pointed}, \ref{fig:FLyn_3_weighted}, and \ref{fig:SF3} it is already evident that $\SF_3$ is not isomorphic to $\Flyn_{3}^{w}\cong \Flyn_{3}^{\bullet}$. We show here that in fact $\SF_n$ is not isomorphic to $\Flyn_{n}^{w}$ or $\Flyn_{n}^{\bullet}$ for $n\ge 3$.

\begin{figure}
    \centering
\begin{tikzpicture}[scale=0.35]

\tikzstyle{every node}=[fill,gray!10,circle, inner sep=0pt,minimum size=80,scale=0.45]
\node (p1) at (13,0.25) {};
\node (q1) at (-4.5,5.25) {};
\node (q2) at (0.5,5.25) {};
\node (q3) at (10.5,5.25) {};
\node (q4) at (15.5,5.25) {};
\node (q5) at (25.5,5.25) {};
\node (q6) at (30.5,5.25) {};

\node (r1) at (-7,10.25) {};
\node (r2) at (-2,10.25) {};
\node (r3) at (3,10.25) {};
\node (r4) at (8,10.25) {};
\node (r5) at (13,10.25) {};
\node (r6) at (18,10.25) {};
\node (r7) at (23,10.25) {};
\node (r8) at (28,10.25) {};
\node (r9) at (33,10.25) {};

\tikzstyle{every node}=[black, draw, circle, inner sep=1, minimum size=15,scale=0.45]
\tikzstyle{every path}=[blue, line width=0.2pt]

\node[rectangle,color=red] (a1) at (13,1){$1$};
\node[rectangle,color=red]  (a2) at (14,0){$2$};
\node[rectangle,color=red]  (a3) at (12,0){$3$};

\node[rectangle,color=red]  (b11) at (0.5,6){$1$};
\node[rectangle,color=red]  (b12) at (1.5,5){$2$};
\node (b13) at (-0.5,5){$3$};
 \draw[thick] (b11) -- (b13) ;

\node[rectangle,color=red]  (c11) at (-4.5,6){$1$};
\node (c12) at (-3.5,5){$2$};
\node[rectangle,color=red]  (c13) at (-5.5,5){$3$};
 \draw[thick] (c11) -- (c12) ;

\node[rectangle,color=red]  (d11) at (10.5,6){$1$};
\node[rectangle,color=red]  (d12) at (11.5,5){$2$};
\node (d13) at (9.5,5){$3$};
 \draw[thick] (d12) -- (d13) ;

\node (b21) at (25.5,6){$1$};
\node[rectangle,color=red]  (b22) at (26.5,5){$2$};
\node[rectangle,color=red]  (b23) at (24.5,5){$3$};
 \draw[thick] (b21) -- (b23) ;

\node (c21) at (15.5,6){$1$};
\node[rectangle,color=red]  (c22) at (16.5,5){$2$};
\node[rectangle,color=red]  (c23) at (14.5,5){$3$};
 \draw[thick] (c21) -- (c22) ;

\node[rectangle,color=red]  (d21) at (30.5,6){$1$};
\node (d22) at (31.5,5){$2$};
\node[rectangle,color=red]  (d23) at (29.5,5){$3$};
 \draw[thick] (d22) -- (d23) ;

\node (e11) at (-7,11){$1$};
\node (e12) at (-6,10){$2$};
\node[rectangle,color=red]  (e13) at (-8,10){$3$};
\draw[thick] (e12)--(e11)--(e13);

\node[rectangle,color=red]  (e21) at (-2,11){$1$};
\node (e22) at (-1,10){$2$};
\node (e23) at (-3,10){$3$};
\draw[thick] (e22)--(e21)--(e23);

\node (e31) at (3,11){$1$};
\node[rectangle,color=red]  (e32) at (4,10){$2$};
\node (e33) at (2,10){$3$};
\draw[thick] (e32)--(e31)--(e33);

\node[rectangle,color=red]  (f11) at (8,11){$1$};
\node (f12) at (9,10){$2$};
\node (f13) at (7,10){$3$};
\draw[thick] (f11)--(f12)--(f13);

\node (f21) at (13,11){$1$};
\node[rectangle,color=red]  (f22) at (14,10){$2$};
\node (f23) at (12,10){$3$};
\draw[thick] (f21)--(f22)--(f23);

\node (f31) at (18,11){$1$};
\node (f32) at (19,10){$2$};
\node[rectangle,color=red]  (f33) at (17,10){$3$};
\draw[thick] (f31)--(f32)--(f33);

\node (g11) at (23,11){$1$};
\node[rectangle,color=red]  (g12) at (24,10){$2$};
\node (g13) at (22,10){$3$};
\draw[thick] (g11)--(g13)--(g12);

\node (g21) at (28,11){$1$};
\node (g22) at (29,10){$2$};
\node[rectangle,color=red]  (g23) at (27,10){$3$};
\draw[thick] (g21)--(g23)--(g22);

\node[rectangle,color=red]  (g31) at (33,11){$1$};
\node (g32) at (34,10){$2$};
\node (g33) at (32,10){$3$};
\draw[thick] (g31)--(g33)--(g32);

\tikzstyle{every path}=[line width=1.5pt]
 \draw (p1) -- (q1) ;
 \draw (p1) -- (q2) ;
 \draw (p1) -- (q3) ;
 \draw (p1) -- (q4) ;
 \draw (p1) -- (q5) ;
 \draw (p1) -- (q6) ;

 \draw (q1) -- (r1) ;
 \draw (q1) -- (r2) ;
 \draw (q2) -- (r2) ;
 \draw (q2) -- (r3) ;
 \draw (q3) -- (r4) ;
 \draw (q3) -- (r5) ;
 \draw (q4) -- (r5) ;
 \draw (q4) -- (r6) ;
 \draw (q5) -- (r7) ;
 \draw (q5) -- (r8) ;
 \draw (q6) -- (r8) ;
 \draw (q6) -- (r9) ;

\end{tikzpicture}

    \caption{$\SF_3$}
    \label{fig:SF3}
\end{figure}

\begin{theorem}\label{theorem:flyn_SF_non_isomorphic}
 For $n\geq 3$, $\Flyn_{n}^{\bullet}$ and $\SF_n$ are not isomoprhic. 
\end{theorem}

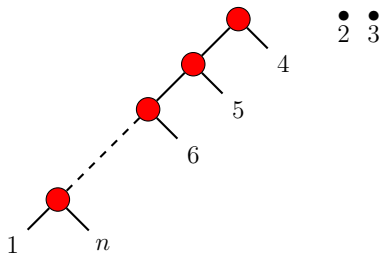
\begin{figure}
     \centering
   \begin{tikzpicture}[every node/.style={draw=black},every 
node/.append style={transform shape},scale=0.8]
     \node[circle,fill=red](one){};
     \node at ($(one)+(-.75,-.75)$)[circle,fill=red](two){};
     \node at ($(two)+(-.75,-.75)$)[circle,fill=red](three){};
     \node at ($(two)+(.75,-.75)$)[circle,draw=none](2son){5};
     \node at ($(three)+(-1.5,-1.5)$)[circle,fill=red](four){};
      \node at ($(three)+(.75,-.75)$)[circle,draw=none](3son){6};
     \node at ($(four)+(-0.75,-.75)$)[draw=none](4son1){1};
     \node at ($(four)+(0.75,-.75)$)[draw=none](4son2){$n$};
     \node at ($(two)+(.75,-.75)$)[circle,draw=none](five){};
     \node at ($(one)+(.75,-.75)$)[circle,draw=none](six){4};
      \node at ($(two)+(2.5,0.8)$)[circle,draw=none]{$\bullet$};
     \node at ($(two)+(3,0.8)$)[circle,draw=none]{$\bullet$};
     \node at ($(two)+(2.5,0.5)$)[circle,draw=none]{2};
     \node at ($(two)+(3,0.5)$)[circle,draw=none]{3};
     \path[thick](one)edge(two)edge(six);
     \path[thick](two)edge(three)edge(2son);
     \path[thick,dashed](three)edge(four);
     \path[thick](three)edge(3son);
     \path[thick](four)edge(4son1)edge(4son2);
\end{tikzpicture}
     \caption{Pointed Lyndon forest used in the proof of Theorem \ref{theorem:flyn_SF_non_isomorphic}.}
     \label{fig:pointed_forest_counterexample}
 \end{figure}

\begin{proof}
Note first that $\SF_n$ is an \emph{uniform} graded poset according to the definition in \cite{Gonzalez2016}. More specifically, if $F\in \SF_n$ is an element of rank $\rho(F)=i$ then the filter $U(F)$ in $\SF_n$ is isomorphic to $\SF_{n-i}$. Indeed, the rules of merging in the filter $U(F)$ are only dependent on the roots of $F$ and any $F\in \SF_n$ of rank $\rho(F)=i$ has $n-i$ roots.

When $n=3$, the posets $\SF_3$ and $\Flyn^\bullet_{3}$ are clearly non-isomorphic as can be appreciated from Figures \ref{fig:FLyn_3-pointed} and \ref{fig:SF3},  so let us assume that $n\ge 4$.
Consider the pointed Lyndon forest $F$ of Figure \ref{fig:pointed_forest_counterexample}. Since the root of the nontrivial tree in $F$ is a Lyndon node with minimal element of the right subtree $4$ that is larger than $2$ and $3$, that the filter $U(F)$ in $\Flyn_{n}^{\bullet}$ is isomorphic to $\Flyn^\bullet_{3}$.  Now, if there is an isomorphism $f:\Flyn_{n}^{\bullet}\rightarrow \SF_n$, this induces an isomorphism $U(F)\cong U(f(F))\cong \SF_3$ since the element $f(F)$ has rank $n-3$, but this is a contradiction.
\end{proof}

The proof of the following theorem follows the same  idea as in Theorem \ref{theorem:flyn_SF_non_isomorphic}
\begin{theorem}
 For $n\geq 3$, $\Flyn_{n}^{w}$ and $\SF_n$ are not isomorphic. 
\end{theorem}

 \begin{remark}
 We should note that the first two authors have found a CW-labeling (a more general version of an EW-labeling) of $\Pi_n^w$ whose corresponding Whitney dual is $\SF_n$. This will be further discussed in a coming work and can be already found in the ArXiv version of~\cite{GonzalezHallam2021}. 
 \end{remark}

 \section*{Acknowledgments}
  We would like to thank Joan Bellier-Mill{\`e}s,  B{\'e}r{\'e}nice Delcroix-Oger, and Eric Hoffbeck, for the helpful conversation and explanation of the context of CL-labelings compatible with isomorphisms coming from their work in  \cite{BelliermillesDelcroixogerHoffbeck2021}.
  
\bibliographystyle{plain}
\bibliography{whitneydualref}
\end{document}